\documentclass[preprint,3p,12pt]{elsarticle}
\usepackage{mathrsfs}
\usepackage{amsmath}
\usepackage{stmaryrd}
\usepackage{bbding}
\usepackage{dcolumn}
\usepackage{graphicx}
\usepackage{amsfonts}
\usepackage{amssymb}
\usepackage{psfrag}
\usepackage{wrapfig}
\usepackage{subfigure}
\usepackage{makeidx}
\usepackage{bm}
\usepackage{epsf}
\usepackage{color}
\usepackage{epsfig}
\usepackage{setspace}
\usepackage{graphicx}
\usepackage{epstopdf}
\usepackage{psfrag}
\usepackage{subfigure}
\usepackage{multirow}
\usepackage{diagbox}
\usepackage{verbatim}
\newtheorem{rmk}{Remark}

\usepackage{makecell}
\usepackage{float}
\usepackage{esint}
\usepackage{comment}
\usepackage{movie15}
\usepackage{hyperref}

\newcommand{\mathsym}[1]{{}}
\newcommand{\unicode}[1]{{}}
\begin{document}
	
\title{A gradient-compression-based compact high-order gas-kinetic scheme on three-dimensional hybrid unstructured mesh}
	
	\author[HKUST1]{Xing Ji}
	\ead{xjiad@connect.ust.hk}
		
	\author[HKUST2]{Wei Shyy}
	\ead{weishyy@ust.hk}
	
	\author[HKUST1,HKUST2,HKUST3]{Kun Xu\corref{cor1}}
	\ead{makxu@ust.hk}

	\address[HKUST1]{Department of Mathematics, Hong Kong University of Science and Technology, Clear Water Bay, Kowloon, Hong Kong}
	\address[HKUST2]{Department of Mechanical and Aerospace Engineering, Hong Kong University of Science and Technology, Clear Water Bay, Kowloon, Hong Kong}
	\address[HKUST3]{Shenzhen Research Institute, Hong Kong University of Science and Technology, Shenzhen, China}
	\cortext[cor1]{Corresponding author}

\begin{abstract}
In this paper, the compact gas-kinetic scheme for compressible flow is extended to hybrid unstructured mesh.
Based on both cell-averaged flow variables and their gradients updated from time accurate gas evolution model at cell interfaces,
a compact third-order least-square-constrained reconstruction can be obtained on unstructured mesh and
a multi-resolution WENO reconstruction is adopted in case of discontinuous solutions.
Moreover, a compression factor for the cell-averaged gradients is proposed to take into account the possible discontinuity in flow variables at cell interface, which significantly improves the robustness of the compact scheme for high-speed flow computation on irregular mesh
and preserves the accuracy.
Numerical tests from incompressible to hypersonic flow are presented to demonstrate the broad applicability of the gradient-compression-based high-order compact scheme.
\end{abstract}


\maketitle

\section{Introduction}
The simulation of compressible flow with complex geometry is of vital importance in the engineering applications of aerospace industry.
The use of unstructured mesh is especially favored because of its geometric flexibility.
The compact gas-kinetic scheme (CGKS) on tetrahedral mesh has been developed recently \cite{ji2021two}.
However, to resolve the viscous boundary layers efficiently, numerical methods based on the hexahedral or prismatic elements with high aspect ratio are needed in aeronautical practice, and present more accurate and stable solutions than those on the tetrahedral elements alone.
Therefore, the development of CGKS on hybrid mesh is necessary for real-world engineering applications.

Computational methods for compressible flow can be generally categorized into compact and non-compact methods in terms of the stencils used.
As a non-compact scheme, the high-order finite volume methods (FVM) with the weighted essentially non-oscillatory (WENO) reconstruction have been developed and applied continuously to a large-scale aeronautical simulation on hybrid unstructured mesh\cite{antoniadis2017assessment}.
The robustness of the schemes can be improved by the extended stencils in reconstruction.
However, they have difficulties in code portability, parallel programming, and boundary treatment.

On the other hand, methods with compact stencils have simple geometry dependency, which brings great mesh adaptability and high scalability.
The development of high-order compact methods becomes a hot topic nowadays.
Two main representatives are the DG \cite{shu2016weno-dg-review} and the FR/CPR methods \cite{Huynh2007FR,yu2014accuracy}.
By updating variables in multiple degrees of freedom (DOFs), these methods can achieve arbitrary spatial order of accuracy with only the targeted cell as the reconstruction stencil.
Successful examples have been demonstrated in large eddy simulation (LES) \cite{wang2017towards} and RANS simulation \cite{yang2019robust} for subsonic flows.
For the flow simulation with discontinuities, these methods usually have less robustness against the traditional high-order FVMs.
In addition, these methods have restricted explicit time steps and high memory consumption \cite{luo2010reconstructed}.
The $P_NP_M$ \cite{dumbser2010arbitrary} and reconstructed-DG (rDG) methods \cite{luo2010reconstructed} target to overcome  these weakness with the release of the compactness of the DG methods.
In comparison with the DG methods, the same order rDG methods can use larger time step and have less memory requirement.

In recent years, a class of high-order CGKS has been developed from the second-order gas-kinetic scheme (GKS) \cite{xu2014direct}.
The CGKS is based on a time-accurate evolution model for the construction of the gas distribution function at a cell interface \cite{pan2016unstructuredcompact}.
The time-dependent solution provides not only the fluxes across a cell interface but also the corresponding flow variables.
As a result, besides updating the cell-averaged flow variables, the cell-averaged gradients of flow variables can be updated simultaneously through the divergence theorem.
The DOFs updated in CGKS are based on the time accurate dynamic evolution solution rather than the weak formulation in the DG/rDG methods.
 Due to their differences, the CGKS can use a larger time step and has better robustness than the same-order DG methods. For example, a CFL number around 0.5 can be taken for the third-order CGKS \cite{ji2021compact}
while it is restricted to be less than 0.33 for the third-order P1P2-rDG scheme.
The P1P2-rDG is claimed to be unstable on tetrahedral mesh with smooth reconstruction.
However, the third-order CGKS is stable with a CFL number of 1 with the same compact stencil \cite{ji2021two}.
Another feature of GKS is to use the two-stage fourth-order temporal discretization method \cite{li2016twostage} or the multi-stage multi-derivative time marching methods \cite{multi-derivative}.
Although the gas-kinetic flux function is more expensive than the time-independent Riemann solvers,
the GKS can achieve fourth-order temporal accuracy with only two stages \cite{Pan2016twostage}, instead of four stages by the fourth-order Runge-Kutta (RK) time discretization.
Overall, the CGKS  turns out to be more efficient in comparison with Riemann-solver-based RK methods \cite{ji2021compact} in serial computation. Moreover, higher parallel efficiency is expected since less communication is required due to the less middle stages.

In this paper, the compact third-order GKS will be extended to mixed-element mesh.
The scheme is linearly stable for smooth flow with unlimited constrained-least-square reconstruction on a compact stencil involving von Neumann neighbors only.
For discontinuous flow, the idea of the multi-resolution WENO reconstruction is adopted. The reconstruction is designed in a hierarchical way,
i.e., the Nth-order of accuracy can be achieved by N central stencils from first-order to Nth-order \cite{zhu2020new}.
In previous work \cite{ji2021two}, the smooth indicators are determined from the cell-averaged conservative flow variables and a two-step reconstruction is proposed to improve the robustness of the scheme.
In the current work, the complexity of the spatial reconstruction is further reduced. Each low-order stencil is simply chosen as one cell in the compact stencil and the smooth indicator is directly obtained from the corresponding cell-averaged slopes.
In this way, there is no extra memory requirement for the sub-stencils and the computational cost for constructing the corresponding low-order polynomials is reduced.
In case of discontinuities in the flow variables at a cell interface, the current gradient update scheme with continuous assumption of flow variables
at cell interface will have error in the updated gradients.
In order to improve the mesh adaptability, the accuracy of updated solution, and the robustness of the scheme for flow simulation with strong shocks,
a cell-averaged gradient  compression factor (CF)  will be proposed to modify the updated slope in the discontinuous flow region.
Based on the CGKS framework, the CF is different from the existing priori or
posteriori limiters \cite{krivodonova2004shock,zhang2017positivity,clain2011high}.
It has the following features:
i) accuracy preserving;
ii) negligible computational cost;
iii) combined with the multi-resolution WENO reconstruction, the scheme can truly reduce to the first-order GKS when numerical discontinuities appear regardless of the local mesh quality.
The resulting CGKS becomes efficient and robust, and easy to program.
Stringent tests including hypersonic flow passing through a space-vehicle validate the robustness of the current compact scheme with complex geometry.

This paper is organized as follows.
The basic framework of the compact high-order GKS on unstructured mesh is presented in Section 2.
In Section 3, the details for the spatial reconstruction on mesh with mixed elements are presented, including the construction of the CF.
Numerical examples from nearly incompressible to hypersonic flows are given in Section 4.
Discussions on the newly designed CF and a concluding remark are given in the last section.

\section{Compact finite volume gas-kinetic scheme}

The 3-D gas-kinetic BGK equation \cite{BGK} is
\begin{equation}\label{bgk}
f_t+\textbf{u}\cdot\nabla f=\frac{g-f}{\tau},
\end{equation}
where $f=f(\textbf{x},t,\textbf{u},\xi)$ is the gas distribution function, which is a function of space $\textbf{x}$, time $t$, particle velocity $\textbf{u}$, and internal variable $\xi$.
 $g$ is the equilibrium state approached by $f$
and $\tau$ is the collision time.

The collision term satisfies the compatibility condition
\begin{equation*}\label{compatibility}
\int \frac{g-f}{\tau} \pmb{\psi} \text{d}\Xi=0,
\end{equation*}
where $\pmb{\psi}=(1,\textbf{u},\displaystyle \frac{1}{2}(\textbf{u}^2+\xi^2))^T$,
$\text{d}\Xi=\text{d}u_1\text{d}u_2\text{d}u_3\text{d}\xi_1...\text{d}\xi_{K}$,
$K$ is the number of internal degrees of freedom, i.e.
$K=(5-3\gamma)/(\gamma-1)$ in 3-D case, and $\gamma$
is the specific heat ratio.

In the continuum flow regime with the smoothness assumption, based on the Chapman-Enskog expansion the gas distribution function can be expressed as \cite{xu2014direct} ,
\begin{align*}
f=g-\tau D_{\textbf{u}}g+\tau D_{\textbf{u}}(\tau
D_{\textbf{u}})g-\tau D_{\textbf{u}}[\tau D_{\textbf{u}}(\tau
D_{\textbf{u}})g]+...,
\end{align*}
where $D_{\textbf{u}}={\partial}/{\partial t}+\textbf{u}\cdot \nabla$.
Different hydrodynamic equations can be derived by truncating on different orders of $\tau$.
With the zeroth-order in truncated distribution function $f=g$,
the Euler equations can be recovered by multiplying $\pmb{\psi}$ on Eq.~\eqref{bgk} and integrating it over the phase space,
\begin{equation*}\label{euler-conservation}
\begin{split}
\textbf{W}_t+ \nabla \cdot \textbf{F}(\textbf{W})=0.
\end{split}
\end{equation*}
With the first-order truncation, i.e.,
\begin{align*} \label{ce-ns}
f=g-\tau (\textbf{u} \cdot \nabla g + g_t),
\end{align*}
the N-S equations can be obtained,
\begin{equation*}\label{ns-conservation}
\begin{split}
\textbf{W}_t+ \nabla \cdot \textbf{F}(\textbf{W},\nabla \textbf{W} )=0,
\end{split}
\end{equation*}
with $\tau = \mu / p$ and $Pr=1$.

The conservative flow variables and their fluxes are the moments of the gas distribution function
\begin{align}\label{point}
\textbf{W}(\textbf{x},t)=\int \pmb{\psi} f(\textbf{x},t,\textbf{u},\xi)\text{d}\Xi,
\end{align}
and
\begin{equation}\label{f-to-flux}
\textbf{F}(\textbf{x},t)=
\int \textbf{u} \pmb{\psi} f(\textbf{x},t,\textbf{u},\xi)\text{d}\Xi.
\end{equation}

\subsection{Compact gas-kinetic scheme on mixed-elements}

For a 3-D polyhedral cell $\Omega_i$, the boundary can be expressed as
\begin{equation*}
\partial \Omega_i=\bigcup_{p=1}^{N_f}\Gamma_{ip},
\end{equation*}
where $N_f$ is the number of cell interfaces for cell $\Omega_i$.
$N_f=4$ for tetrahedron, $N_f=5$ for prism and pyramid, $N_f=6$ for hexahedron.

\begin{figure}[htbp]	
	\centering
	\subfigure[Triangluar face]{
		\label{triangle}
		\includegraphics[width=0.48\textwidth]
		{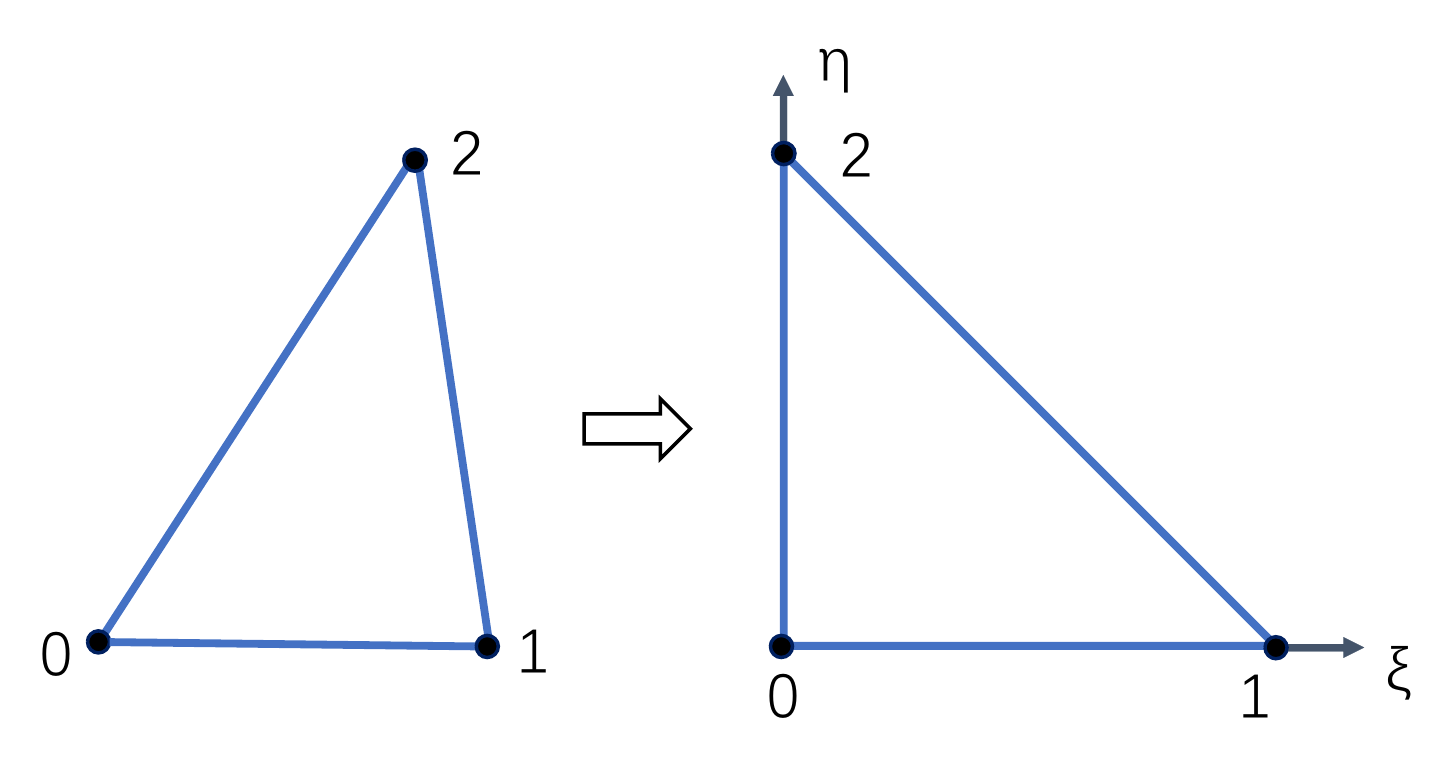}}
	\subfigure[Quadrilateral face]{
		\label{tetra}
		\includegraphics[width=0.48\textwidth]
		{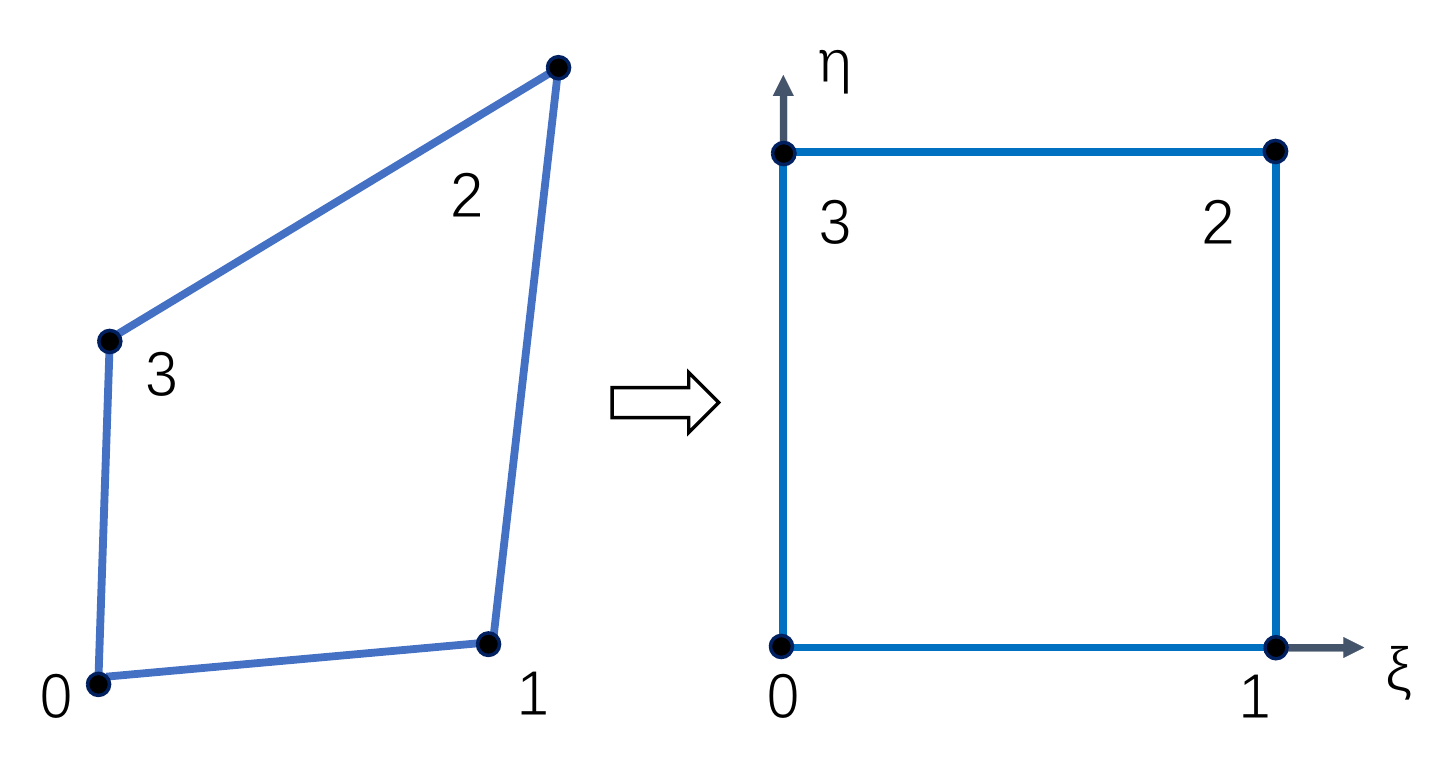}}
	\caption{The controlling points and isoparametric transformation of the cell interfaces.}
	\label{face}
\end{figure}

The semi-discretized form of FVM for conservation laws can be written as
\begin{equation}\label{semidiscrete}
\frac{\text{d} \textbf{W}_{i}}{\text{d}t}=\mathcal{L}(\textbf{W}_i)=-\frac{1}{\left| \Omega_i \right|} \sum_{p=1}^{N_f} \int_{\Gamma_{ip}}
\textbf{F}(\textbf{W}(\textbf{x},t))\cdot\textbf{n}_p \text{d}s,
\end{equation}
with
\begin{equation*}\label{f-to-flux-in-normal-direction}
\textbf{F}(\textbf{W}(\textbf{x},t))\cdot \textbf{n}_p=\int\pmb{\psi}  f(\textbf{x},t,\textbf{u},\xi) \textbf{u}\cdot \textbf{n}_p \text{d}\Xi,
\end{equation*}
where $\textbf{W}_{i}$ is the cell averaged values over cell $\Omega_i$, $\left|
\Omega_i \right|$ is the volume of $\Omega_i$, $\textbf{F}$ is the interface fluxes, and $\textbf{n}_p=(n_1,n_2,n_3)^T$ is the unit vector representing the outer normal direction of $\Gamma_{ip}$.
To evaluate the surface integral of fluxes,
the iso-parametric transformation is used, which can be written as
\begin{align*}
\textbf{X}(\xi, \eta)= \sum_{l=0}^{N_v} \textbf{x}_l \phi_l (\xi, \eta),
\end{align*}
where $\textbf{x}_l$ is the location of the mth vertex for each element
and $\phi_l$ is the base function \cite{wang2017thesis}.
In this work, the linear element is considered, and a schematic for the transformation is shown in Fig.~\ref{face}.

After the transformation, the Gaussian quadrature points can be determined and $\textbf{F}_{ip}(t)$ can be approximated by the numerical quadrature
	\begin{equation*}\label{fv-3d-general-quadrature}
\sum_{p=1}^{N_f} \int_{\Gamma_{ip}}
\textbf{F}(\textbf{W}(\textbf{x},t))\cdot\textbf{n}_p \text{d}s =  \left|\Gamma_{ip}\right| \sum_{k=1}^{M} \omega_k
	\textbf{F}(\textbf{x}_{p,k},t)\cdot\textbf{n}_p.
	\end{equation*}
To meet the requirement of a third-order spatial accuracy,
three Gaussian points are used for a triangular face and four Gaussian points are used for a quadrilateral face.
The details can be found in \cite{pan2020high,ji2021two}.
In the computation, the fluxes are obtained under the local coordinate. Details can be found in \cite{ji2021two}.

\subsection{Gas-kinetic solver}
Based on the integral solution of BGK equation \cite{xu2014direct}, a second-order time accurate gas distribution function at a local Gaussian point $\textbf{x}=(0,0,0)$ is constructed as
\begin{align}\label{2nd-flux}
f(\textbf{0},t,\textbf{u},\xi)
=&(1-e^{-t/\tau_n}) g^{c}+[(t+\tau)e^{-t/\tau_n}-\tau]a_{x_i}^{c}u_i g^{c}\nonumber
+(t-\tau+\tau e^{-t/\tau_n})A^{c}  g^{c}\nonumber\\
+&e^{-t/\tau_n}g^l[1-(\tau+t)a_{x_i}^{l}u_i-\tau A^l]H(u_1)\nonumber\\
+&e^{-t/\tau_n}g^r[1-(\tau+t)a_{x_i}^{r}u_i-\tau A^r] (1-H(u_1)).
\end{align}
The superscript $l,r$ represents the initial gas distribution function $f_0$, with a possible discontinuity at the left and right sides of a cell interface.
The superscript $c$ is the corresponding
equilibrium state $g$ in space and time.
The integral solution basically states a physical process from the particle free transport in $f_0$ in the kinetic scale
to the hydrodynamic flow evolution in the integral of $g$ term.
The flow evolution at the cell interface depends on the ratio of time step
to the  local particle collision time $\Delta t/\tau$.

The $g^k,~k=l,r$ has a form of a Maxwellian distribution
\begin{align*}
g^k=\rho^k (\frac{\lambda^k}{\pi})e^{-\lambda^k((u_i-U_i^k)^2+\xi^2)},
\end{align*}
which can be determined from the
 macroscopic variables $\textbf{W}
^l, \textbf{W}
^r$ through spatial reconstruction
\begin{align}\label{wlr_to_glr}
\int\pmb{\psi} g^{l}\text{d}\Xi=\textbf{W}
^l,\int\pmb{\psi} g^{r}\text{d}\Xi=\textbf{W}
^r.
\end{align}
The spatial and temporal microscopic derivatives are denoted as
\begin{align*}
a_{x_i} \equiv  (\partial g/\partial x_i)/g=g_{x_i}/g,
A \equiv (\partial g/\partial t)/g=g_t/g,
\end{align*}
which is determined by the spatial derivatives of macroscopic flow
variables and the compatibility condition as follows
\begin{align}\label{dwlr_to_dglr}
&\langle a_{x_1}\rangle =\frac{\partial \textbf{W} }{\partial x_1}=\textbf{W}_{x_1},
\langle a_{x_2}\rangle =\frac{\partial \textbf{W} }{\partial x_2}=\textbf{W}_{x_2},
\langle a_{x_3}\rangle =\frac{\partial \textbf{W} }{\partial x_3}=\textbf{W}_{x_3},\nonumber\\
&\langle A+a_{x_1}u_1+a_{x_2}u_2+a_{x_3}u_3\rangle=0,
\end{align}
where $\left\langle ... \right\rangle$ are the moments of a gas distribution function defined by
\begin{align*}
\langle (...) \rangle  = \int \pmb{\psi} (...) g \text{d} \Xi .
\end{align*}
Similarly, the
equilibrium state $ g^{c}$ and its derivatives $a_{x_i}^c, A_{x_i}^c$ are determined by the corresponding $\textbf{W}^c, \textbf{W}^c_{x_i}$.
The construction of the $\textbf{W}^c, \textbf{W}^c_{x_i}$ will be introduced in the next section.
The details for calculation of each microscopic term from macroscopic quantities can refer to \cite{ji2019high}.

For smooth flow, the time dependent solution in Eq.~\eqref{2nd-flux} can be simplified as \cite{GKS-2001}
\begin{align}\label{2nd-smooth-flux}
f(\textbf{0},t,\textbf{u},\xi)= g^{c}-\tau (a^c_{x_{i}} u_i +A^c)g^{c}+A^c g^{c}t,
\end{align}	
under the assumptions of $g^{l,r}=g^c$, $a^{l,r}_{x_i}=a^c_{x_i}$.
The above gas-kinetic solver for smooth flow has less numerical dissipations than the complete GKS solver in Eq.~\eqref{2nd-flux}.

In  smooth flow region, the collision time is determined by
\begin{align*}
\tau=\mu/p,
\end{align*}
where $\mu$ is the dynamic viscosity coefficient and $p$ is the pressure at the cell interface.
In order to properly capture the un-resolved discontinuities, additional numerical dissipation is needed.
The physical collision time $\tau$ in the exponential function part can be replaced by a numerical collision time $\tau_n$.
For the
inviscid flow, the collision time $\tau_n$ is modified as
\begin{align*}
\tau_n=\varepsilon \Delta t+C\displaystyle|\frac{p_l-p_r}{p_l+p_r}|\Delta
t,
\end{align*}
where $\varepsilon=0.01$ and $C=1$.
For the viscous flow, the collision time is related to the viscosity coefficient,
\begin{align*}
\tau_n=\frac{\mu}{p}+C \displaystyle|\frac{p_l-p_r}{p_l+p_r}|\Delta t,
\end{align*}
where $p_l$ and $p_r$ denote the pressure on the left and right
sides of the cell interface.

\subsection{Direct evolution of the cell averaged first-order spatial derivatives} \label{slope-section}

As shown in Eq.~\eqref{2nd-flux}, a time evolution solution at a cell interface is provided by the gas-kinetic solver, which is distinguished from the Riemann solvers with a constant solution.
Recall Eq.~\eqref{point}, the conservative variables at the Gaussian point  $\textbf{x}_{p,k}$ can be updated through the moments $\pmb{\psi}$
of the gas distribution function,
\begin{align}\label{point-interface}
\textbf{W}_{p,k}(t^{n+1})=\int \pmb{\psi} f^n(\textbf{x}_{p,k},t^{n+1},\textbf{u},\xi) \text{d}\Xi,~ k=1,...,M.
\end{align}

Then the cell-averaged first-order derivatives within each element at $t^{n+1}$ is given through the divergence theorem,
\begin{equation}\label{gauss-formula}
\begin{aligned}
\overline{W}_x^{n+1}
&=\frac{1}{ \Delta V}\int_{V} \nabla \cdot (\overline{W}(t^{n+1}),0,0) \text{d}V
=\frac{1}{ \Delta V}\int_{\partial V} (1,0,0) \cdot \textbf{n}  \overline{W}(t^{n+1}) \text{d}S \\
&=\frac{1}{ \Delta V}\int_{\partial V} \overline{W}(t^{n+1}) n_1   \text{d}S
=\frac{1}{ \Delta V} \sum_{p=1}^{N_f}\sum_{k=1}^{M_p} \omega_{p,k} W^{n+1}_{p,k} (n_1)_{p,k} \Delta S_p,
\\
\overline{W}_y^{n+1}
&=\frac{1}{ \Delta V}\int_{V} \nabla \cdot (0,\overline{W}(t^{n+1}),0) \text{d}V
= \frac{1}{ \Delta V}\int_{\partial V} (0,1,0) \cdot \textbf{n}  \overline{W}(t^{n+1}) \text{d}S \\
& = \frac{1}{ \Delta V}\int_{\partial V} \overline{W}(t^{n+1}) n_2   \text{d}S
=\frac{1}{ \Delta V} \sum_{p=1}^{N_f}\sum_{k=1}^{M_p} \omega_{p,k} W^{n+1}_{p,k} (n_2)_{p,k} \Delta S_p,
\\
\overline{W}_z^{n+1}
&=\frac{1}{ \Delta V}\int_{V} \nabla \cdot (0,0,\overline{W}(t^{n+1})) \text{d}V
=\frac{1}{ \Delta V}\int_{\partial V} (0,0,1) \cdot \textbf{n}  \overline{W}(t^{n+1}) \text{d}S \\
&=\frac{1}{ \Delta V}\int_{\partial V} \overline{W}(t^{n+1}) n_3 \text{d}S
=\frac{1}{ \Delta V} \sum_{p=1}^{N_f}\sum_{k=1}^{M_p} \omega_{p,k} W^{n+1}_{p,k} (n_3)_{p,k} \Delta S_p,
\end{aligned}
\end{equation}
where $\textbf{n}_{p,k}=((n_{1})_{p,k},(n_{2})_{p,k},(n_{3})_{p,k})$ is the outer unit normal direction at each Gaussian point $\textbf{x}_{p,k}$.

\subsection{Two-stage temporal discretization}

The two-stage fourth-order (S2O4) temporal discretization is
adopted here as that in the previous CGKS \cite{zhao2020compact,ji2021compact}.
Following the definition of Eq.~\eqref{semidiscrete},
a fourth-order time-accurate solution for cell-averaged conservative flow variables $\textbf{W}_i$ are updated by
\begin{equation}\label{s2o4}
\begin{aligned}
\textbf{W}_i^*&=\textbf{W}_i^n+\frac{1}{2}\Delta t\mathcal
{L}(\textbf{W}_i^n)+\frac{1}{8}\Delta t^2\frac{\partial}{\partial
	t}\mathcal{L}(\textbf{W}_i^n), \\
\textbf{W}_i^{n+1}&=\textbf{W}_i^n+\Delta t\mathcal
{L}(\textbf{W}_i^n)+\frac{1}{6}\Delta t^2\big(\frac{\partial}{\partial
	t}\mathcal{L}(\textbf{W}_i^n)+2\frac{\partial}{\partial
	t}\mathcal{L}(\textbf{W}_i^*)\big),
\end{aligned}
\end{equation}
where
$\mathcal{L}(\textbf{W}_i^n)$ and $\frac{\partial}{\partial t}\mathcal{L}(\textbf{W}_i^n)$ are given by
\begin{equation*} \label{flux-operator}
\begin{aligned}
\mathcal{L}(\textbf{W}_i^n)&= -\frac{1}{\left| \Omega_i \right|} \sum_{p=1}^{N_f}
\sum_{k=1}^{M} \omega_{p,k} \textbf{F}(\textbf{x}_{p,k},t_n)\cdot \textbf{n}_{p,k}\Delta S_p,\\
\frac{\partial}{\partial t}\mathcal{L}(\textbf{W}_i^n)&= -\frac{1}{\left| \Omega_i \right|} \sum_{p=1}^{N_f}\sum_{k=1}^{M} \omega_{p,k}
\partial_t \textbf{F}(\textbf{x}_{p,k},t_n)\cdot \textbf{n}_{p,k}\Delta S_p, \\
\frac{\partial}{\partial t}\mathcal{L}(\textbf{W}_{i}^*)&=-\frac{1}{\left| \Omega_i \right|} \sum_{p=1}^{N_f}\sum_{k=1}^{M} \omega_{p,k}
\partial_t \textbf{F}(\textbf{x}_{p,k},t_*)\cdot \textbf{n}_{p,k}\Delta S_p.
\end{aligned}
\end{equation*}
The proof for the fourth-order accuracy in time is shown in \cite{li2016twostage}.
The time dependent gas distribution function at a cell interface is updated in a similar way,
\begin{equation}\label{step-du}
\begin{split}
&f^*=f^n+\frac{1}{2}\Delta tf_t^n,\\
&f^{n+1}=f^n+\Delta tf_{t}^*.
\end{split}
\end{equation}
Thus, $f^*$ and $f^{n+1}$ are fully determined by Eq.~\eqref{2nd-flux} or Eq.~\eqref{2nd-smooth-flux} and the macroscopic flow variables and their fluxes at the cell interface can be obtained simultaneously by Eq.~\eqref{point} and Eq.~\eqref{f-to-flux}. The details can be found in \cite{zhao2020compact}.
A fourth-order temporal accuracy for the Euler equations can be achieved for the conservative flow variables on arbitrary mesh by Eq.~\eqref{s2o4} and Eq.~\eqref{step-du}. The complete proofs are given in \cite{li2016twostage,zhao2020compact}.

\section{Compact third-order reconstruction}

In this section, the details for the construction of the compact reconstruction for smooth and discontinuous flow are presented.
Especially, a special treatment, namely the cell-averaged slope compression factor, is introduced and improves significantly the robustness of CGKS for supersonic and hypersonic flow simulation under irregular mesh.

\subsection {Smooth reconstruction}
For a piecewise smooth function $Q( \textbf{x} )$ over cell $\Omega_{0}$, a
polynomial $P^r(\textbf{x})$ with degree $r$ can be constructed to
approximate $Q(\textbf{x})$ as follows
\begin{equation*}
P^r(\textbf{x})=Q(\textbf{x})+O(\Delta h^{r+1}),
\end{equation*}
where $\Delta h \sim |\Omega_{0}|^{\frac{1}{3}}$ is the equivalent cell size.
In order to achieve a third-order accuracy and satisfy conservative property,
the following quadratic polynomial over cell $\Omega_{0}$ is needed
\begin{equation*}\label{p2-def}
P^2(\textbf{x})= \overline{Q}_{0}+\sum_{|k|=1}^2a_kp^k(\textbf{x}),
\end{equation*}
where $\overline{Q}_{0}$ is the cell averaged value of $Q(\textbf{x})$ over cell $\Omega_{0}$, $k=(k_1,k_2,k_3)$, $|k| = k_1+k_2+k_3$.
The $p^k(\textbf{x})$ are basis functions, which are given by
\begin{align}\label{base}
\displaystyle p^k(\textbf{x})=x_1^{k_1}x_2^{k_2}x_3^{k_3}-\frac{1}{\left| \Omega_{0} \right|}\displaystyle\iiint_{\Omega_{0}}x_1^{k_1}x_2^{k_2}x_3^{k_3} \text{d}V.
\end{align}

The volume integral in Eq.\eqref{base} for a hexahedron $\Omega_{0}$ can be evaluated by the iso-parametric transformation described in \cite{ji2021compact}.
Other types of elements, i.e., tetrahedron, pyramid and prism, can be treated as the special cases of a hexahedron with some vertices merging together.

The quadratic polynomial $P^2(\textbf{x})$ on $\Omega_{0}$ is constructed on
the compact stencil $S_2$ including  $\Omega_{0}$ and all its von Neumann neighbors, $\Omega_{m}, m=1,...,N_f$,
where the averages of $Q(\textbf{x})$ and averaged derivatives of $Q(\textbf{x})$ over each cell are known.

The following values on $S_2$ are used to obtain $P^2(\textbf{x})$,
\begin{itemize}
	\item cell averages $\overline{Q}$ for cell $0,...,N_f$,
	\item cell averages of the $x$-direction partial derivative $\overline{Q}_{x_1}$ for cell $1,...,N_f$;
	\item cell averages of the $y$-direction partial derivative $\overline{Q}_{x_2}$ for cell $1,...,N_f$;
	\item cell averages of the $z$-direction partial derivative $\overline{Q}_{x_3}$ for cell $1,...,N_f$.
\end{itemize}

The polynomial $P^2(\textbf{x})$ is required to exactly satisfy
\begin{align*} \label{large-stenci-condition-1}
\iiint_{\Omega_{m}}P^2(\textbf{x})\text{d}V=\overline{Q}_{m}\left| \Omega_{m}\right|,
\end{align*}
where $\overline{Q}_m$ is the cell averaged value over $\Omega_{m},~ m=1,...,N_f$,
with the following condition satisfied in a least-square sense
\begin{equation*}\label{large-stenci-condition-2}
\begin{split}
\iiint_{\Omega_{m}}
\frac{\partial}{\partial x_1} P^2(\textbf{x})\text{d}V=(\overline{Q}_{x_1})_m|\Omega_{m}|,\\
\iiint_{\Omega_{m}}
\frac{\partial}{\partial x_2} P^2(\textbf{x})\text{d}V=(\overline{Q}_{x_2})_m|\Omega_{m}|,\\
\iiint_{\Omega_{m}}
\frac{\partial}{\partial x_3} P^2(\textbf{x})\text{d}V=(\overline{Q}_{x_3})_m|\Omega_{m}|,
\end{split}	
\end{equation*}
where  $\overline{Q}_{x_i}, i=1,2,3$  are the cell averaged directional derivatives over $\Omega_{m}$ in a global coordinate,
respectively.
The constrained least-square method is used to solve the above linear system \cite{li2014efficient}.
The left and right states $W^{l,r}$ provided by the reconstructed $P^2(\textbf{x})$ yield a linearly stable third-order CGKS on hybrid mesh, as validated in Section \ref{test-case}.

\subsection{Multi-resolution WENO procedure}
In order to deal with discontinuity, the multi-resolution WENO reconstruction \cite{zhu2020new} is applied.
Define three polynomials
\begin{equation}
\begin{split}
p_{2}(\textbf{x})&=\frac{1}{\gamma_{2,2}} P^2(\textbf{x})-\sum_{\ell=0}^{1} \frac{\gamma_{\ell, 2}}{\gamma_{2,2}} p_{\ell}(\textbf{x}),\\
p_{1}(\textbf{x})&=\frac{1}{\gamma_{1,1}} P^{1}(\textbf{x})-\frac{\gamma_{0,1}}{\gamma_{1,1}} P^{0}(\textbf{x}),\\
p_{0}(\textbf{x})&= P^0(\textbf{x}).
\end{split}
\end{equation}
For a third-order reconstruction, the second-order polynomial $P^{2}(\textbf{x})$ can be rewritten as
\begin{equation} \label{rewrite_expression-3rd}
P^2(\textbf{x})=\gamma_{2,2}p_{2}+\gamma_{1,2}p_{1}+\gamma_{0,2}p_{0},
\end{equation}
with arbitrary positive coefficients $\gamma_{m,n}$ satisfying $\gamma_{0,2}+\gamma_{1,2}+\gamma_{2,2}=1, \gamma_{0,1}+\gamma_{1,1}=1 $.
The coefficients are chosen as $\gamma_{2,2}:\gamma_{1,2}:\gamma_{0,2} = 100:1:6$, and
$\gamma_{1,1}:\gamma_{0,1}=1:6$ as suggested in \cite{zhu2020new}.

The first-order polynomial $P^1(\textbf{x})$ is determined solely from the targeted cell $\Omega_{0}$
\begin{itemize}
	\item cell averages $\overline{Q}$ and cell averages of the $x_i$-direction partial derivatives $\overline{Q}_{x_i}$, i=1,2,3 for $\Omega_{0}$.
\end{itemize}
Thus, the $P^1(\textbf{x})$ becomes
\begin{align*} \label{p1-stenci-condition}
P^1(\textbf{x})=\overline{Q}_{0} + \overline{Q}_{x_i} x_i,~i=1,2,3.
\end{align*}

The zeroth-order polynomial $P^0(\textbf{x})$ is simply determined by the cell-averaged conservative variables on the targeted cell $\Omega_{0}$ itself,
i.e. $P^0(\textbf{x}) = \overline{Q}_0$.

The smoothness indicators $\beta_{j}, j=1,2$ are defined as
\begin{equation}\label{smooth-indicator}
\beta_j=\sum_{|\alpha|=1}^{r_j}|\Omega|^{ \frac{2}{3}|\alpha|-1}\iiint_{\Omega}\big(D^{\alpha}P_j(\textbf{x})\big)^2
\text{d} V,
\end{equation}
where $\alpha$ is a multi-index and $D$ is the derivative operator, $r_1=1$, $r_2=2$.
The smoothness indicators in Taylor series at $(x_0,y_0)$ have the order
\begin{align*}
\beta_2&=O\{|\Omega_0|^{\frac{2}{3}}[1+O(|\Omega_0|^{\frac{2}{3}})]\}=O(|\Omega_0|)^{\frac{2}{3}} = O(h^2),\\
\beta_1&=O\{|\Omega_0|^{\frac{2}{3}}[1+O(|\Omega_0|^{\frac{1}{3}})]\}=O(|\Omega_0|)^{\frac{2}{3}} = O(h^2).
\end{align*}
Assuming a suitable $\beta_0$
\begin{align*}
\beta_0&=O\{|\Omega_0|^{\frac{2}{3}}[1+O(|\Omega_0|^{\frac{1}{3}})]\}=O(|\Omega_0|)^{\frac{2}{3}} = O(h^2),
\end{align*}
 a global smoothness indicator $\sigma$ similar to that in \cite{zhu2020new} can be defined
\begin{equation*}
\sigma^{3rd} = (\frac{1}{2}(|\beta_2-\beta_1|+|\beta_2-\beta_0|))^{\frac{4}{3}} = O(|\Omega_0|^{\frac{4}{3}}) = O(h^4).
\end{equation*}
Then, the corresponding non-linear weights are given by
\begin{equation}\label{non-linear-weight}
\begin{split}
&\omega_{m,n}=\gamma_{m,n}(1+\frac{\sigma}{\epsilon+\beta_{m}}), \\
&\bar{\omega}_{m,n}=\frac{\omega_{m,n}}{\sum \omega_{m,n}} = {\gamma}_{m,n} + O(h^2),
\end{split}
\end{equation}
where $m=0,1,2$ when $n=2$; $m=0,1$ when $n=1$, and $\epsilon$ takes $10^{-16}$ to avoid zero in the denominator.

Replacing $\gamma_{m,n}$ by the normalized non-linear weights $\bar{\omega}_{m,n}$ in Eq.~\eqref{rewrite_expression-3rd}, the final reconstructed polynomials are given by
\begin{equation} \label{final_weno_expression-3rd}
R^{3rd}(\textbf{x})=\bar{\omega}_{2,2}p_{2}+\bar{\omega}_{1,2}p_{1}+\bar{\omega}_{0,2}p_{0}.
\end{equation}
As a result, the non-linear reconstruction meets the requirement for a third-order accuracy  $R(\textbf{x})=P(\textbf{x})+O(h^3)$.
If any of these values yields negative density or pressure,  the first-order reconstruction is used instead.
The desired non-equilibrium states at Gaussian points can be obtained from the weighted polynomials
\begin{equation*} \label{final-weno-value}
Q^{l,r}_{p,k}=R^{l,r}(\textbf{x}_{p,k}), ~(Q^{l,r}_{x_i})_{p,k}= \frac{\partial R^{l,r}}{\partial {x_i}}(\textbf{x}_{p,k}).
\end{equation*}

In order to improve the robustness of the compact scheme, a two-step reconstruction has been designed and the smooth indicator of the zeroth-order polynomial $P^{0}(\textbf{x})$ as a non-linear combination of the first-order biased sub-stencil on tetrahedron mesh \cite{ji2021two}.
In this paper, however, these sub stencils are simply chosen as each neighboring cell
\begin{itemize}
	\item cell averages $\overline{Q}_j$ and cell averages of the $x_i$-direction partial derivatives $\overline{Q}_{x_i,j}$, i=1,2,3 for cell $j=1,...,N_f$.
\end{itemize}
The smooth indicators for each stencils are
\begin{equation*}\label{p0-beta}
\begin{split}
\beta_{0,j} = |\Omega_0|^{\frac{2}{3}} (b_{1,j}^2 + b_{2,j}^2 + b_{3,j}^2)
            =|\Omega_0|^{\frac{2}{3}}(\overline{Q}_{x_1,j}^2 + \overline{Q}_{x_2,j}^2 +            \overline{Q}_{x_3,j}^2).
\end{split}
\end{equation*}	

In this way, the sub-stencils are only single-cell involved. Thus, the smooth indicators are less affected by the local bad geometry.
In addition, the WENO procedure becomes more simple and less memory-consuming.
For the current reconstruction, only the coefficients for $P^{2}(\textbf{x})$ are stored.
Specially, a coefficient matrix with dimension $9\times M$ is needed, where $M=16$  for tetrahedron, $M=24$ for hexahedron.
The detailed formulation for the smooth indicator $\beta_0$ is given as
\begin{equation*}\label{p0-non-linear-weight}
\begin{split}
& \sigma^{1st} = [\frac{1}{ \frac{1}{2} N_f (N_f-1)}(\sum |\beta_{0,j}-\beta_{0,k}|)]^{\frac{4}{3}},\\
&\omega^{1st}_{j}=1+\frac{\sigma^{1st}}{\epsilon+\beta_{j}}, \\
&\bar{\omega}^{1st}_{j}=\frac{\omega_{j}}{\sum \omega_{j}},	\\
&\beta_{0} =\min (\sum \bar{\omega} ^{1st}_{j} \beta _{0,j},\beta _{0,0}),
\end{split}
\end{equation*}	
where $j,k=1,...,N_f$ and $j>k$.

However, such a choice is not robust enough for high-speed flow.
To improve the robustness of the spatial reconstruction for the CGKS, the CF will be introduced later.

\subsection{Reconstruction of the boundary cells} \label{boundary-recon}

One ghost cell is created for each boundary face by mirror symmetry of the corresponding inner cells. The cell-averaged quantities can be assigned according to the boundary condition.
Then the reconstruction for the inner cell can be determined.
After obtaining the inner state (assume as $\tilde{\textbf{W}}^r$) at a boundary Gaussian point, a ghost state (assume as $\tilde{\textbf{W}}^l$) can be assigned according to boundary condition under local coordinates.
There is a possible discontinuity between $\tilde{\textbf{W}}^l$ and $\tilde{\textbf{W}}^r$.
The ghost state setting for the Maxwell slip isothermal wall is given as follows.
	\begin{itemize}
		\item
		The slip wall assumption allows a discontinuity in velocities at the cell interface, and the corresponding gas distribution function is $f = f_0^r|_{u_1<0} + {f^l_0}|_{u_1>0}$, where the inner non-equilibrium state is
		\begin{equation*}
		\begin{aligned}
		f_0^r=e^{-t/\tau}g^r[1-\tau(a_{x_i}^{r}u_i+A^r)-ta^{r}_{x_i}u_i].
		\end{aligned}
		\end{equation*}
		A time-independent Maxwellian distribution  $f^l_0(\textbf{0},\textbf{u},t)=\rho^l (\frac{\lambda^l}{\pi})e^{-\lambda^l(\textbf{u}^2+\xi^2)}$ is assumed in the ghost state with zero macroscopic velocities ${U}_{i}^l=0$, a fixed temperature $T^l= 1/{(2R\lambda^l)}$, where $R$ is gas constant, and zero derivatives of the conservative variables   $\partial \textbf{W} ^l=0 $.
		\item Then,  $\rho^l$ is determined by no penetration constraint through the solid wall, which is obtained by solving the zero mass flux
		$ \int u_1f \text{d}\Xi = \int u_1 (f_0^r|_{u_1<0} + {f^l_0}|_{u_1>0}) \text{d}\Xi=0$.
		\item The conservative variables $\textbf{W}^{n+1}$ and the fluxes $\textbf{F}^{n+1}$ are given by the above $f$ at the boundary.
	\end{itemize}
The settings for the slipwall, non-slip adiabatic wall, and non-slip isothermal wall are described in \cite{ji2021two}.

\subsection{Reconstruction of equilibrium state}
The reconstructions for the non-equilibrium states  have the same spatial order of accuracy and can be used to get the equilibrium state
$ g^{c},g_{x_i}^{c}$ directly by a suitable average of $g^{l,r},g_{x_i}^{l,r}$.
To be consistent with the construction of  $ g^{c}$,
a kinetic-based weighting method is adopted
\begin{align}\label{wc_to_gc}
&\int\pmb{\psi} g^{c}\text{d}\Xi=\textbf{W}^c=\int_{u>0}\pmb{\psi}
g^{l}\text{d}\Xi+\int_{u<0}\pmb{\psi} g^{r}\text{d}\Xi, \nonumber \\
&\int\pmb{\psi} g^{c}_{x_i}\text{d}\Xi=\textbf{W}_{x_i}^c=\int_{u>0}\pmb{\psi}
g_{x_i}^{l}\text{d}\Xi+\int_{u<0}\pmb{\psi} g_{x_i}^{r}\text{d}\Xi.
\end{align}
The data for this method has compact support.
In programming, this procedure is included inside the subroutine of the gas distribution function, since it is performed at the local coordinate. Thus, it is also cache-friendly.
This method has been validated in the previous CGKS \cite{ji2021compact,ji2021two}.
In this way, all components of the microscopic slopes in Eq.~\eqref{2nd-flux} can be determined.

\subsection{Cell-averaged gradient compression factor}
A necessary condition for the divergence theorem  is the smoothness of the conservative variables inside the targeted cell.
When a discontinuity exists within the targeted cell, the cell averaged gradients obtained by Eq.~\eqref{gauss-formula} are not reliable.
Therefore, one idea to improve the robustness of the CGKS is to make the absolute value of the cell averaged gradients be small enough near discontinuities.
Then, the WENO reconstruction in Eq.~\eqref{final_weno_expression-3rd} will approach to the first order of accuracy at next step.
In this paper, the CF $\alpha _c$ at targeted cell $\Omega_{0}$ is defined as
\begin{equation}\label{cf}
\alpha _c
= \prod \limits_{i=0}^{N_f} \prod \limits_{k=0}^{M_i}  \alpha_{i,k},~~\alpha _c \in (0,1],
\end{equation}
where $\alpha_{i,k}$ is the gradient compression factor for the kth Gaussian point at the interface i around $\Omega_{0}$.
Then, the updated slope is then modified as
\begin{equation*}
\tilde{\textbf{W}}_{x_i}^{n+1}=\alpha _c \textbf{W}_{x_i}^{n+1},
\end{equation*}

The gradient compression factor $\alpha$ at a Gaussian point is defined as
\begin{equation}\label{cf-face}
\alpha = \frac{1}{1 + [\delta Q /(\delta \overline{Q}+\epsilon)]^{K_s} F},
\end{equation}
where 
$\overline{Q}$ is the cell averaged quantities,
$\delta Q$ is the absolute difference of the  left and right values at a cell interface,
\begin{equation*}
\delta Q = |Q^l - Q^r|.
\end{equation*}
$F$ is given as
\begin{equation*}
\begin{split}
&F=[C_1 D_p+C_2(D_{Ma,2}+D_{Ma,3})]^{K_t},~~F \in [0,\infty),\\
&D_p = |\frac{p^l-p^r}{p^l+p^r}|,\\
&D_{Ma,2}= |\frac{{Ma_2}^l-{Ma_2}^r}{{Ma_2}^l+{Ma_2}^r+\epsilon}|,\\
&D_{Ma,3}= |\frac{{Ma_3}^l-{Ma_3}^r}{{Ma_3}^l+{Ma_3}^r+\epsilon}|,
\end{split}
\end{equation*}
where $D_p$ is the related pressure difference, $Ma_2$ and $Ma_3$ are the Mach differences in two tangential directions.

$Q$ can be  density, pressure, or entropy, and $K_s \geq 2$.
The other parameters are suggested as $C_1=0.5 \sim 2 $, $C_2=0.1 \sim 1 $, $K_t=2 \sim 6$.
The parameters chosen in this paper are: $Q$ is the density, $K_s=2$, $C_1=1.5$, $C_2=0.2$, $K_t=4$.

\begin{rmk}
	When flow around $\Omega_{0}$ is smooth,
	\begin{equation*}
	\begin{split}
	&[\delta Q /(\delta \overline{Q})+\epsilon)]^{K_s} = O(\Delta x)^{ p K_s},\\
	&F=[C_1 O(\Delta x)^{ p+1}+C_2O(\Delta x)^{ p+1}]^{K_t} = O(\Delta x)^{ (p+1) K_t},
	\end{split}
	\end{equation*}
where $p$ is the order of the reconstructed polynomial.
	Thus, recall Eq.~\eqref{cf-face} and Eq.~\eqref{cf},
	\begin{equation*}
	\begin{split}
	&\alpha_{i,k} \to 1 + O(\Delta x)^{p K_s+ (p+1) K_t},\\
	&\alpha _c
	= \prod \limits_{i=0}^n \prod \limits_{k=0}^{M_i}  \alpha_{i,k} \to 1 + O(\Delta x)^{p K_s+ (p+1) K_t},
	\end{split}
	\end{equation*}
	which suggests $\alpha _c \sim 1$ if the flow is smooth.
	
	When flow around $\Omega_{0}$ contains discontinuity,
	\begin{equation*}
	\begin{split}
	&[\delta Q /(\delta \overline{Q})+\epsilon)]^{K_s} \to O(1),\\
	& D_p \to 1,~~~D_{Ma,2} \to 1,~~~D_{Ma,3} \to 1,\\
	&F\to(C_1+C_2)^{K_t},
	\end{split}
    \end{equation*}
As a result,
	\begin{equation*}
	\begin{split}
	&\alpha_{i,k} \to \frac{1}{1 + |O(1)|},\\
	&\alpha _c
	= \prod \limits_{i=0}^n \prod \limits_{k=0}^{M_i}  \alpha_{i,k}
	\to \frac{1}{1 + |O(1)|}.
	\end{split}
	\end{equation*}
	When the discontinuity is strong, the final $\alpha _c$ can approach to $0$.
\end{rmk}

Ideally, it is hoped that the CGKS with the CF can give both accurate and robust results for the flow region from the subsonic to the hypersonic one.
For this reason, four representative test cases are used to evaluate its performance in the following section.
\begin{enumerate}
	\item  Sin-wave accuracy test, where the flow is smooth and the analytical solution exists.
	\item Shu-Osher test cases, which contains shock-vortex interactions and unsteady linear wave propagations.
	\item Flow over a sphere with Ma=2, Re=300. The case is steady and quantitative results including the drag force coefficient can be used to evaluate the influence of the CF.
	\item Hypersonic viscous flow over a space vehicle model with Ma=10. This demo is adopted to demonstrate the robustness of the CGKS.
\end{enumerate}

\section{Numerical examples} \label{test-case}

\begin{figure}[hbt!]	
	\centering
	\includegraphics[width=0.96\textwidth]
	{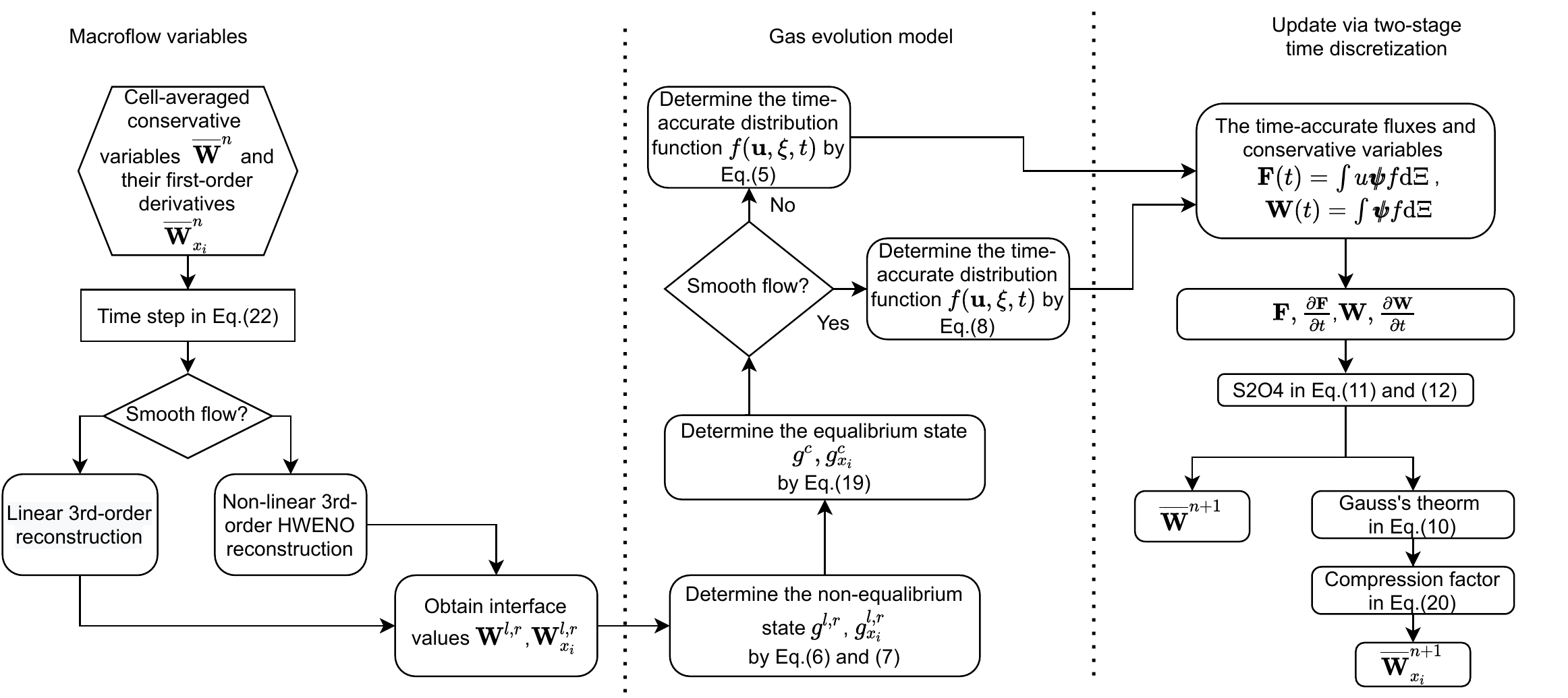}
	\caption{ The brief algorithm of the CGKS.}
	\label{flowchart}
\end{figure}

In this section, numerical tests will be presented to validate the proposed scheme.
The time step is determined by
\begin{align}\label{cfl-condition}
\Delta t = C_{CFL} \mbox{Min} ( \frac{ \Delta r_i}{||\textbf{U}_i||+(a_s)_i}, \frac{ (\Delta r_i)^2}{3\nu _i}),
\end{align}
where $C_{CFL}$ is the CFL number, and $||\textbf{U}_i||$, $(a_s)_i$, and $\nu _i= (\mu /\rho) _i$ are the magnitude of velocities, sound speed, and kinematic viscosity coefficient for cell i. The $\Delta r_i$ is taken as
\begin{align*}
\Delta r_i = \frac{3|\Omega_i|}{\sum|\Gamma_{ip}|},
\end{align*}
for a tetrahedron or pyramid,
and
\begin{align*}
\Delta r_i = \frac{|\Omega_i|}{\max|\Gamma_{ip}|},
\end{align*}
for a hexahedron or prism. The CFL number is taken as 0.5 if no specified.
An algorithm flowchart of the CGKS is given in Fig.~\ref{flowchart}.
In the present work, the WENO reconstruction based on the conservative variables and the complete flux in Eq.~\eqref{2nd-flux} are adopted without using the CF if no specified.

\subsection{3-D sinusoidal wave propagation}
\begin{figure}[htp]	
	\centering
	\includegraphics[width=0.44\textwidth]
	{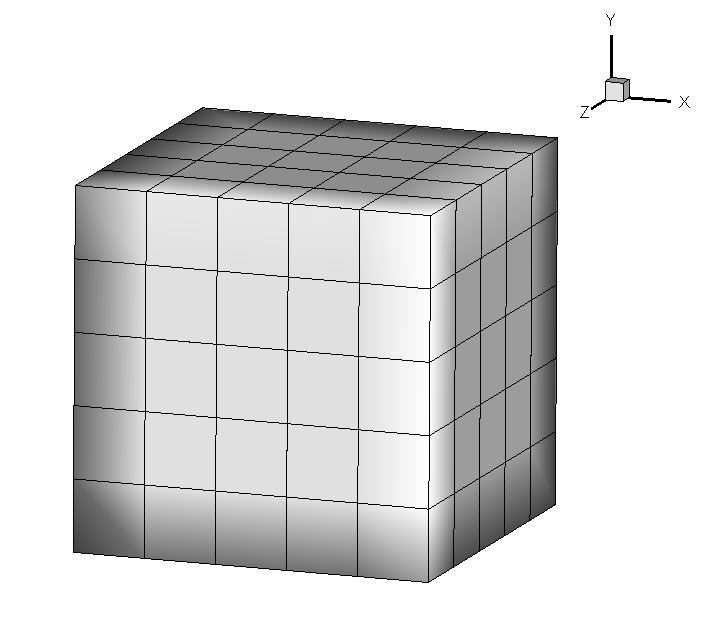}
	\includegraphics[width=0.44\textwidth]
	{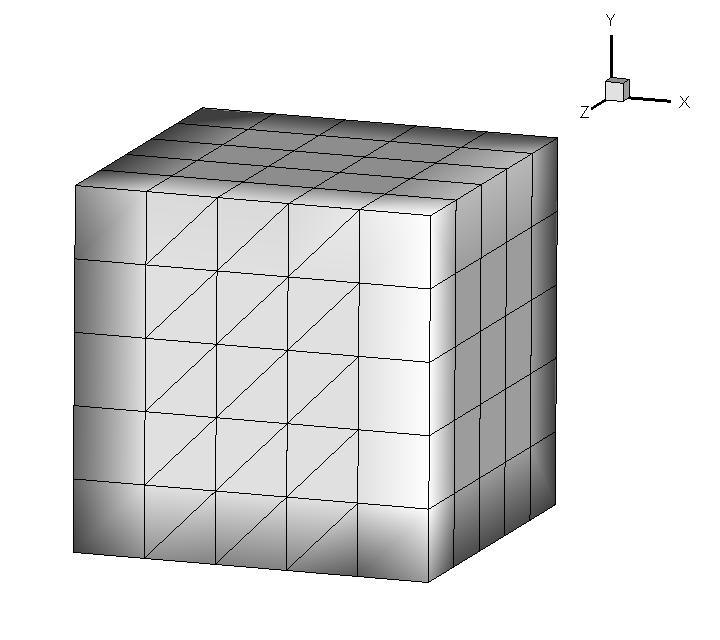}		
	\caption{Mesh sample for the 3D sin-wave
		propagation. Left: hexahedral mesh. Right: Hybrid mesh.}
	\label{sinwave-mesh}
\end{figure}

The initial
condition for the advection of density perturbation is given as
\begin{align*}
\rho(x,y,z)=1+0.2\sin(\pi (x+y+z)),\\
 \textbf{U}(x,y,z)=(1,1,1),~~~  p(x,y,z)=1,
\end{align*}
within a cubic domain $[0, 2]\times[0, 2]\times[0, 2]$.
A series of sequentially refined hexahedral meshes and hybrid meshes are used in the test, as shown in Fig.~\ref{sinwave-mesh}.
With the periodic boundary condition in all directions, the analytic
solution is
\begin{align*}
\rho(x,y,z,t)=1+0.2\sin(\pi(x+y+z-t)),\\
\textbf{U}(x,y,z)=(1,1,1),~~~~~~~~p(x,y,z,t)=1.
\end{align*}
The flow is inviscid and the collision time $\tau$ is $0$.
The $L^1$, $L^2$ and $L^{\infty}$ errors and the corresponding orders with linear weights at $t=2$ under both meshes are given in Tab.~\ref{3d-accuracy-linear-hex} and Tab.~\ref{3d-accuracy-linear-hy}.
The results with non-linear Z-type weights for uniform meshes are given in Tab.~\ref{3d-accuracy-weno-hex}.
Expected accuracy is achieved for the above cases.
Then, the results for the CGKS with the CF are listed in  Tab.~\ref{3d-accuracy-cf-hex}.
In comparison with the results only without the CF, slightly larger absolute errors are observed in Tab.~\ref{3d-accuracy-cf-hy}.
The same conclusion can be drawn for the cases under hybrid meshes, as shown in Tab.~\ref{3d-accuracy-cf-hy}.
The third-order accuracy is kept for both cases.

\begin{table}[htp]
	\small
	\begin{center}
		\def\temptablewidth{1\textwidth}
		{\rule{\temptablewidth}{1pt}}
		\begin{tabular*}{\temptablewidth}{@{\extracolsep{\fill}}c|cc|cc|cc}	
			Mesh number & $L^1$ error & Order & $L^2$ error & Order& $L^{\infty}$ error & Order  \\
			\hline
$5^3$ & 8.572800e-02 & ~ & 9.508770e-02 & ~ & 1.324228e-01 & ~\\
$10^3$ & 2.199962e-02 & 1.96 & 2.441038e-02 & 1.96 & 3.420145e-02 & 1.95 \\
$20^3$ & 3.083322e-03 & 2.83 & 3.431975e-03 & 2.83 & 5.052936e-03 & 2.76 \\
$40^3$ & 3.948944e-04 & 2.96 & 4.377654e-04 & 2.97 & 6.581416e-04 & 2.94 \\
		\end{tabular*}
		{\rule{\temptablewidth}{0.1pt}}
	\end{center}
	\vspace{-4mm} \caption{\label{3d-accuracy-linear-hex}
Accuracy test for the 3D sin-wave
propagation by the linear third-order compact reconstruction. Uniform hexahedral mesh.  }
\end{table}

\begin{table}[htp]
	\small
	\begin{center}
		\def\temptablewidth{1\textwidth}
		{\rule{\temptablewidth}{1pt}}
		\begin{tabular*}{\temptablewidth}{@{\extracolsep{\fill}}c|cc|cc|cc}
			
			mesh number & $L^1$ error & Order & $L^2$ error & Order& $L^{\infty}$ error & Order  \\
			\hline
			$1.6 \times 5^3$ & 9.697442e-02 & ~ & 1.080211e-01 & ~ & 1.512287e-01 & ~ \\
			$1.6 \times 10^3$  & 2.821966e-02 & 1.78 & 3.151973e-02 & 1.77 & 4.422013e-02 & 1.77 \\
			$1.6 \times 20^3$  & 4.036691e-03 & 2.81 & 4.496854e-03 & 2.81 & 6.650040e-03 & 2.73 \\
			$1.6 \times 40^3$  & 5.168948e-04 & 2.97 & 5.729851e-04 & 2.97 & 8.596728e-04 & 2.95 \\
		\end{tabular*}
		{\rule{\temptablewidth}{0.1pt}}
	\end{center}
	\vspace{-4mm} \caption{\label{3d-accuracy-linear-hy} Accuracy test for the 3D sin-wave
		propagation by the linear third-order compact reconstruction. Hybrid mesh.   }
\end{table}

\begin{table}[htp]
	\small
	\begin{center}
		\def\temptablewidth{1\textwidth}
		{\rule{\temptablewidth}{1pt}}
		\begin{tabular*}{\temptablewidth}{@{\extracolsep{\fill}}c|cc|cc|cc}
			
			mesh number & $L^1$ error & Order & $L^2$ error & Order& $L^{\infty}$ error & Order  \\
			\hline
			$5^3$ & 9.307715e-02 & ~ & 1.026267e-01 & ~ & 1.431789e-01 & ~ \\
			$10^3$  & 1.782439e-02 & 2.38 & 2.040047e-02 & 2.33 & 3.581868e-02 & 2.00 \\
			$20^3$  & 2.988300e-03 & 2.58 & 3.592716e-03 & 2.51 & 9.255122e-03 & 1.96 \\
			$40^3$  & 4.108034e-04 & 2.86 & 5.008456e-04 & 2.84 & 1.114581e-03 & 3.05 \\
		\end{tabular*}
		{\rule{\temptablewidth}{0.1pt}}
	\end{center}
	\vspace{-4mm} \caption{\label{3d-accuracy-weno-hex} Accuracy test for the 3D sin-wave
		propagation by the third-order compact WENO reconstruction with $d_0:d_1:d_2=100:1:6$. Uniform hexahedral mesh.   }
\end{table}

\begin{table}[htp]
	\small
	\begin{center}
		\def\temptablewidth{1\textwidth}
		{\rule{\temptablewidth}{1pt}}
		\begin{tabular*}{\temptablewidth}{@{\extracolsep{\fill}}c|cc|cc|cc}
			
			mesh number & $L^1$ error & Order & $L^2$ error & Order& $L^{\infty}$ error & Order  \\
			\hline
			$5^3$ & 9.177411e-02 & ~ & 1.014747e-01 & ~ & 1.421968e-01 & ~ \\
			$10^3$  & 1.783722e-02 & 2.36 & 2.041512e-02 & 2.31 & 3.586210e-02 & 1.99 \\
			$20^3$  & 2.827416e-03 & 2.65 & 3.717858e-03 & 2.48 & 1.140440e-02 & 1.65 \\
			$40^3$  & 4.073909e-04 & 2.80 & 5.042275e-04 & 2.88 & 1.228363e-03 & 3.21 \\
		\end{tabular*}
		{\rule{\temptablewidth}{0.1pt}}
	\end{center}
	\vspace{-4mm} \caption{\label{3d-accuracy-cf-hex} Accuracy test for the 3D sin-wave
		propagation by the third-order compact WENO reconstruction with $d_0:d_1:d_2=100:1:6$ and the CF. Uniform hexahedral mesh.   }
\end{table}

\begin{table}[htp]
	\small
	\begin{center}
		\def\temptablewidth{1\textwidth}
		{\rule{\temptablewidth}{1pt}}
		\begin{tabular*}{\temptablewidth}{@{\extracolsep{\fill}}c|cc|cc|cc}
			
			mesh number & $L^1$ error & Order & $L^2$ error & Order& $L^{\infty}$ error & Order  \\
			\hline
			$1.6 \times 5^3$ & 9.835342e-02 & ~ & 1.105509e-01 & ~ & 1.534561e-01 & ~ \\
			$1.6 \times 10^3$  & 4.301722e-02 & 1.19 & 4.875222e-02 & 1.18 & 7.870885e-02 & 0.96 \\
			$1.6 \times 20^3$  & 4.335145e-03 & 3.31 & 5.664078e-03 & 3.11 & 1.174496e-02 & 2.74 \\
			$1.6 \times 40^3$  & 5.491532e-04 & 2.98 & 6.424909e-04 & 3.14 & 1.412450e-03 & 3.06 \\
		\end{tabular*}
		{\rule{\temptablewidth}{0.1pt}}
	\end{center}
	\vspace{-4mm} \caption{\label{3d-accuracy-cf-hy} Accuracy test for the 3D sin-wave
		propagation by the third-order compact WENO reconstruction with $d_0:d_1:d_2=100:1:6$ and the CF. Hybrid mesh.   }
\end{table}

\subsection{Shu-Osher problem}
\begin{figure}[htp]	
	\centering
	\includegraphics[width=0.44\textwidth]
	{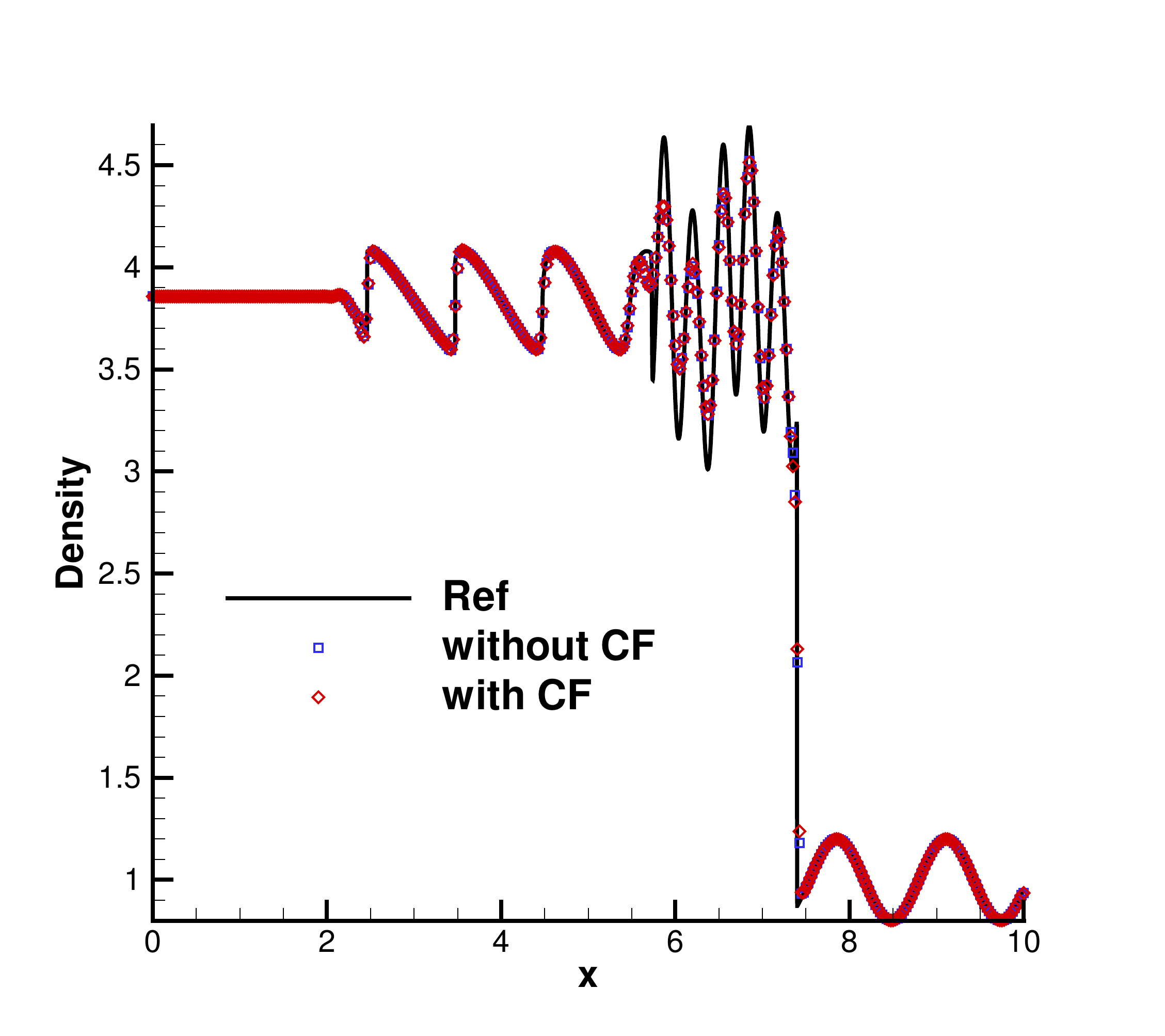}
	\includegraphics[width=0.44\textwidth]
	{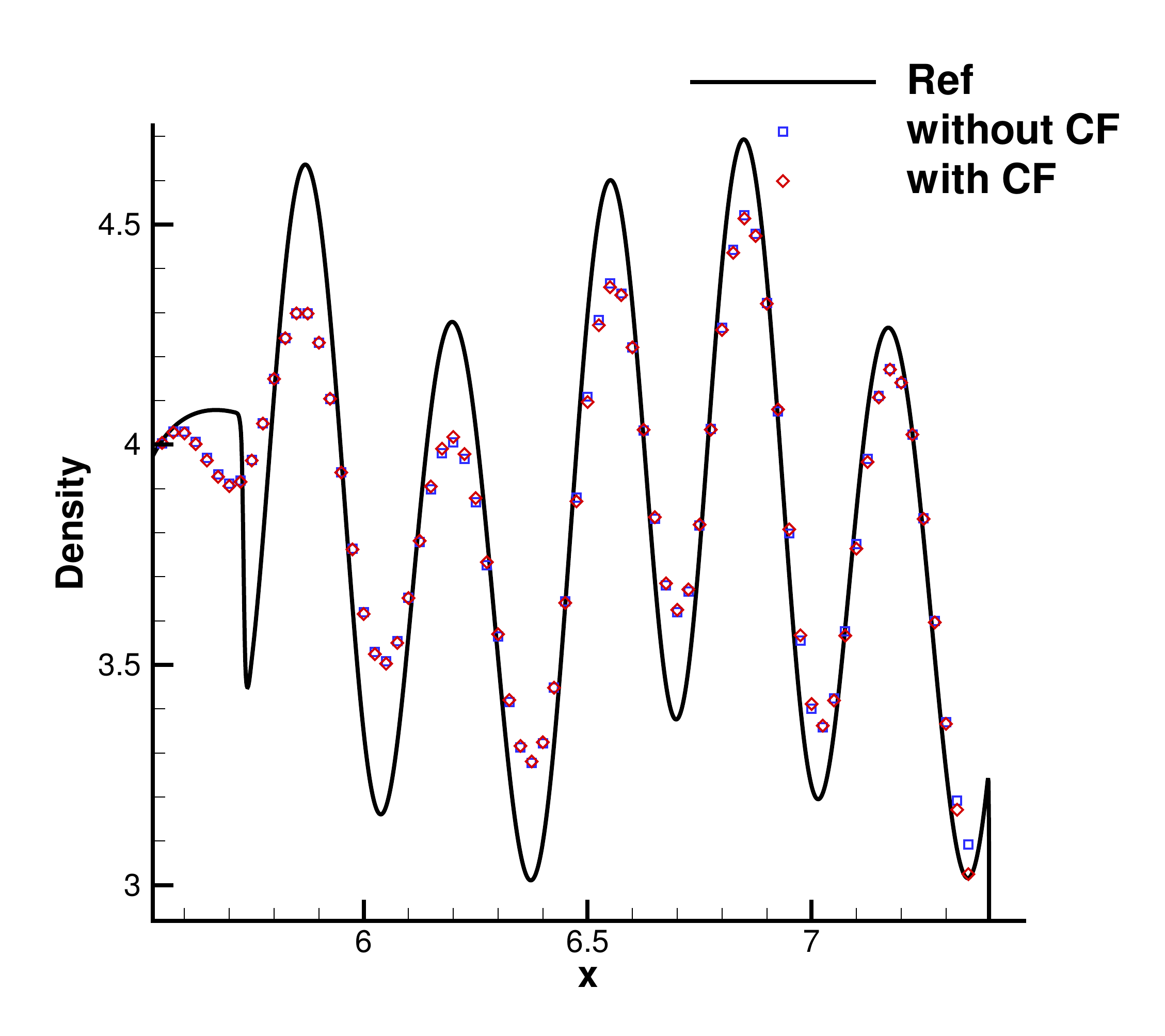}	
	\caption{Shu-Osher problem. Mesh number: $400 \times 2 \times 2$. }
	\label{shuosher-3d-mesh400}
\end{figure}

\begin{figure}[htp]	
	\centering
	\includegraphics[width=0.44\textwidth]
	{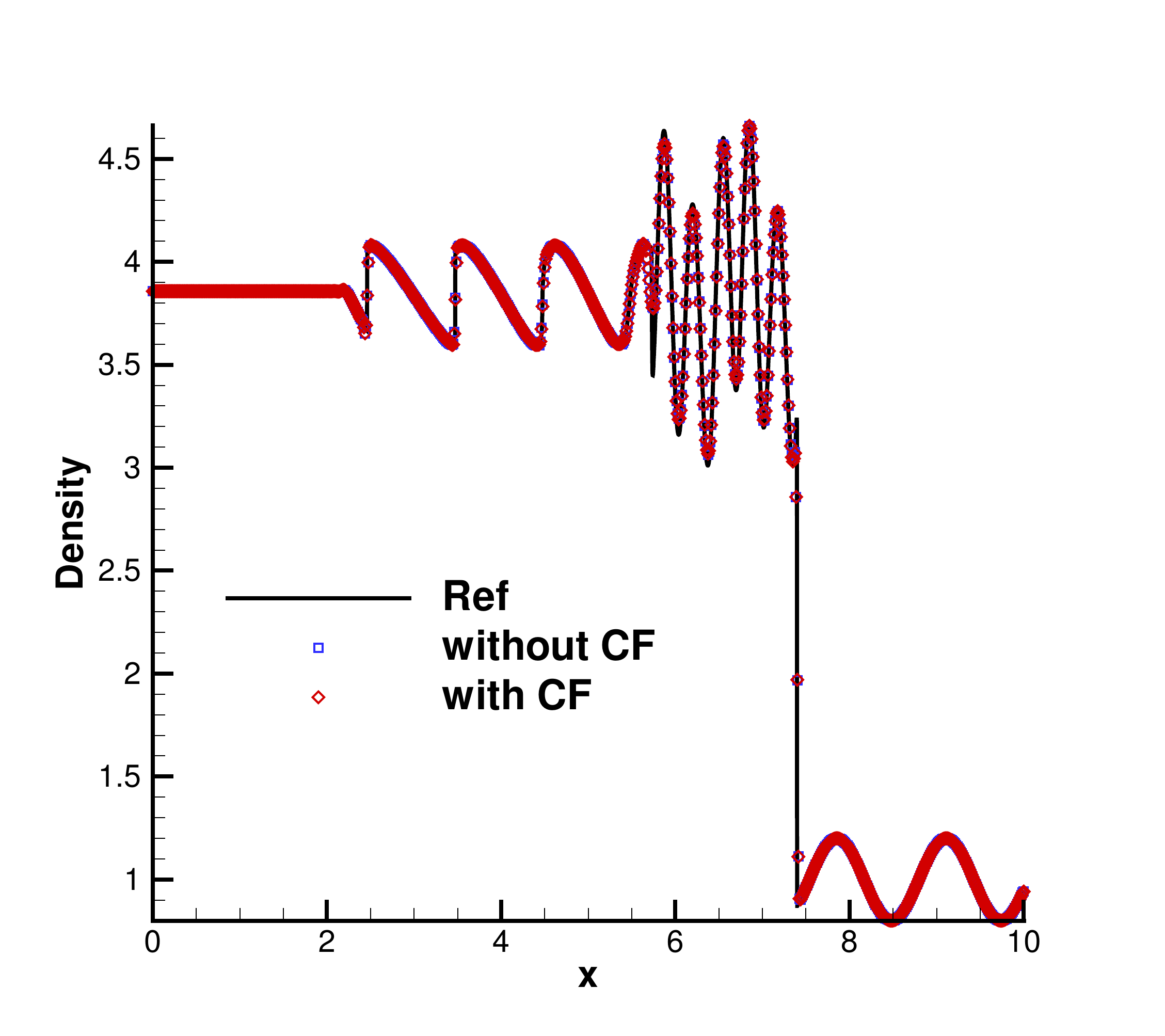}
	\includegraphics[width=0.44\textwidth]
	{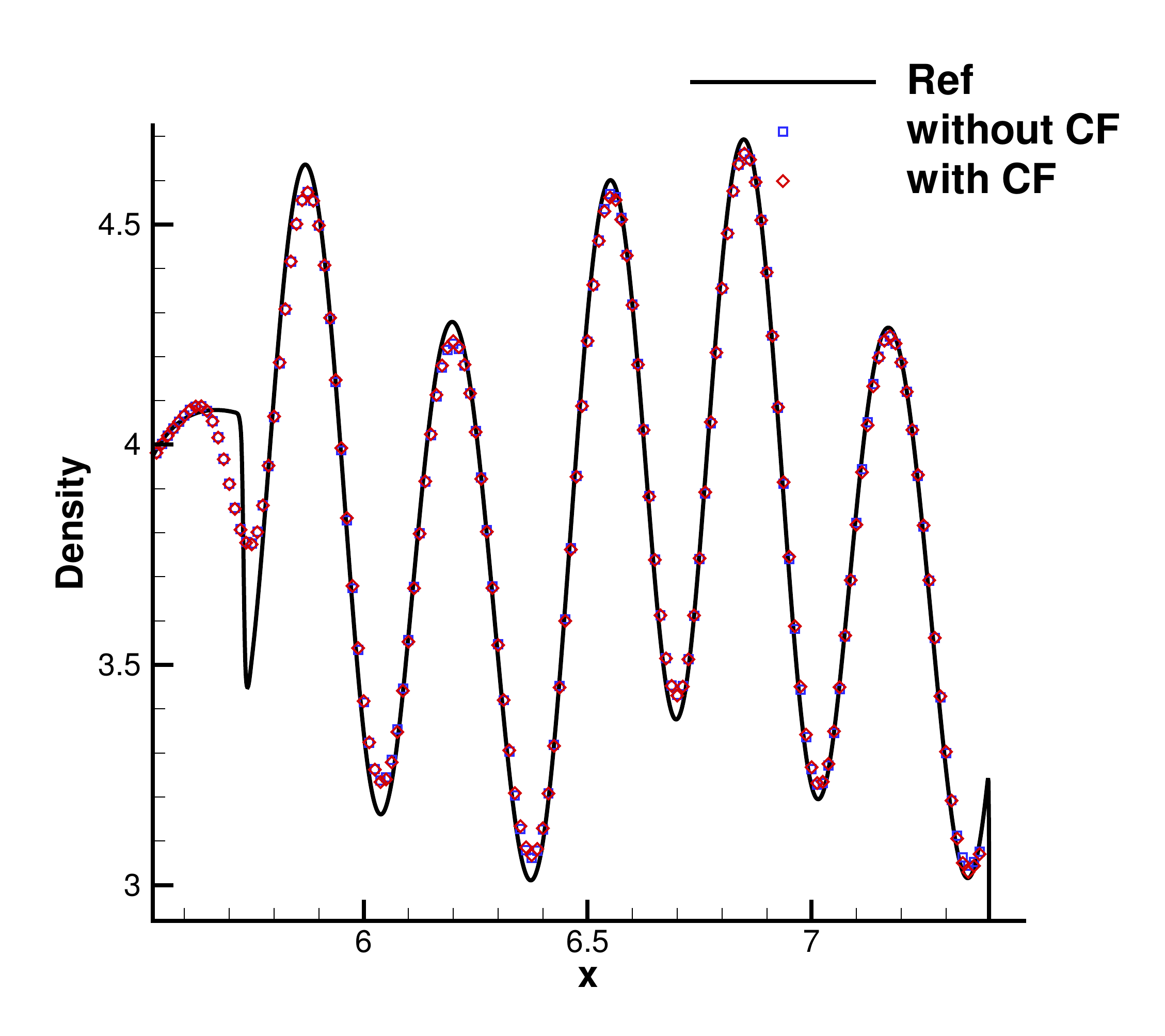}		
	\caption{Shu-Osher problem. Mesh number: $800 \times 2 \times 2$. }
	\label{shuosher-3d-mesh800}
\end{figure}

\begin{figure}[htp]	
	\centering
	\includegraphics[width=0.44\textwidth]
	{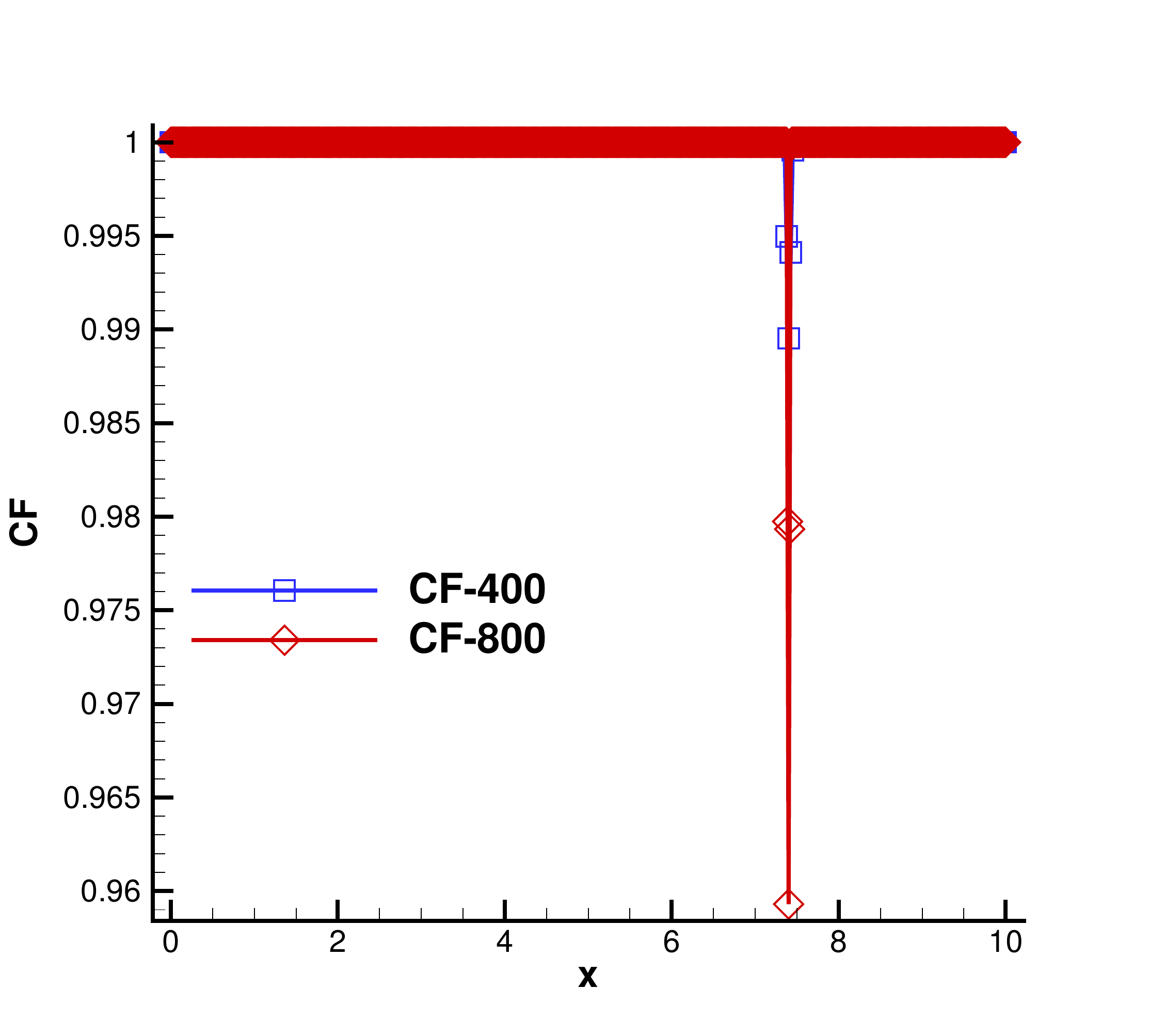}
	\includegraphics[width=0.44\textwidth]
	{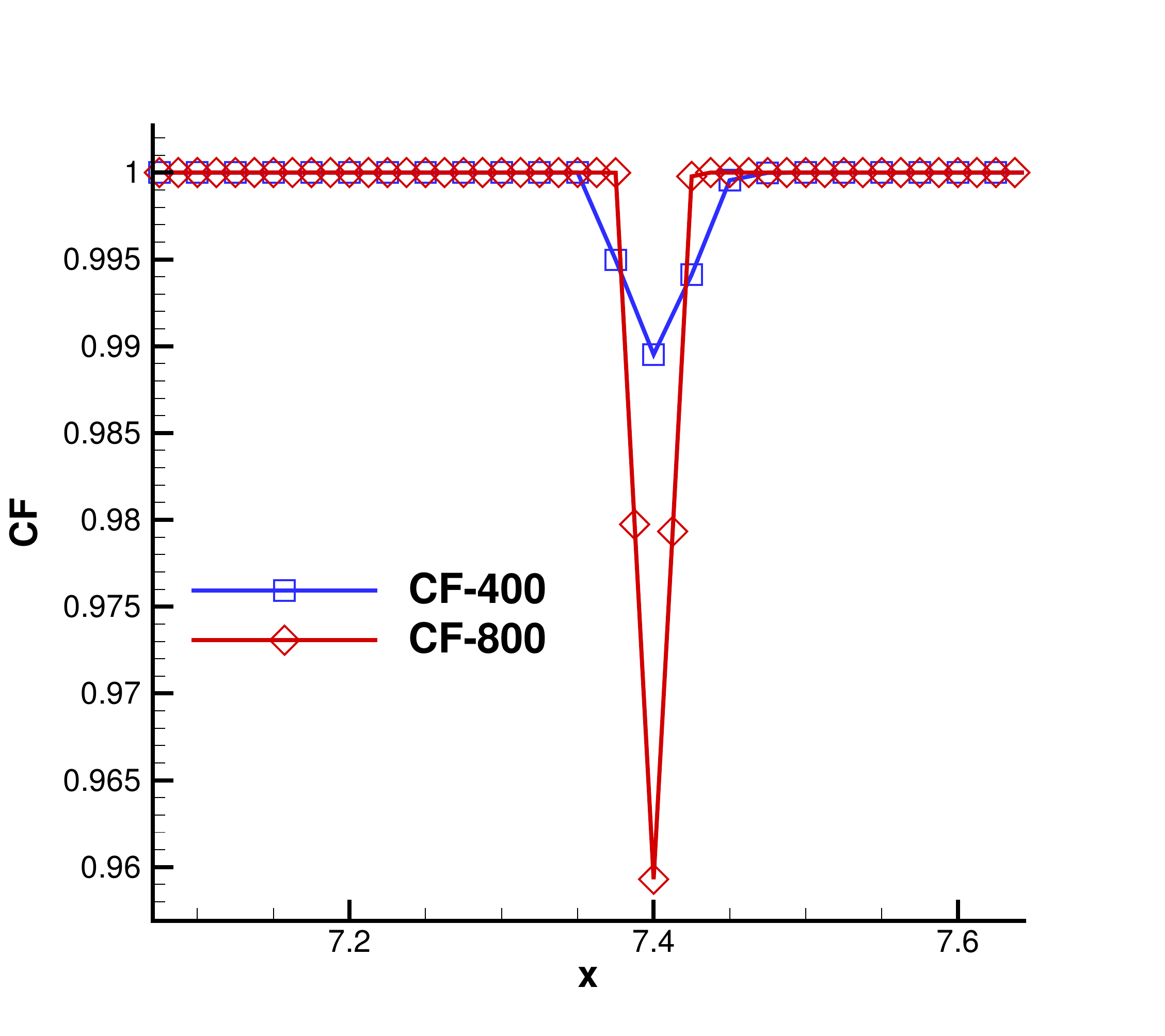}		
	\caption{The CF distributions and their local enlargements for Shu-Osher problem.}
	\label{shuosher-3d-cf}
\end{figure}

\begin{figure}[htp]	
	\centering
	\includegraphics[width=0.44\textwidth]
	{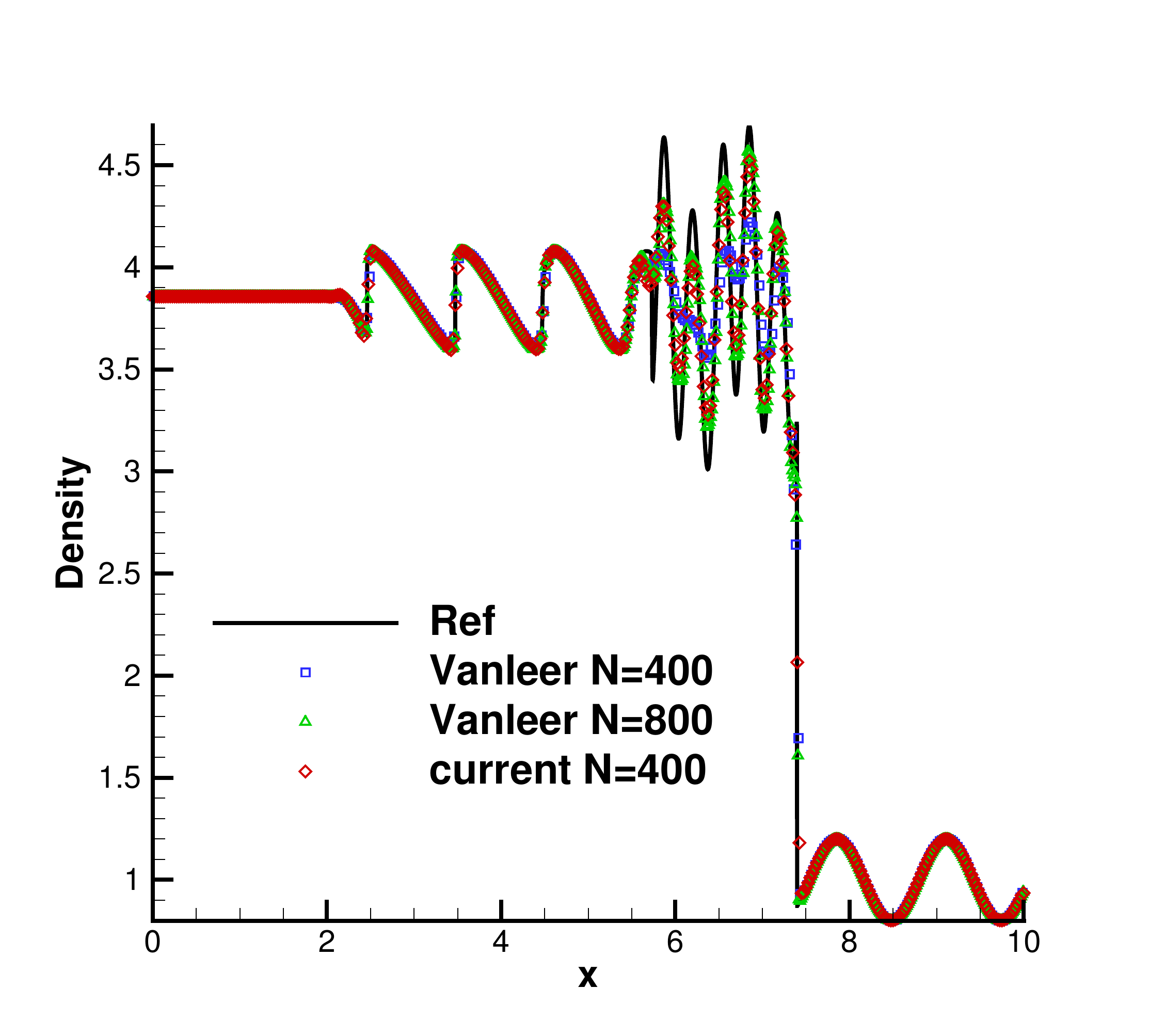}
	\includegraphics[width=0.44\textwidth]
	{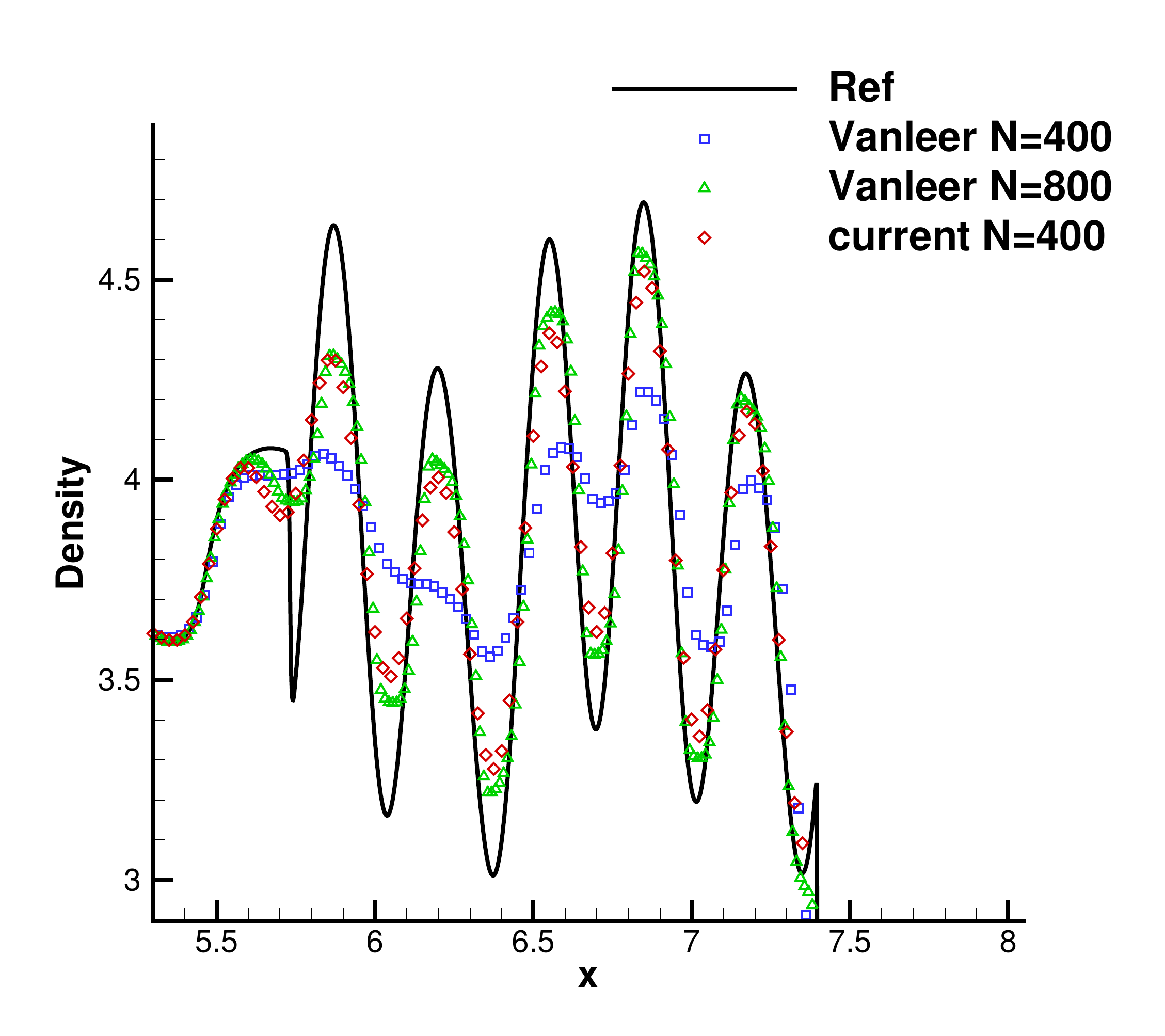}		
	\caption{Shu-Osher problem. In comparison with 1-D second-order GKS.}
	\label{shuosher-3d-compare}
\end{figure}

The initial condition for the Shu-Osher problem \cite{shu1989efficient} is
\begin{align*}
(\rho,U,p) =\begin{cases}
(3.857134, 2.629369, 10.33333), &  0<x \leq 1,\\
(1 + 0.2\sin (5x), 0, 1),  &  1 <x<10.
\end{cases}
\end{align*}
The flow is one-dimensional along the x-axis, and two uniform hexahedral meshes with a fixed length $L=10$ in x-direction are used in the computation.
The mesh sizes are $\Delta x = 1/40$ and $1/80$ respectively.
The fixed wave profile is extended on the right while the non-reflecting boundary condition is given on the left.
The computed density profiles and local enlargements for the Shu-Osher problem at $t = 1.8$ with both meshes are plotted in Fig.~\ref{shuosher-3d-mesh400} and Fiq.~\ref{shuosher-3d-mesh800}. The CGKS with/without the CF can both resolve the linear wave nicely.
To get a better understanding of the behavior of the CF in this test, the distributions for the CF at $t=1.8$ in each cell are plotted in Fig.~\ref{shuosher-3d-cf}.
It can be observed that the CF only takes effect near the normal shock.
The results are also compared with the 1-D second-order GKS based on the van Leer limiter. A similar resolution can be obtained by the current scheme with only half of the total mesh points used by a second-order method, as shown in Fig.~\ref{shuosher-3d-compare}.

\subsection{Laminar boundary layer}
\begin{figure}[htp]	
	\centering	
	\includegraphics[width=0.4\textwidth]
	{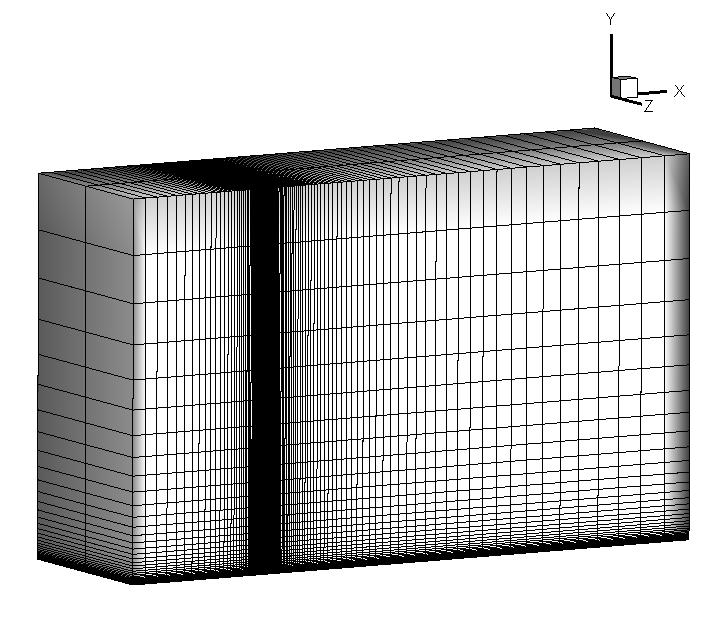}
	\includegraphics[width=0.4\textwidth]
	{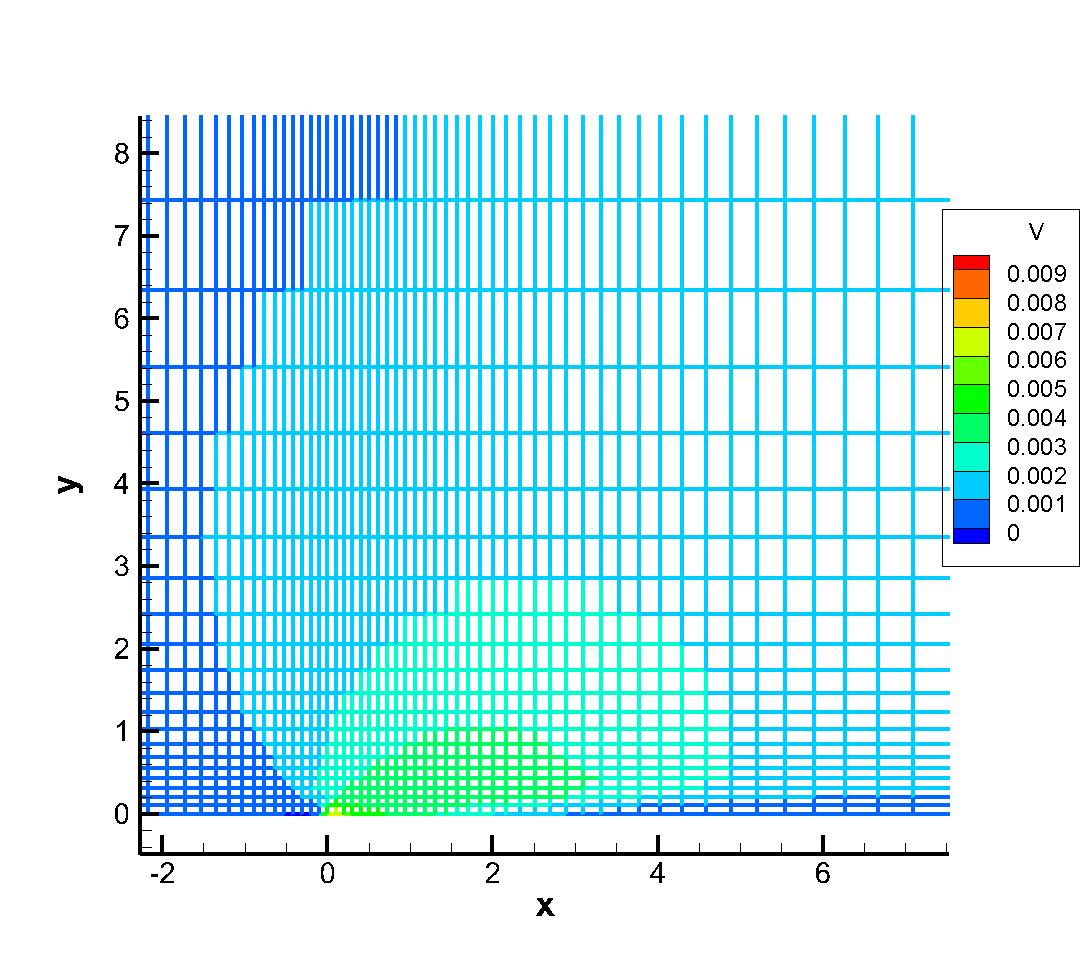}
	\vspace{-4mm} \caption{\label{flate-plate-mesh}
		Laminar boundary layer. Left: mesh with $120 \times 35 \times 2$  cells. Right: local enlargement colored by V velocities.}
\end{figure}

\begin{figure}[htp]	
	\centering
	\includegraphics[width=0.3\textwidth]
	{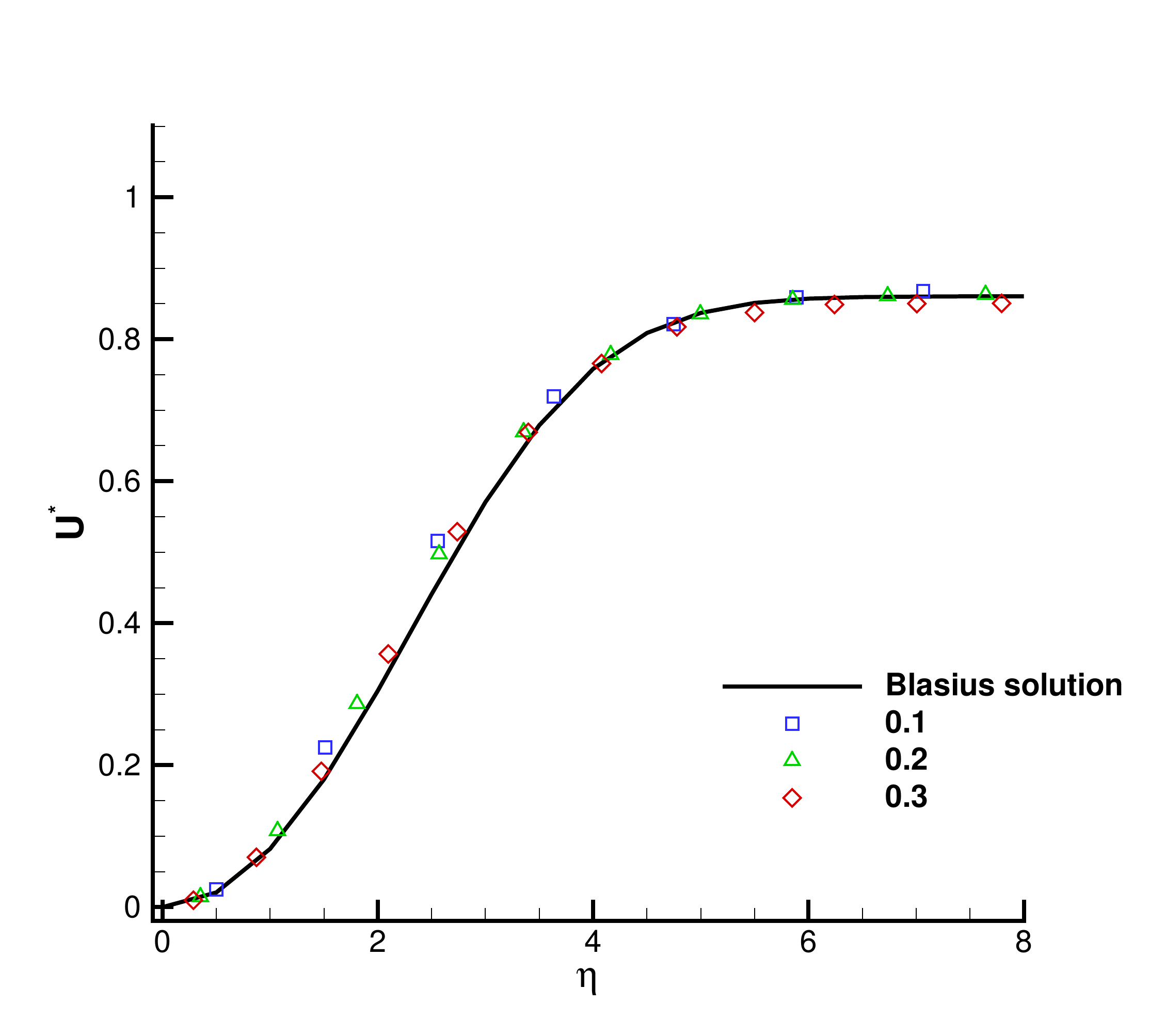}	
	\includegraphics[width=0.3\textwidth]
	{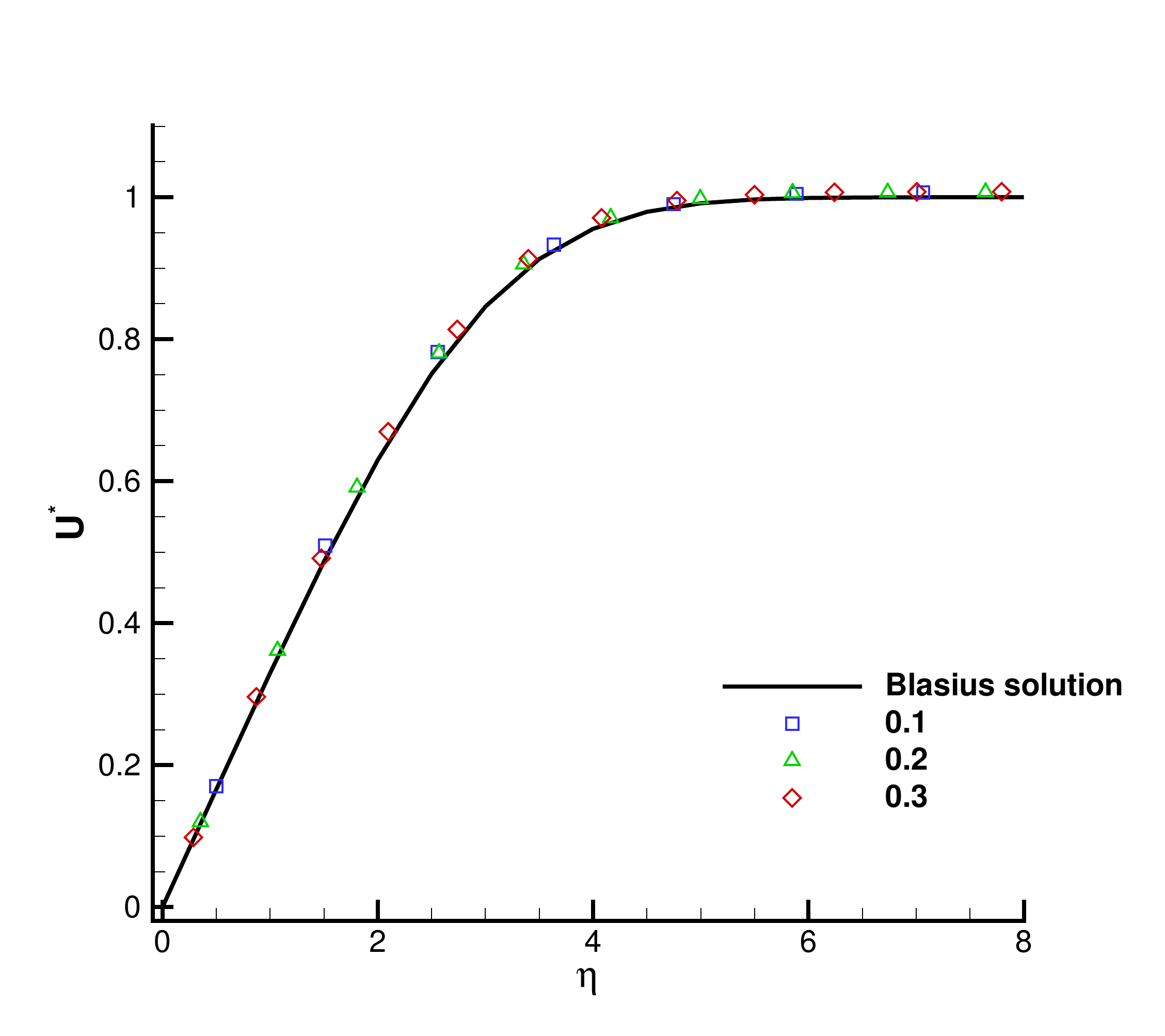}
	\includegraphics[width=0.3\textwidth]
	{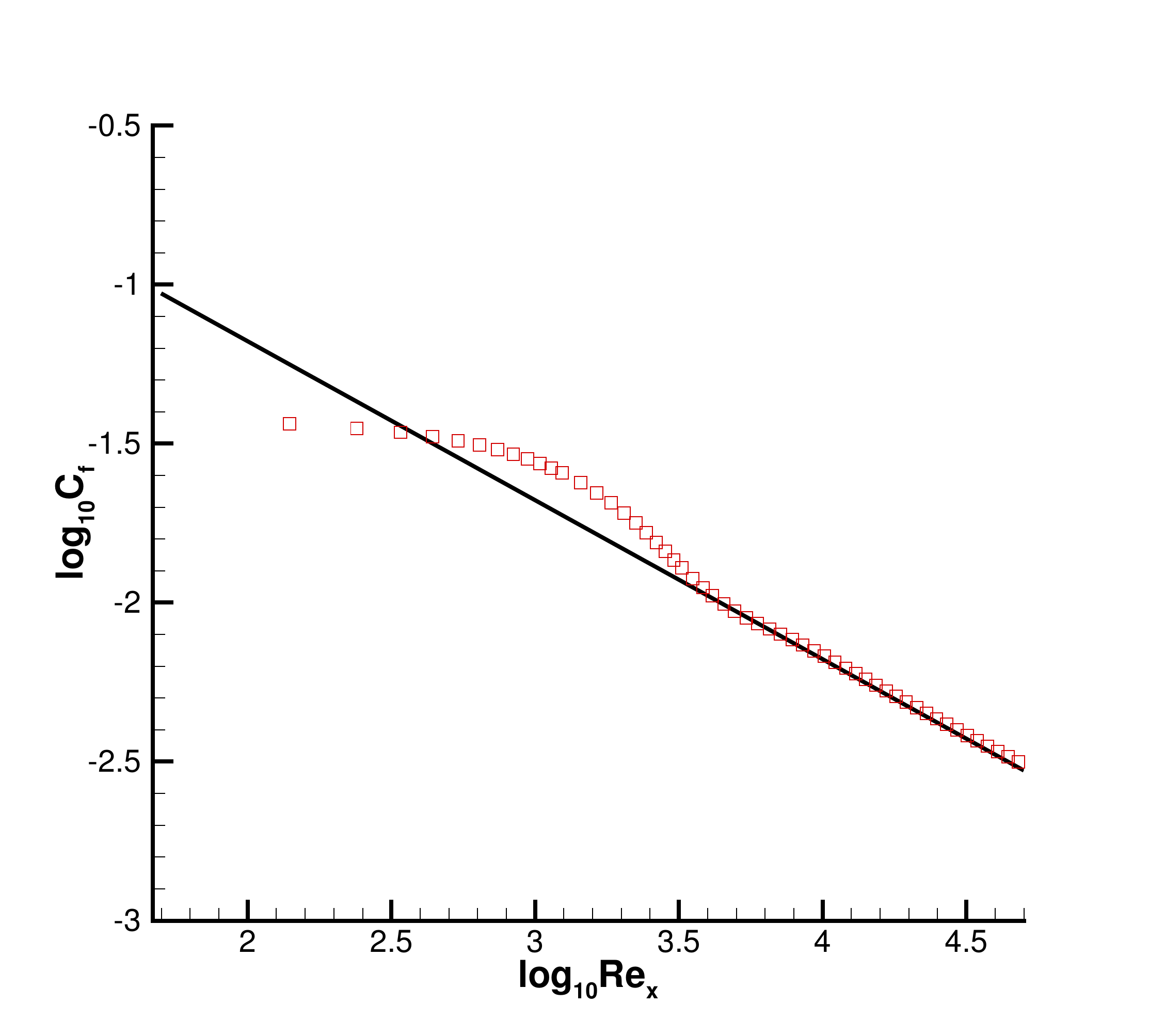}
	\vspace{-4mm} \caption{\label{flate-plate-line}
		Laminar boundary layer. The velocities profiles and skin-frication coefficient.}
\end{figure}

A laminar boundary layer over a flat plate with incoming Mach number $Ma=0.15$ is simulated. The Reynolds number $Re=U_{\infty}L/{\nu}=10^5$, where the characteristic length $L=100$.
The computational domain is shown in Fig.\ref{flate-plate-mesh}, where the flat plate is placed at $x>0$ and $y=0$.
Total $120 \times 35 \times 2$  hexahedral cells are used in a cuboid domain $[-30,100]\times[0,80]\times[0,0.2]$ with a cell height $h=0.1$ adjacent to the boundary.
The adiabatic non-slip boundary condition is imposed on the plate and symmetric slip boundary condition is set in the front of the plate.
The non-reflecting boundary condition based on the Riemann invariants is adopted for the other boundaries, where the free stream is set as $\rho_{\infty}=1, p_{\infty}=1/\gamma$.
Since the flow is nearly incompressible, the smooth reconstruction and the simplified solver in Eq.~\eqref{2nd-smooth-flux} are adopted to further reduce the numerical dissipation.
The non-dimensional velocity U and V are given in Fig.\ref{flate-plate-line} at three selected locations. The wall distributions of the skin-fraction coefficients $C_f$ are also plotted, where the local Reynolds number$Re_x$ and the $C_f$ are defined as
\begin{align*}
Re_x = \frac{x}{L}Re,~~~
C_f = \frac{\tau_{wall}}{\frac{1}{2}\rho_{\infty} U^2_{\infty}},
\end{align*}
where $\tau_{wall}$ is the skin shear stress.
The numerical results agree well with the Blasius solutions with a few mesh points at $x/L=0.1$.

\subsection{Subsonic viscous flow around a cylinder at Re=40}
\begin{figure}[htp]	
	\centering
	\includegraphics[width=0.4\textwidth]
	{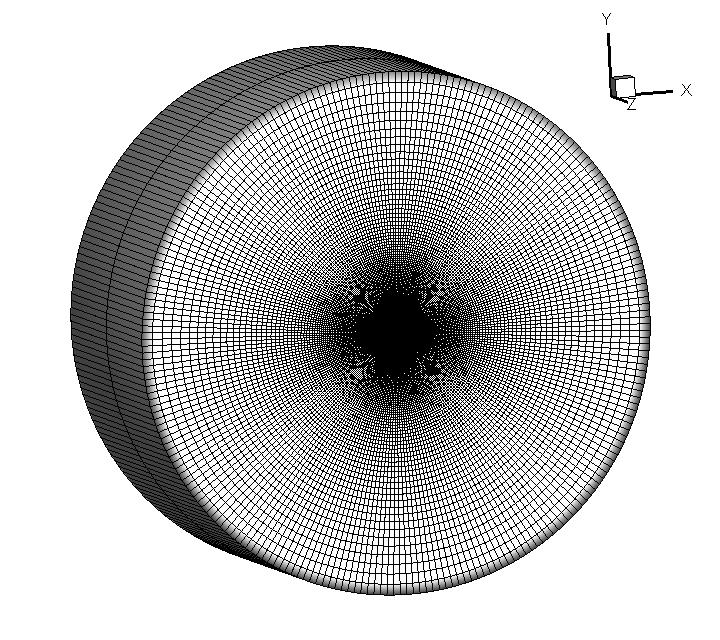}
	\includegraphics[width=0.4\textwidth]
	{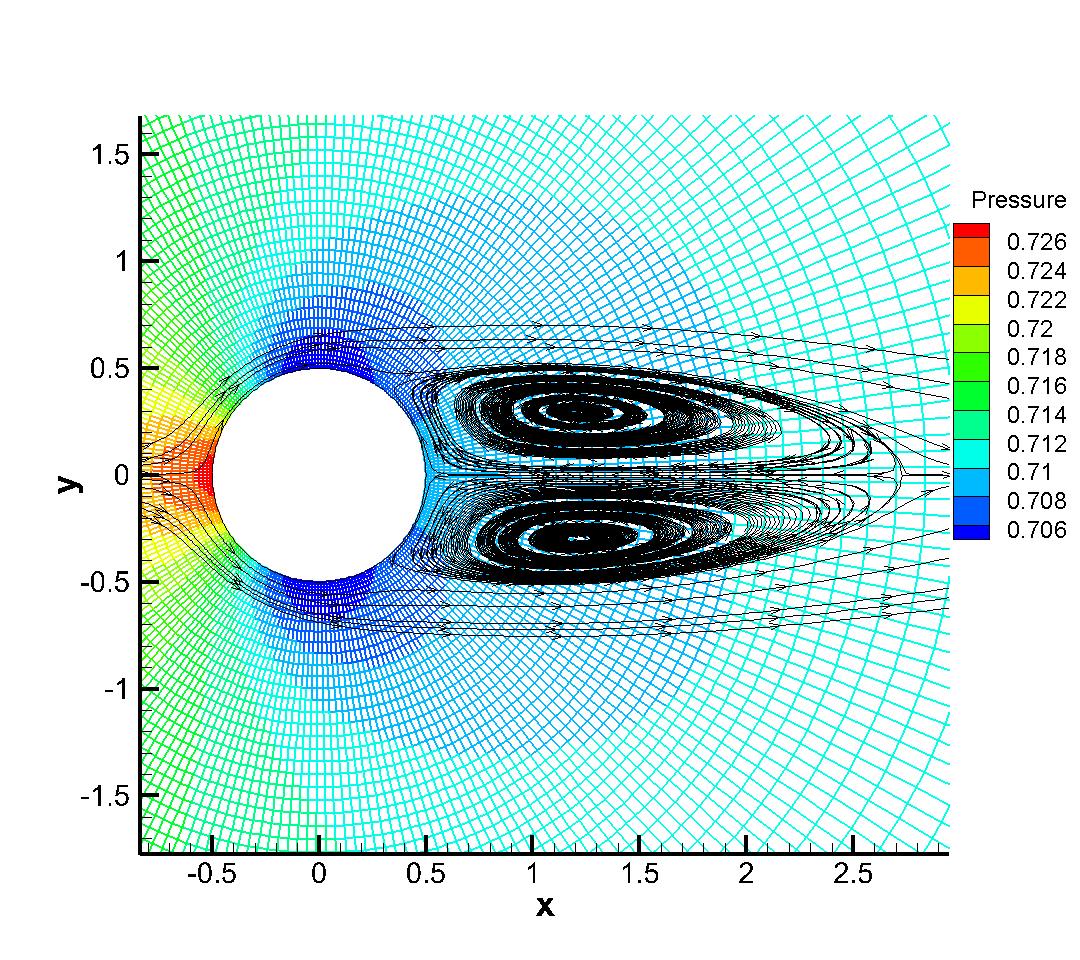}
	\vspace{-4mm} \caption{\label{cylinder-re40-mesh}
		Circular cylinder: Re=40. Left: mesh with $241\times 114\times2 $ cells.
		Right: local mesh distribution around cylinder colored by pressure and streamline. }
\end{figure}

\begin{figure}[htp]	
	\centering
	\includegraphics[width=0.4\textwidth]
	{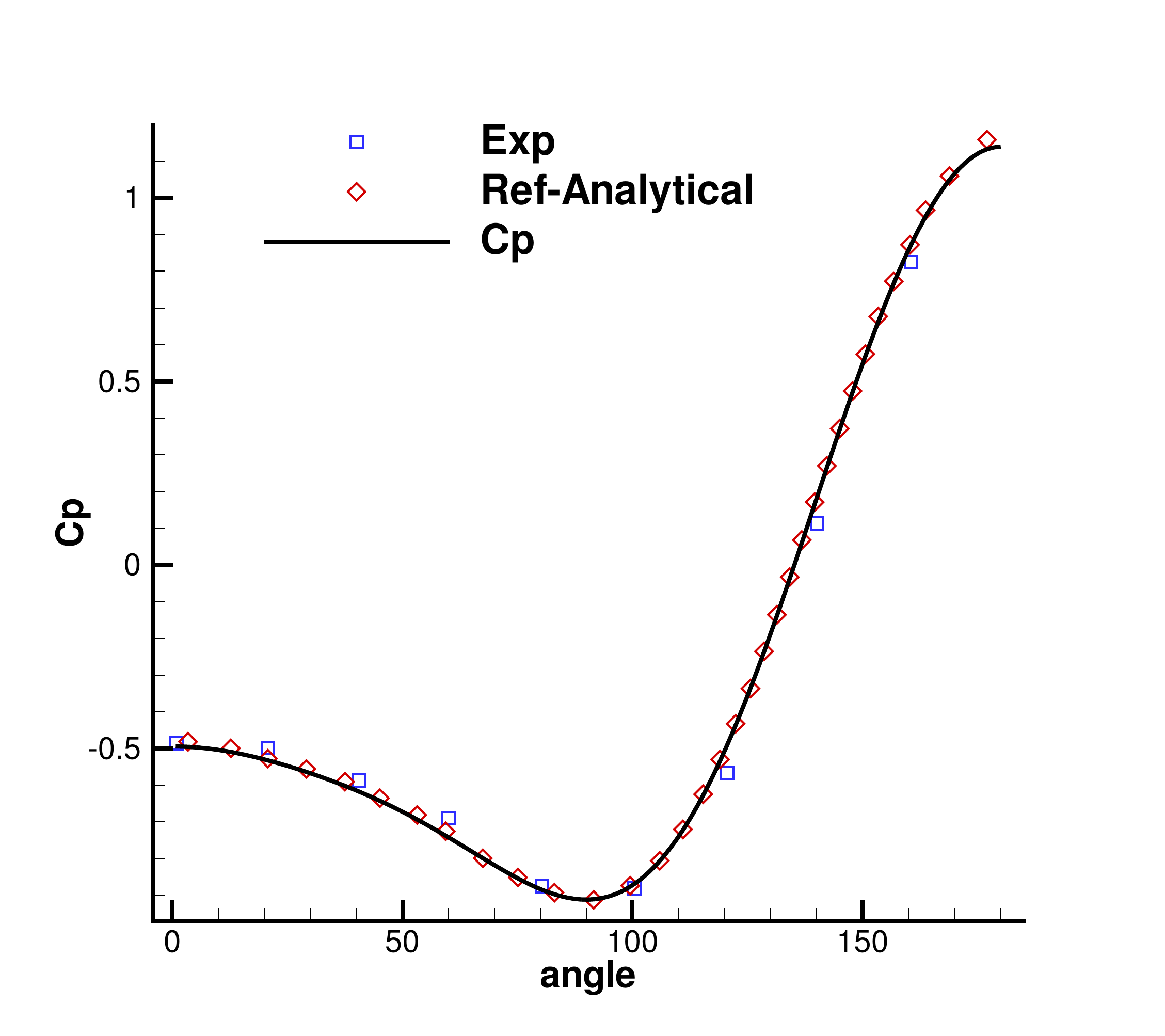}	
	\includegraphics[width=0.4\textwidth]
	{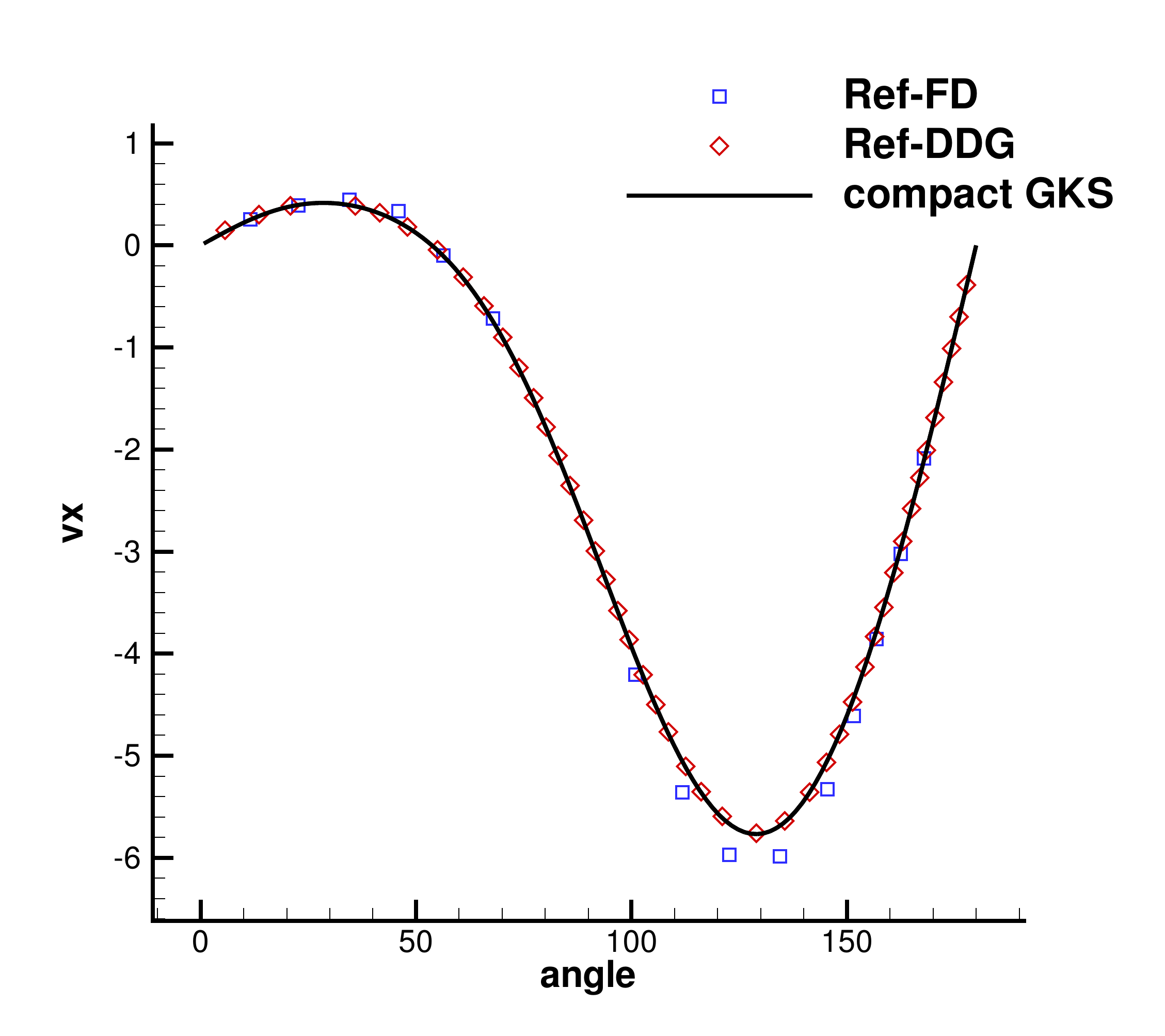}
	\vspace{-4mm} \caption{\label{cylinder-re40-line}
		Circular cylinder: Re=40. Left: surface pressure coefficient distribution. Right: surface local tangential velocity gradient  distributions.}
\end{figure}

An incoming flow with Mach number $Ma=0.15$ and Reynolds number $Re=40$ based on the diameter of the cylinder $D=1$ 
around a circular cylinder is simulated.
The flow is smooth and the same setting as the boundary layer is adopted in the computation.
The computational domain is shown in Fig.~\ref{cylinder-re40-mesh}, where total $241 \times 114 \times 2$  hexahedral cells are used in a cylindrical domain $D_{mesh}=96.0$, $H_{mesh}=0.1$ with a near wall size $h=1/96$.
A steady and symmetrical separation bubble is located at the wake of the cylinder.
Quantitative results including the drag and lift coefficients $C_D,~C_L$, the wake length $L$, and the separation angle $\theta$, etc are listed in Tab.~\ref{cylinder-re-40-cd-cl}, which agree well with the experimental and numerical references \cite{tritton1959experiments,coutanceau1977experimental,zhang2019direct}.
Furthermore, the quantities on the cylinder surface are extracted, including
the surface pressure coefficient $C_p=\frac{p-p_{\infty}}{\frac{1}{2}\rho_{\infty}U_{\infty}^2}$.
and the non-dimensional local tangential velocity gradient$ \frac{2U_{\infty}}{D} \frac{\partial U_{\tau}}{\partial \eta}$, as shown in Fig.~\ref{cylinder-re40-line}.
The $C_p$ from the current CGKS matches nicely with the experimental data \cite{grove1964experimental} and the analytical solution \cite{bharti2006steady}.
The tangential velocity gradient obtained by the current scheme is compared with those by the finite difference method \cite{braza1986Numerical} and the direct DG method \cite{zhang2019direct}.

\begin{table}[htp]
	\footnotesize
	\begin{center}
		\def\temptablewidth{1.0\textwidth}
		{\rule{\temptablewidth}{1pt}}
		\begin{tabular*}{\temptablewidth}{@{\extracolsep{\fill}}c|c|c|c|c|c|c}
			
			Case & $C_D$ &$C_L$ & $L$ & Vortex Height & Vortex Width & $\theta$ \\
			\hline
			Experiment \cite{tritton1959experiments} & 1.46 - 1.56 & -- & --   &  -- & -- & --\\
			Experiment \cite{coutanceau1977experimental} & -- & -- & 2.12   &  0.297 & 0.751 & 53.5$^\circ$\\	
			DDG\cite{zhang2019direct} & 1.529 & -- & 2.31   &  -- & -- & --\\					
			Current & 1.525 & 3.3e-14 & 2.22 & 0.296 & 0.714 & 53.3$^\circ$\\	
		\end{tabular*}
		{\rule{\temptablewidth}{0.1pt}}
	\end{center}
	\vspace{-4mm} \caption{\label{cylinder-re-40-cd-cl} Comparison of results for steady flow past a circular cylinder Re=40.}
\end{table}

\subsection{Flow passing through a sphere from subsonic to hypersonic cases}
Viscous flow over a sphere with a wide range of Mach numbers are tested below to validate the capability of the CGKS in different flow regions. The Reynolds number based on the diameter of the sphere $D=1$.
The far-field condition is set at outside boundary of the domain with the free stream condition
\begin{equation*}
\begin{split}
(\rho,U,V,W,p)_{\infty} =(1,Ma,0,0,\frac{1}{\gamma}),
\end{split}
\end{equation*}
with $\gamma=1.4$.
For the subsonic cases, the smooth reconstruction and the simplified solver in Eq.~\eqref{2nd-smooth-flux} are adopted to achieve a higher resolution.

\noindent{\sl{(a) Subsonic case: Re=118, Ma=0.2535.}}

\begin{figure}[htp]	
	\centering	
	\includegraphics[width=0.29\textwidth]{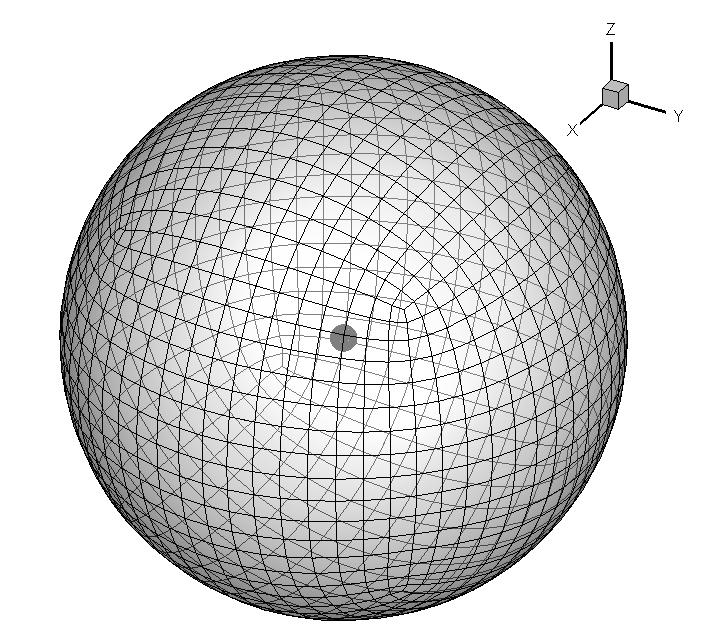}
	\includegraphics[width=0.29\textwidth]{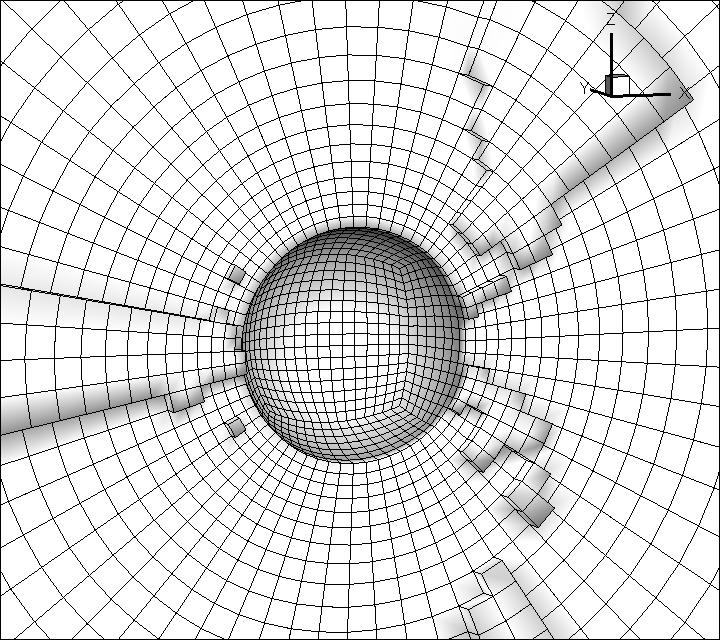}
	\includegraphics[width=0.34\textwidth]{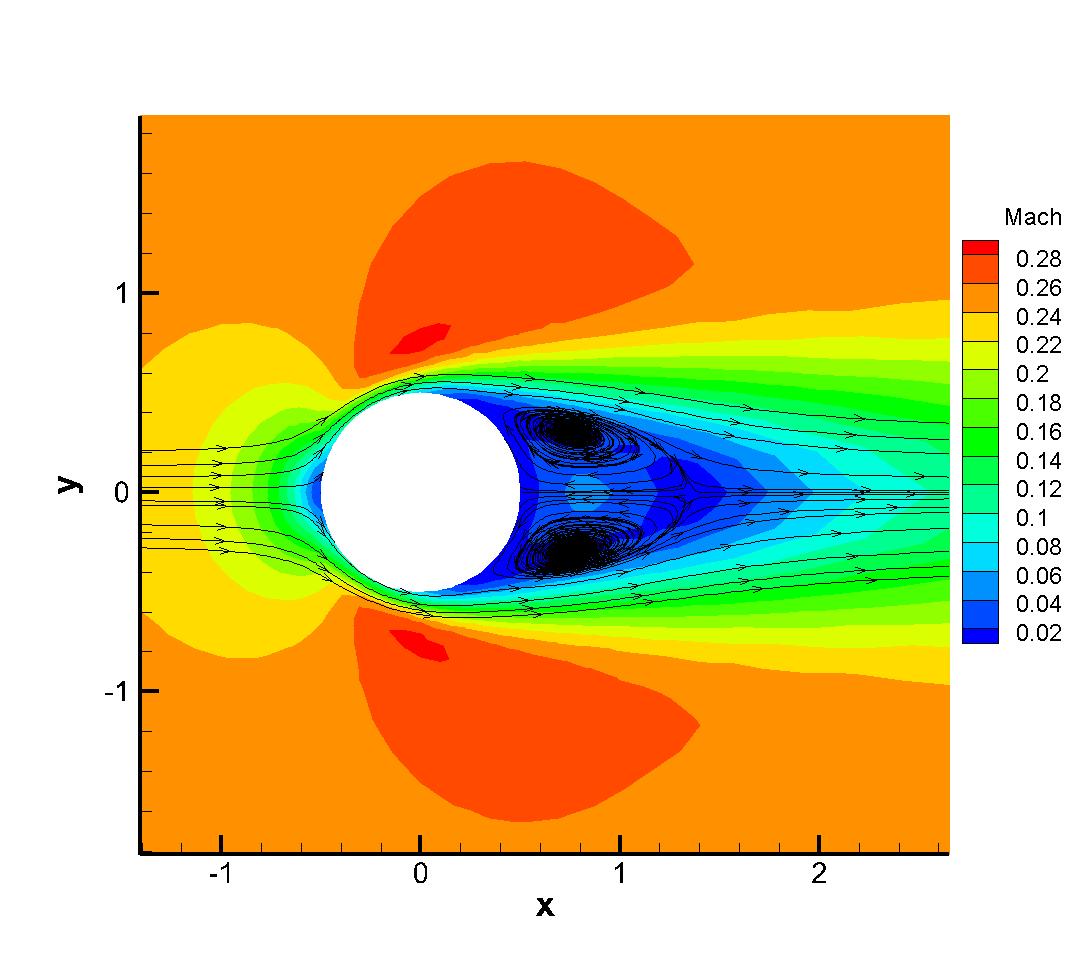}
	\caption{\label{sphere-re118-mesh}
		Flow passing through a sphere. Ma=0.2535. Re=118. Mesh number: 50,688.}
\end{figure}

A low-speed viscous flow passing through a sphere is tested first.
In such case, a drag coefficient $C_D=1$ was reported from the experiment in \cite{taneda1956experimental}.
The surface of the sphere is set as non-slip and adiabatic.
A hexahedral mesh is used and the first mesh off the wall has the size $ h \approx 4.5 \times 10^{-2} D$, as shown in Fig.~\ref{sphere-re118-mesh}.
The Mach magnitude contour and the streamline around the sphere are also given in Fig.~\ref{sphere-re118}.
The quantitative results are given in Table \ref{sphere-re118}, including the drag coefficient $C_D$,  the separation angle $\theta$, and the closed wake length $L$, as defined in \cite{ji2021compact}. The drag coefficient is very close to those by other methods even with a much coarser mesh.

\begin{table}[htp]
	\small
	\begin{center}
		\def\temptablewidth{1.0\textwidth}
		{\rule{\temptablewidth}{1pt}}
		\begin{tabular*}{\temptablewidth}{@{\extracolsep{\fill}}c|c|c|c|c|c}
			Scheme & Mesh number & $C_D$  & $\theta$  & $L$ &$C_L$\\
			\hline
			Experiment \cite{taneda1956experimental}	&-- & 1.0  & 151 & 1.07 & -- \\ 	
			Third-order DDG \cite{cheng2017parallel} & 160,868 & 1.016 & 123.7 & 0.96 & --\\
			Fourth-order VFV \cite{wang2017thesis}  & 458,915 & 1.014 & --& -- & 2.0e-5\\
			Current & 50,688 & 1.016  & 124.7 & 0.86 & 3.6e-3\\
		\end{tabular*}
		{\rule{\temptablewidth}{0.1pt}}
	\end{center}
	\vspace{-4mm} \caption{\label{sphere-re118} Quantitative comparisons among different compact schemes  for the viscous flow over a sphere.}
\end{table}

\noindent{\sl{(b) Subsonic case: Re=300, Ma=0.3.}}

\begin{figure}[htp]	
	\centering	
	\includegraphics[height=0.28\textwidth]{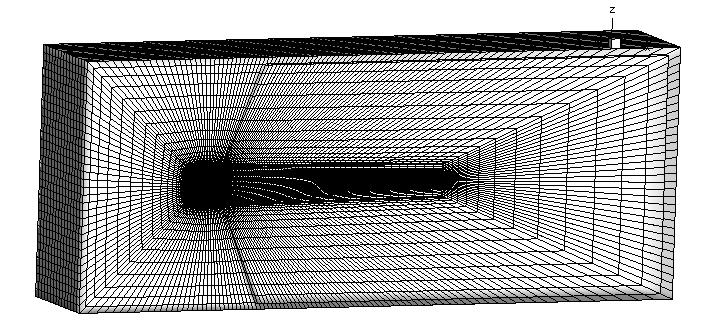}
	\includegraphics[height=0.28\textwidth]{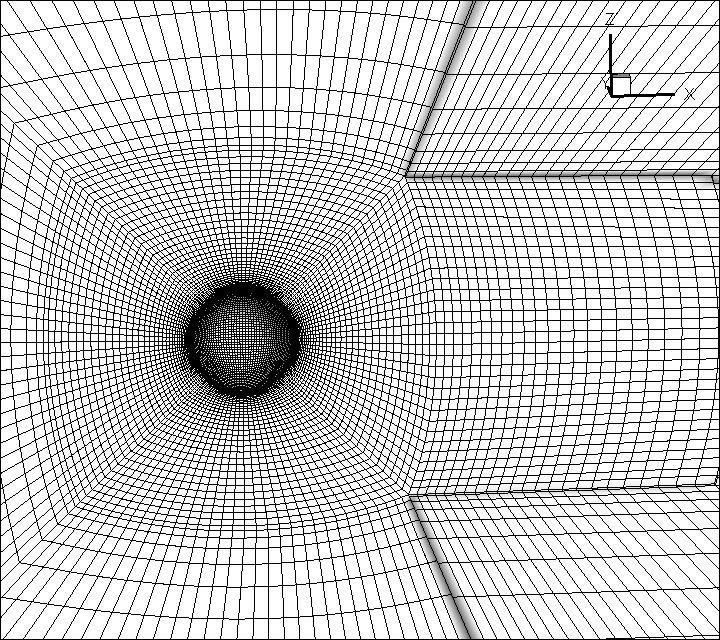}
	\caption{\label{sphere-ma0x3-re300-mesh}
		Flow passing through a sphere. Ma=0.3. Re=300. Mesh number: 479,232.}
\end{figure}

\begin{figure}[htp]	
	\centering	
	\includegraphics[height=0.28\textwidth]{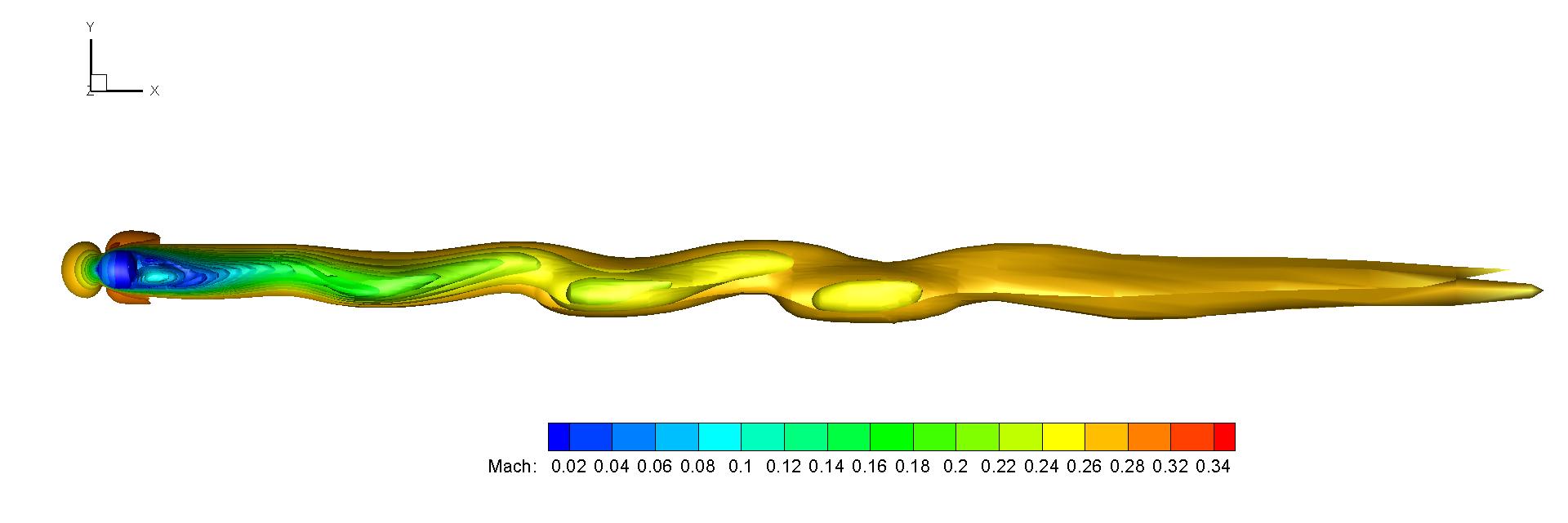}
	\includegraphics[height=0.28\textwidth]{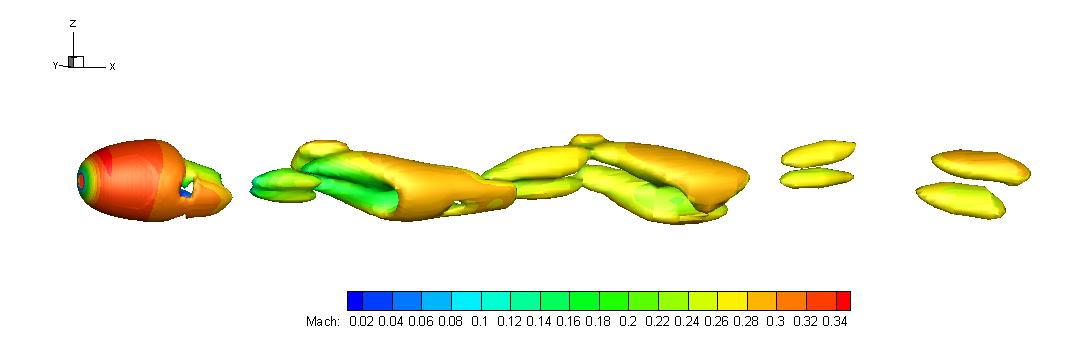}
	\includegraphics[height=0.28\textwidth]{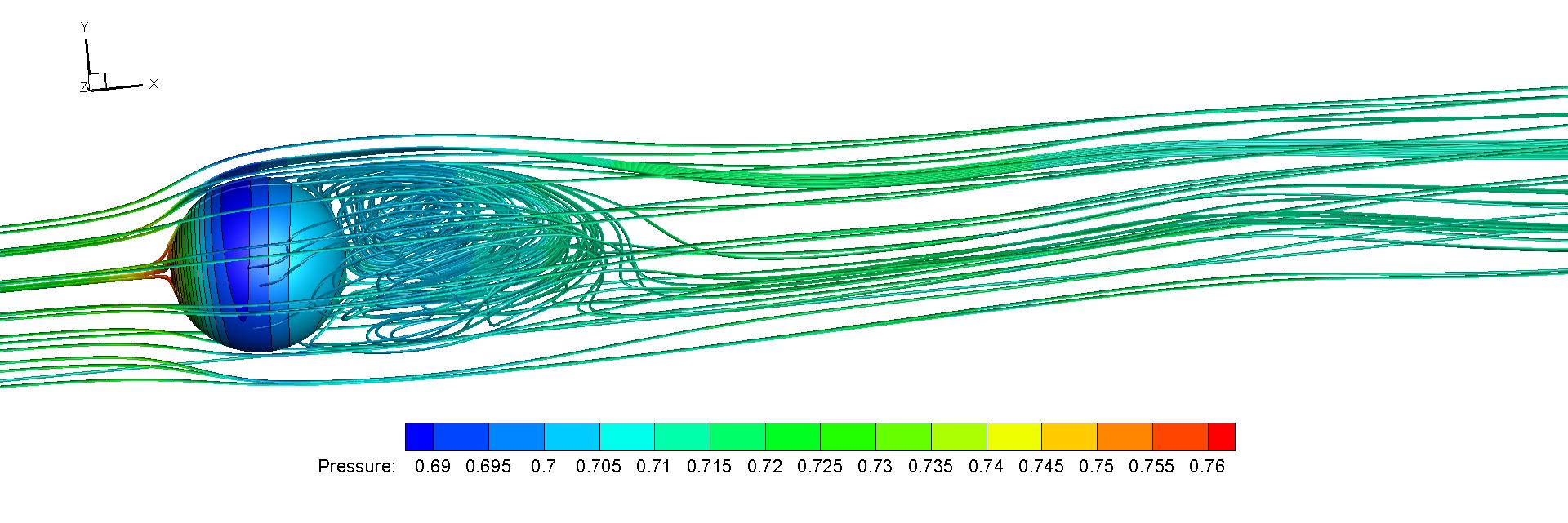}
	\caption{\label{sphere-ma0x3-re300-contour}
		Flow passing through a sphere. Ma=0.3. Re=300. Top: Iso-surface of the Mach number. Middle: Iso-surface of the Q criterion $Q=5\times10^{-4}$ colored by Mach number. Bottom: 3-D streamline colored by pressure.}
\end{figure}

\begin{figure}[htp]	
	\centering	
	\includegraphics[height=0.28\textwidth]{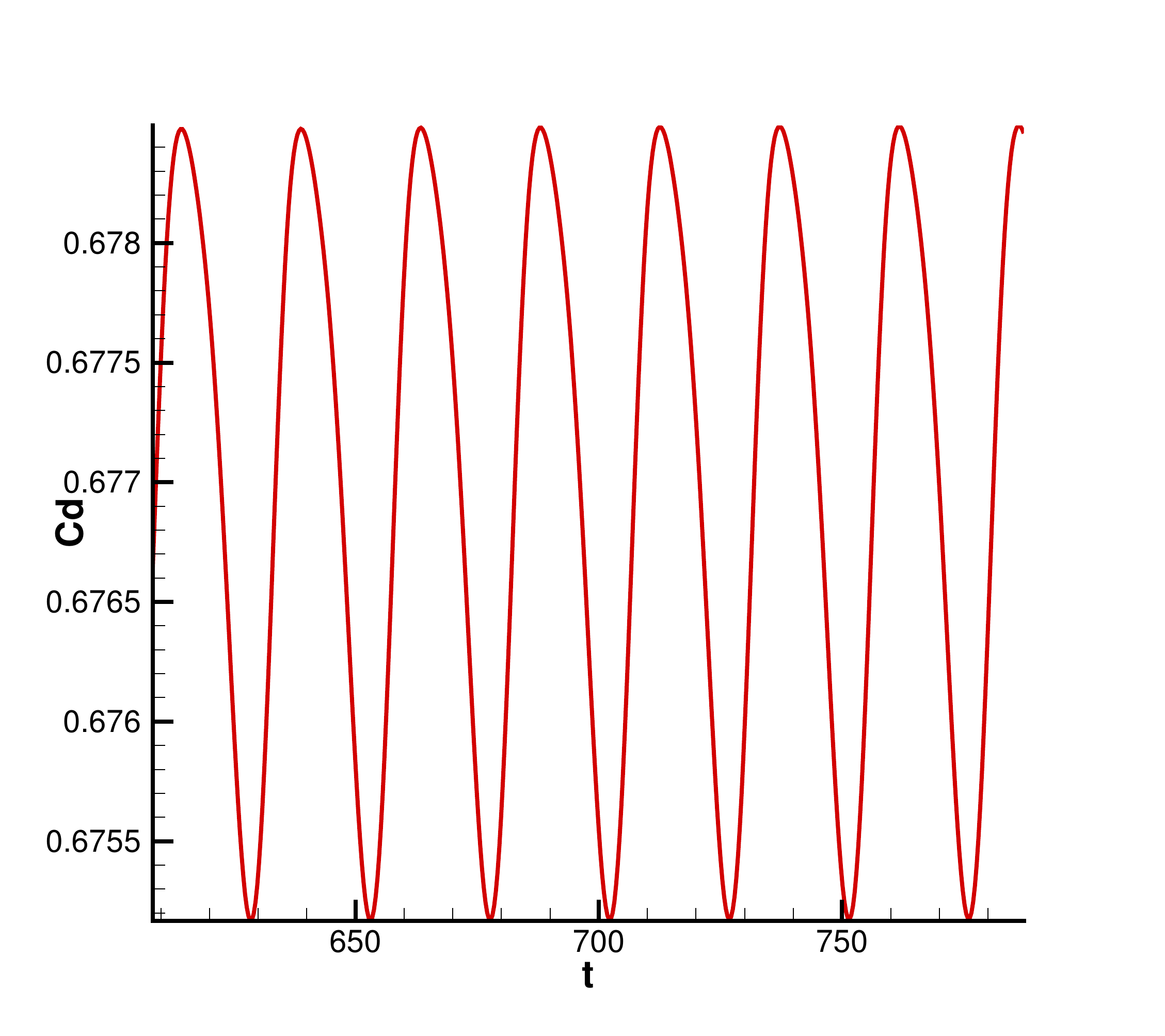}
	\includegraphics[height=0.28\textwidth]{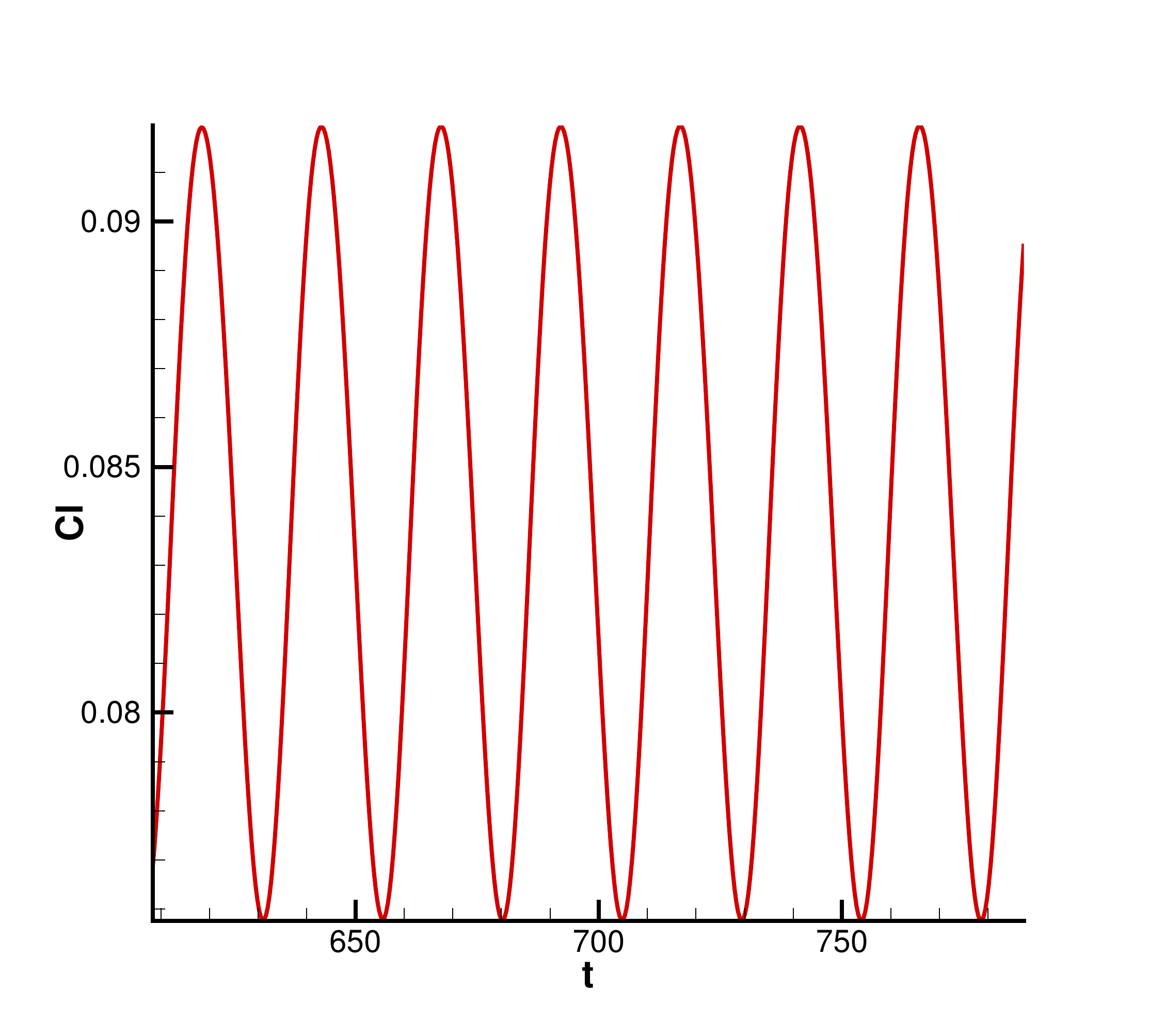}
	\caption{\label{sphere-ma0x3-re300-cd}
		The time history of the $C_D$ and $C_L$. Ma=0.3. Re=300.}
\end{figure}

The flow is unsteady in this case and the hairpin vortex structure will be formed in the wake region of the sphere.
A hexahedral mesh is used and the whole computation domain has a dimension $[-10,40]\times[-10,10]\times[-10,10]$ with a near wall size  $h \approx \frac{1}{100} D$ along the radial direction and $h \approx \frac{1}{128} D$ along the circumferential direction,
as shown in Fig.~\ref{sphere-ma0x3-re300-mesh}.
The current scheme can resolve the vortex shedding nicely,
as shown in Fig.~\ref{sphere-ma0x3-re300-contour}.
The drag and lift coefficients change in a single frequency mode with time increasing, as shown in Fig.~\ref{sphere-ma0x3-re300-cd}.
The frequency and averaged drag coefficient agree well with the reference data, as listed in Tab.~\ref{sphere-ma0x3-re300}.

\begin{table}[htbp]
	\small
	\begin{center}
		\def\temptablewidth{1.0\textwidth}
		{\rule{\temptablewidth}{1pt}}
		\begin{tabular*}{\temptablewidth}{@{\extracolsep{\fill}}c|c|c|c|c|c|c}
			Scheme & Mesh Number & $C_D$  & $\Delta C_D$  & $|C_L|$  & $\Delta C_L$ & St\\
			\hline
			Third-order k-exact \cite{li2014efficient} &2,065,612 & 0.674  & 0.003 & 0.055 & 0.013&0.133 \\ 	
			Current & 479,232 & 0.677  & 0.003 & 0.084 & 0.016 & 0.135 \\	
		\end{tabular*}
		{\rule{\temptablewidth}{0.1pt}}
	\end{center}
	\vspace{-4mm} \caption{\label{sphere-ma0x3-re300} Quantitative comparisons between the current scheme and the reference solution for the supersonic viscous flow over a sphere.}
\end{table}

\noindent{\sl{(c) Transonic case: Re=300, Ma=0.95.}}

\begin{figure}[htp]	
	\centering	
	\includegraphics[height=0.28\textwidth]{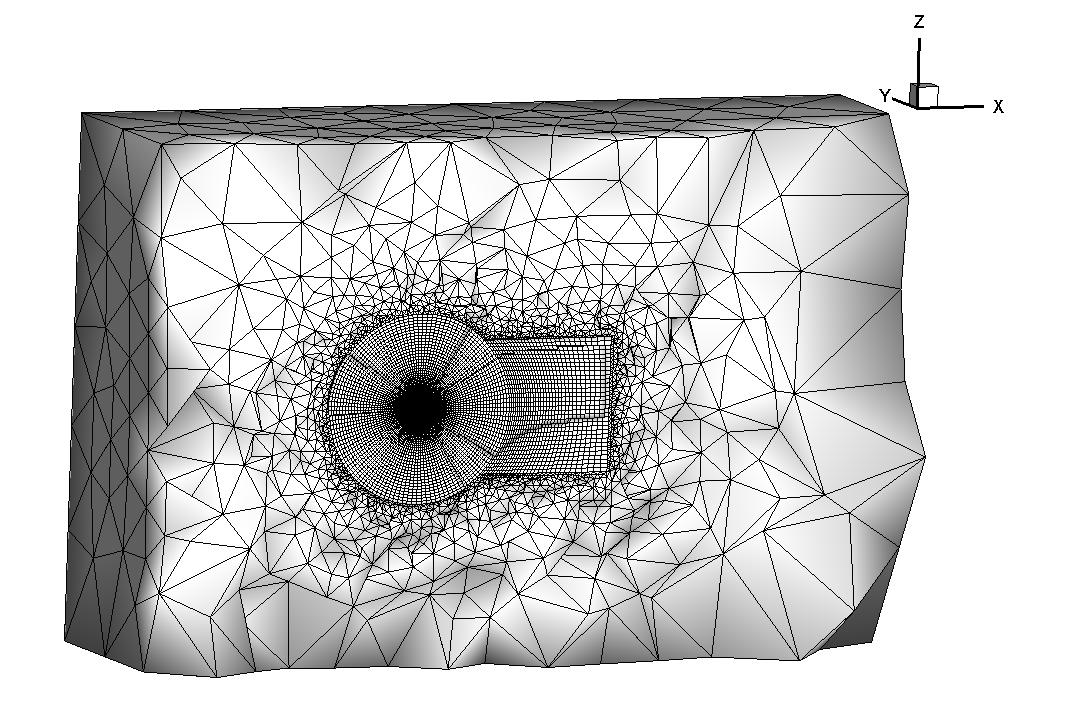}
	\includegraphics[height=0.28\textwidth]{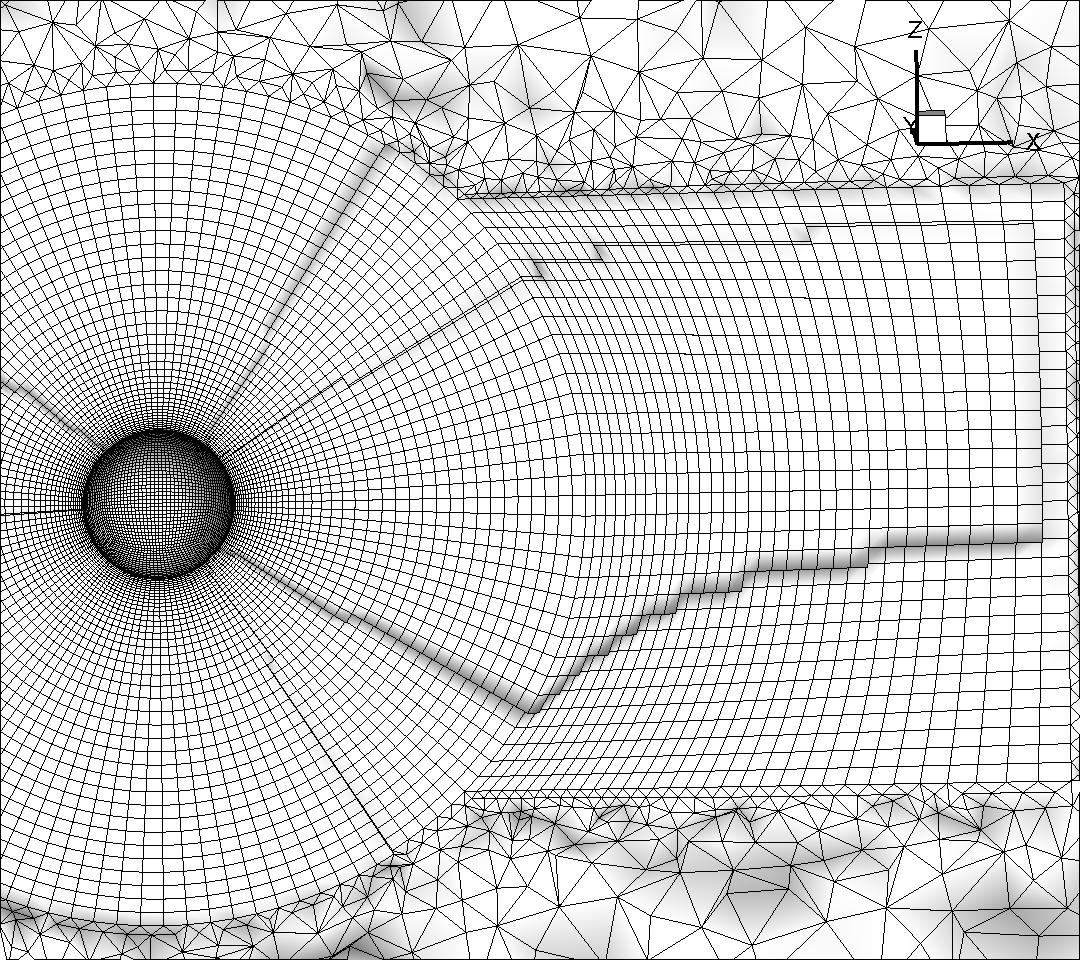}
	\caption{\label{sphere-ma0x95-re300-mesh}
		Flow passing through a sphere. Ma=0.95. Re=300. Mesh number: 515,453.}
\end{figure}

\begin{figure}[htp]	
	\centering	
	\includegraphics[height=0.28\textwidth]{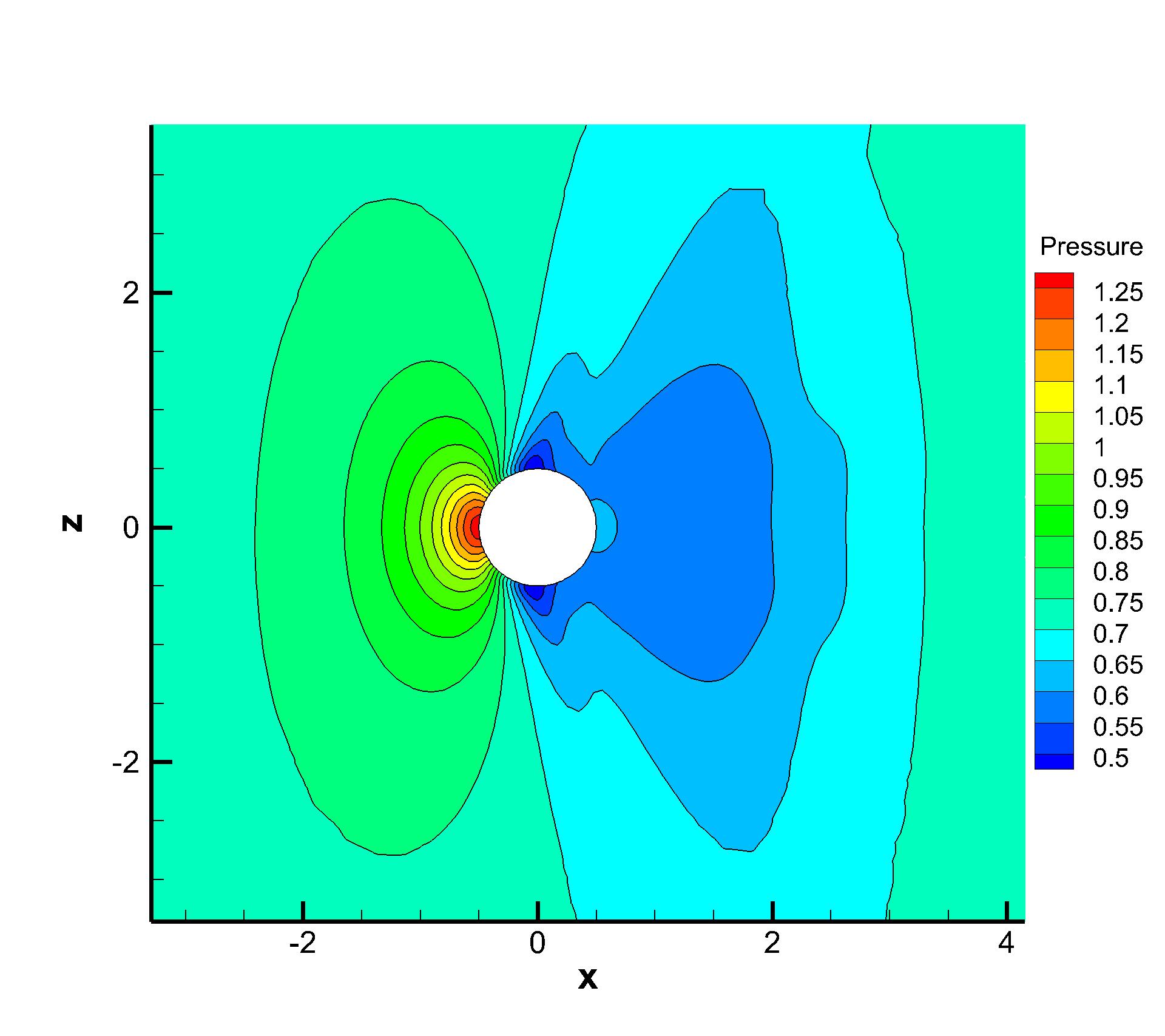}
	\includegraphics[height=0.28\textwidth]{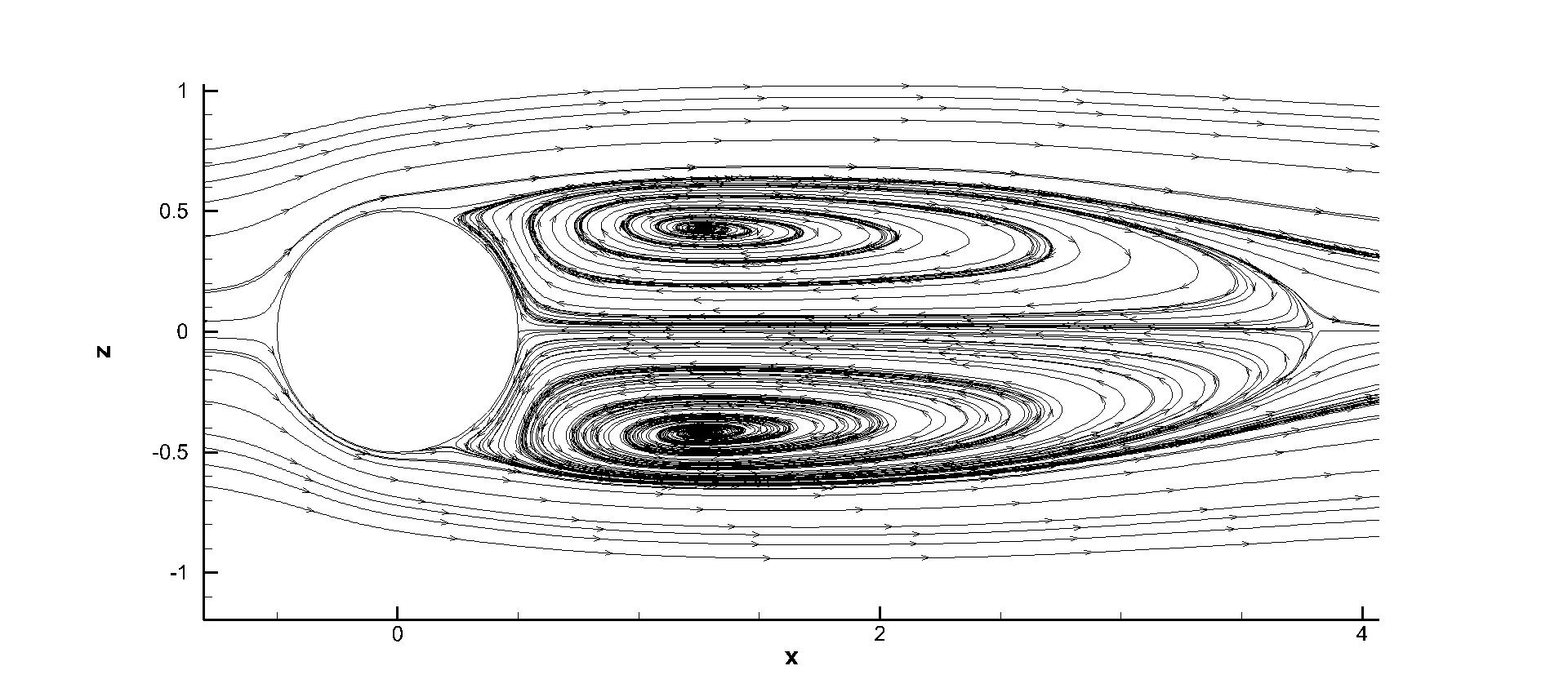}
	\caption{\label{sphere-ma0x95-re300-contour}
		Flow passing through a sphere. Ma=0.95. Re=300.}
\end{figure}

A hybrid unstructured mesh with total  $515,453$ cells is used in the computation, as shown in Fig.~\ref{sphere-ma0x95-re300-mesh}.
The first grid off the wall is $1 \times 10^{-2} D$ while 128 cells are distributed along the circumferential direction.
The non-reflective boundary is adopted on the outside boundary with a dimension $[-8,16]\times [-8,8]\times [-8,8]$.
The mesh is refined at the wake of the sphere where a long separation bubble is formed, as shown in Fig.~\ref{sphere-ma0x95-re300-contour}.
The pressure distribution is also given in Fig.~\ref{sphere-ma0x95-re300-contour}, where the weak shock can be observed.
Quantitative results are compared with the benchmark solutions in
\cite{Nagata2016sphere}, where the $C_d$ and $\theta$ agree well with each other, as shown in Tab.~\ref{sphere-ma0x95-re300}.

\begin{table}[htp]
	\small
	\begin{center}
		\def\temptablewidth{0.8\textwidth}
		{\rule{\temptablewidth}{1pt}}
		\begin{tabular*}{\temptablewidth}{@{\extracolsep{\fill}}c|c|c|c|c}
			Scheme & Mesh Number & Cd  & $\theta$  &L \\
			\hline
			WENO6 \cite{Nagata2016sphere} 	&909,072 & 0.968  & 111.5 & 3.48  \\ 	
			Current & 515,453 & 0.950  & 112.7 & 3.30 \\	
		\end{tabular*}
		{\rule{\temptablewidth}{0.1pt}}
	\end{center}
	\vspace{-4mm} \caption{\label{sphere-ma0x95-re300} Quantitative comparisons between the CGKS and the benchmark solution \cite{Nagata2016sphere} for the transonic flow over a sphere. Ma=0.95. Re=300.}
\end{table}

\noindent{\sl{(d) Supersonic case: Re=300, Ma=2.0.}}

\begin{figure}[htp]	
	\centering
	\includegraphics[height=0.4\textwidth]
	{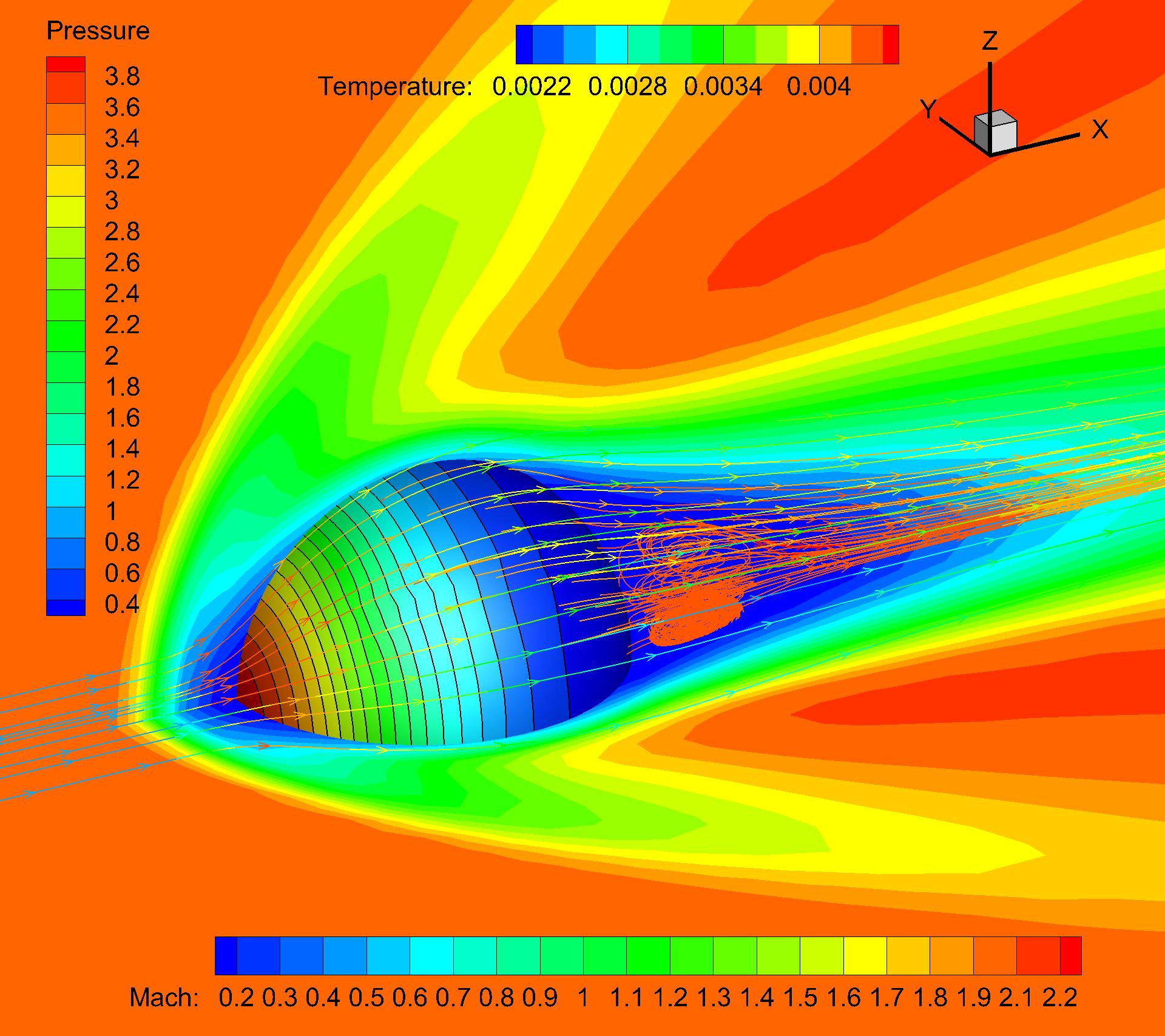}
	\includegraphics[height=0.4\textwidth]
    {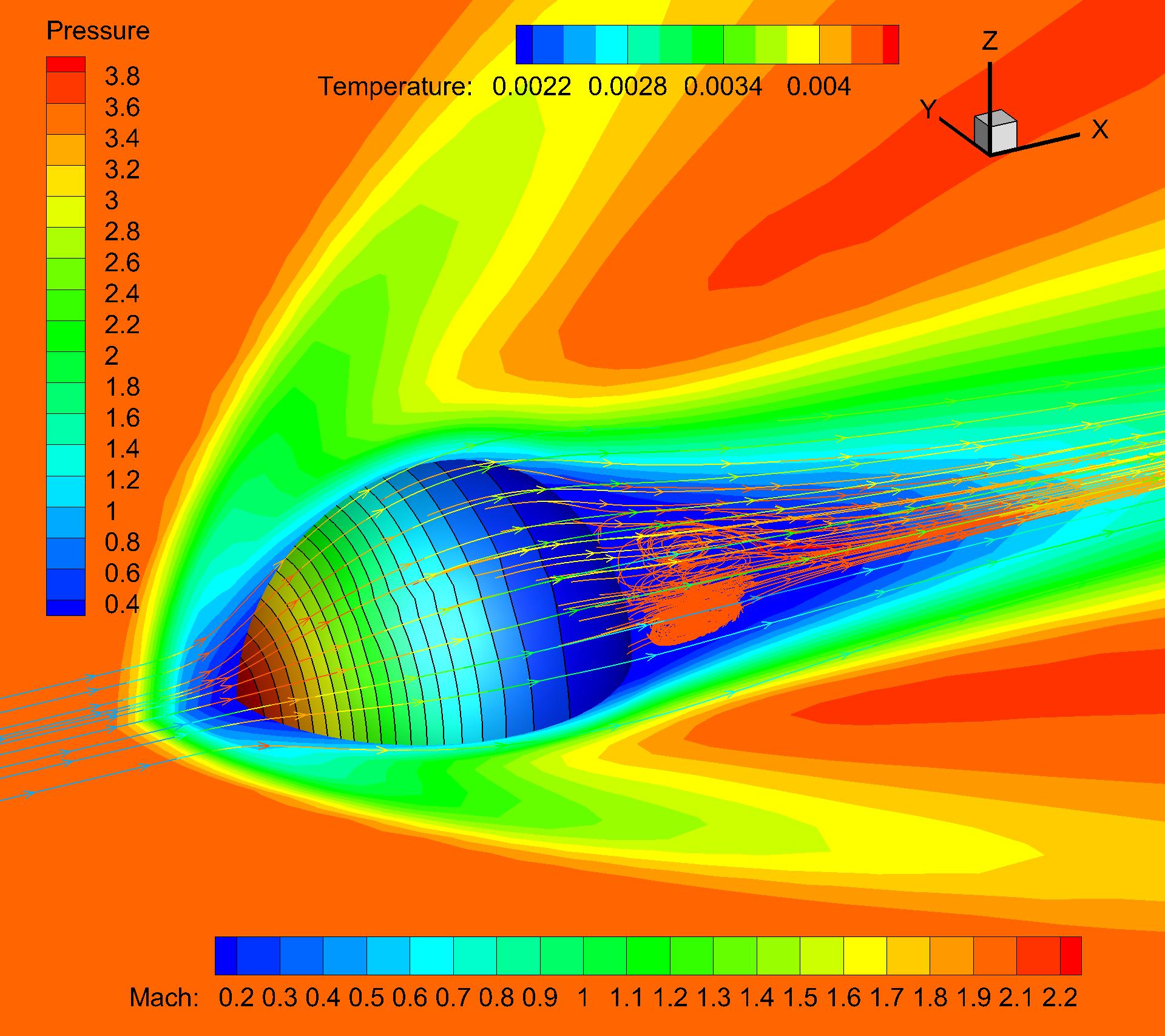}
	\caption{\label{sphere-ma2-ma}
		Flow passing through a sphere. Ma=1.2. Re=300. Left: without the CF. Right: with the CF.}
\end{figure}

\begin{figure}[htp]	
	\centering
	\includegraphics[height=0.2\textwidth]
	{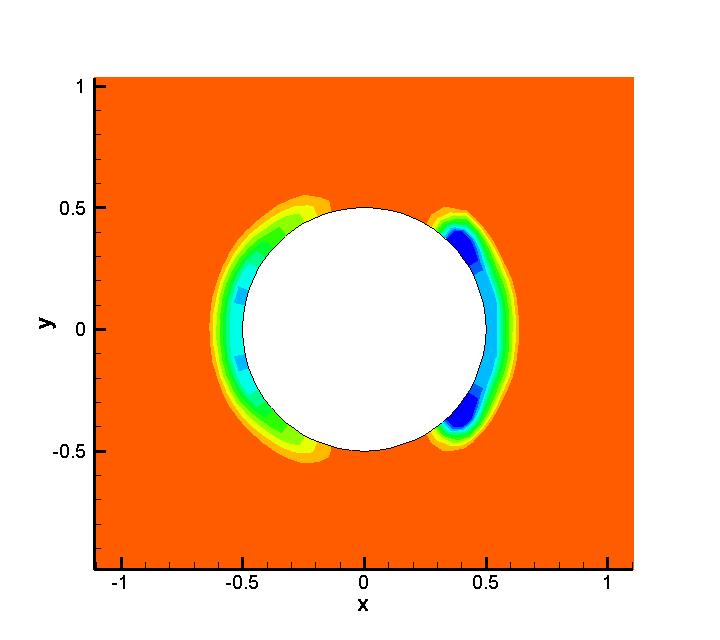}
			\includegraphics[height=0.2\textwidth]
	{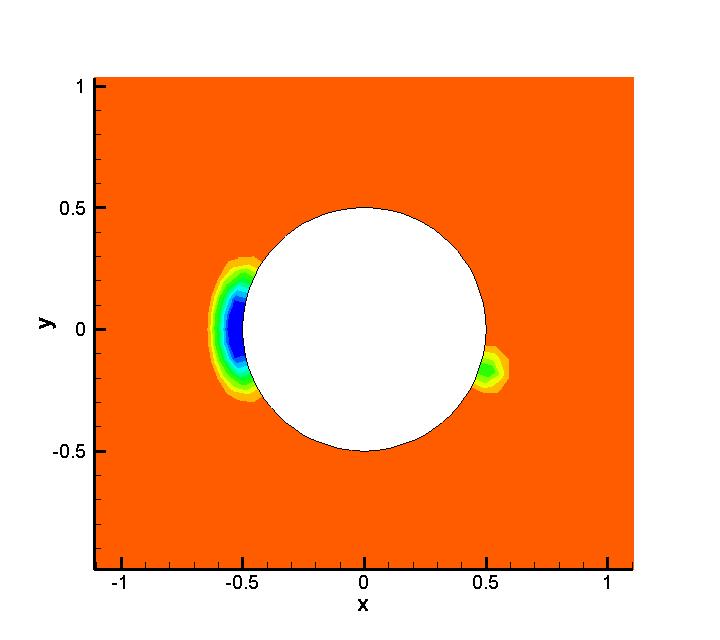}
		\includegraphics[height=0.2\textwidth]
	{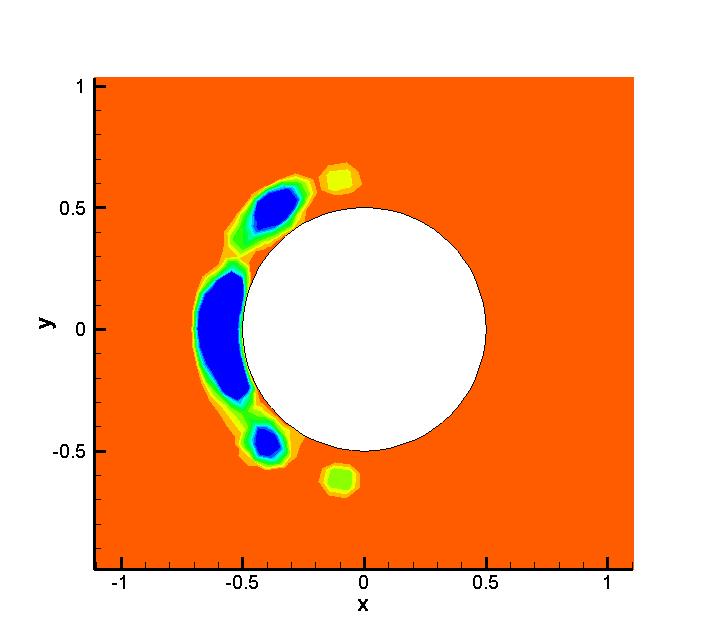}
			\includegraphics[height=0.2\textwidth]
	{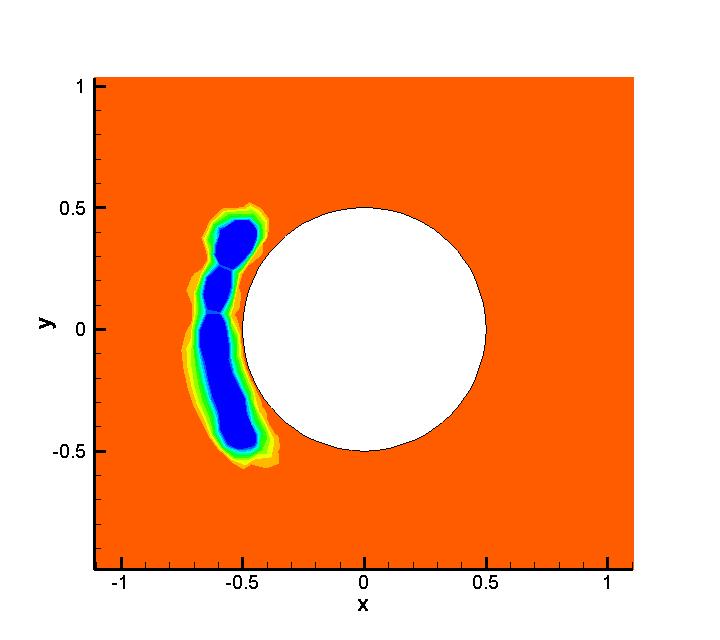}
	\caption{\label{sphere-ma2-cf-1}
	The CF distributions at step 2, 10, 20, 30 (from left to right).}
\end{figure}

\begin{figure}[htp]	
	\centering
	\includegraphics[height=0.2\textwidth]
	{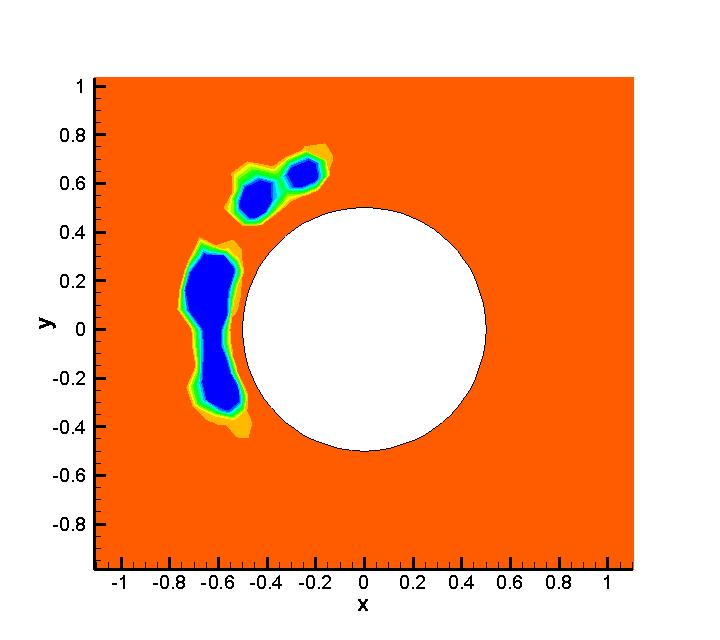}
	\includegraphics[height=0.2\textwidth]
	{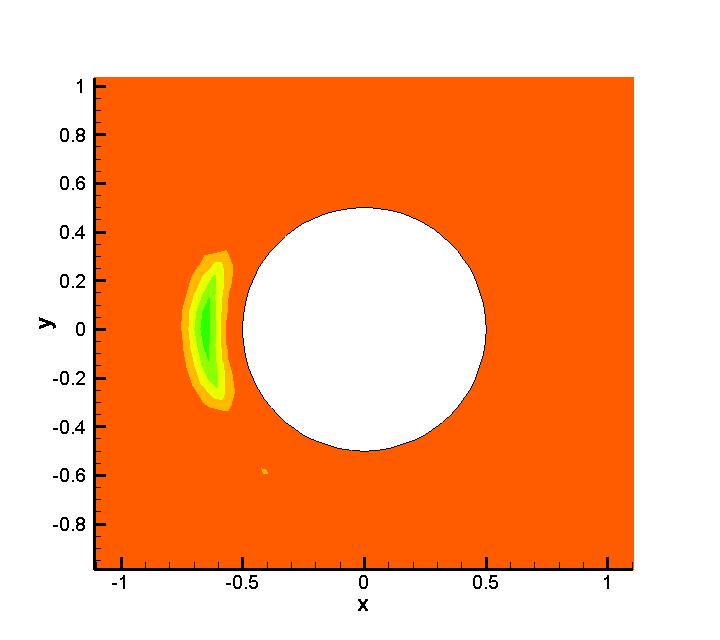}
	\includegraphics[height=0.2\textwidth]
	{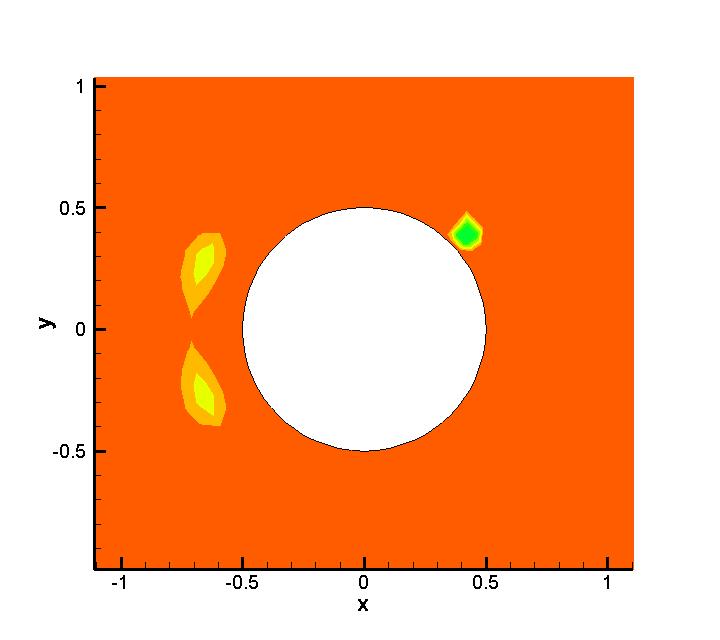}
	\includegraphics[height=0.2\textwidth]
	{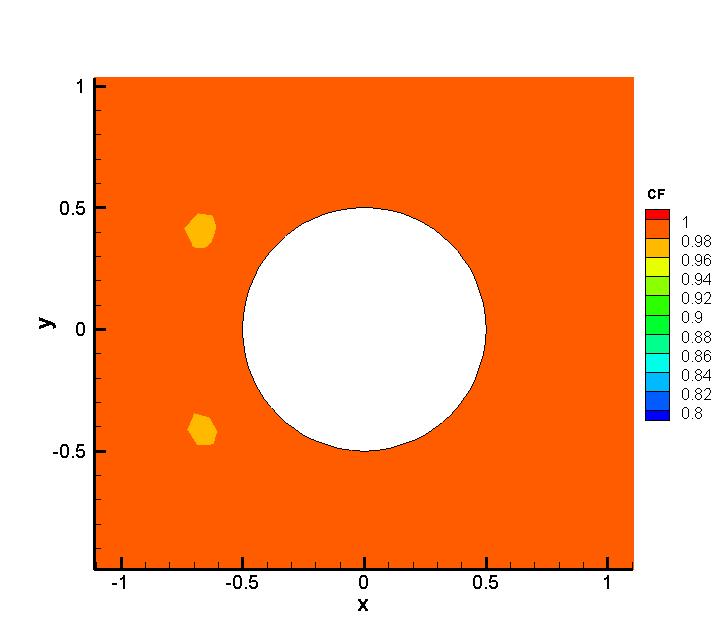}
	\caption{\label{sphere-ma2-cf-2}
		The CF distributions at step 40, 50, 100, 1000 (from left to right).}
\end{figure}

To evaluate the effect of the CF for the supersonic flow, a viscous flow around a sphere with $Ma=2.0$ is tested.
The non-slip adiabatic boundary condition is imposed on the surface of the sphere.
The same computational mesh for case (a) is used here.
To pass this case, an initial field calculated by the first-order kinetic method \cite{xu2014direct} has to be used for the CGKS without the CF.
The numerical results obtained by the CGKS with/without the CF are shown in Fig.~\ref{sphere-ma2-ma}.
Almost identical contours and streamlines are obtained.
Quantitative comparisons are listed in Tab.~\ref{sphere-ma2}.
Very close results are obtained by the CGKS with/without the CF.
And they have good agreements with those given by Nagata et al. \cite{Nagata2016sphere}.
Based on these observations, it is suggested that the current CF can improve the robustness of the CGKS while keeping the same level of accuracy in supersonic region.
The CF distributions are different time steps are given in Fig.~\ref{sphere-ma2-cf-1} and Fig.~\ref{sphere-ma2-cf-2}.
An interesting dynamical change of the CF can be observed.
At the very early step (Step 2), the slopes are modified near the tail of the sphere, where a rarefaction wave is formed instantly.
Then, the region for a small CF gradually moves to the front of the sphere, where the bow shock formed.
Finally, there is almost no region with the CF less than 0.98 when approaching to the steady state.

\begin{table}[htp]
	\small
	\begin{center}
		\def\temptablewidth{1.0\textwidth}
		{\rule{\temptablewidth}{1pt}}
		\begin{tabular*}{\temptablewidth}{@{\extracolsep{\fill}}c|c|c|c|c|c}
			Scheme & Mesh Number & Cd  & $\theta$  &L & Shock stand-off\\
			\hline
			WENO6 \cite{Nagata2016sphere} 	&909,072 & 1.386  & 150.9 & 0.38 & 0.21 \\ 	
			CGKS without CF & 50,688 & 1.368  & 148.5 & 0.45 & 0.28-0.31 \\
			CGKS with CF & 50,688 & 1.368  & 149.2 & 0.45 & 0.28-0.31 \\	
		\end{tabular*}
		{\rule{\temptablewidth}{0.1pt}}
	\end{center}
	\vspace{-4mm} \caption{\label{sphere-ma2} Quantitative comparisons between the current scheme and the reference solution for the flow over a sphere with $Ma=2$, $Re=300$.}
\end{table}

\noindent{\sl{(e) Hypersonic case: Re=300, Ma=5.0.}}

\begin{figure}[htp]	
	\centering
	\includegraphics[width=0.32\textwidth]
	{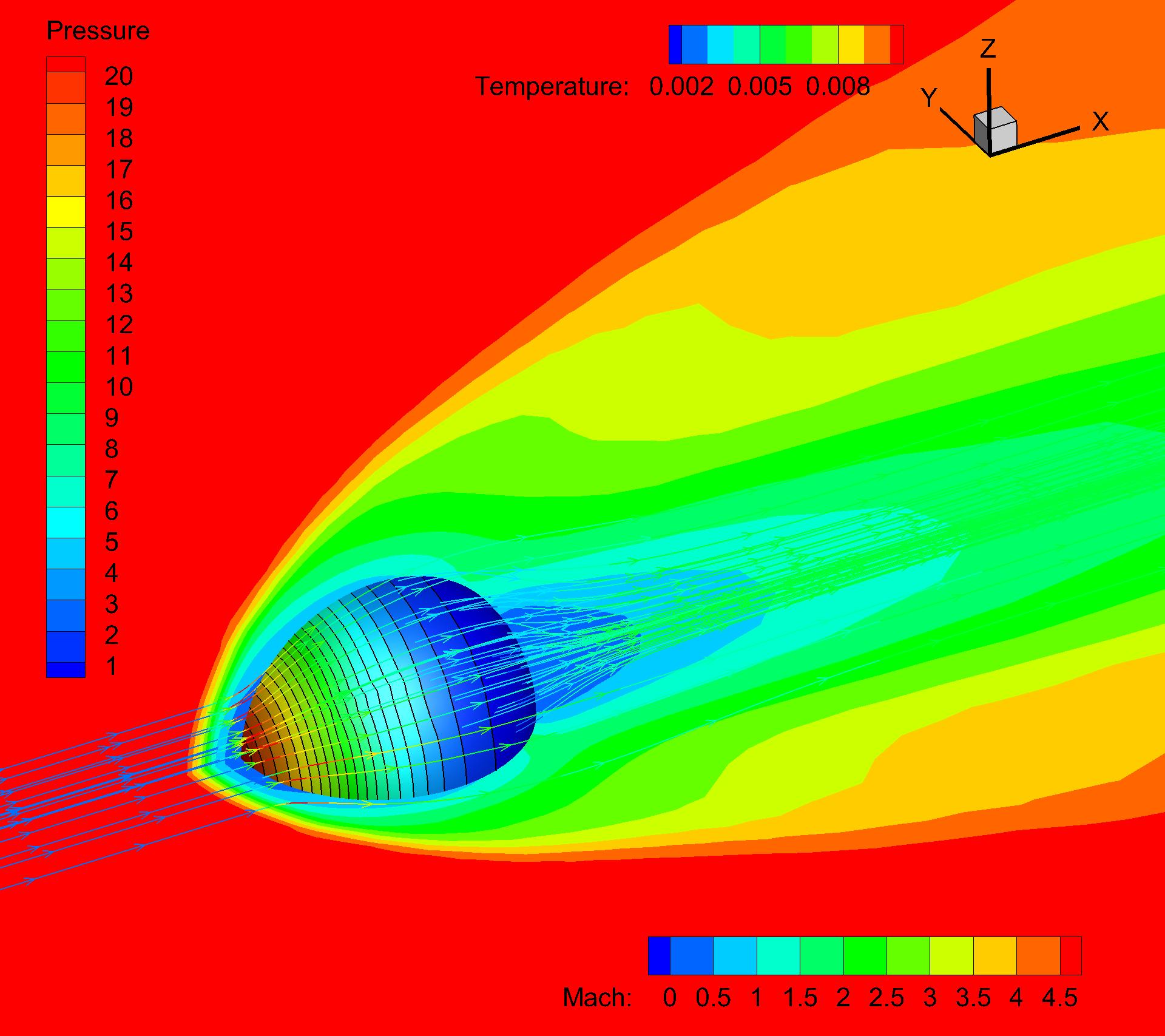}
	\includegraphics[width=0.32\textwidth]
	{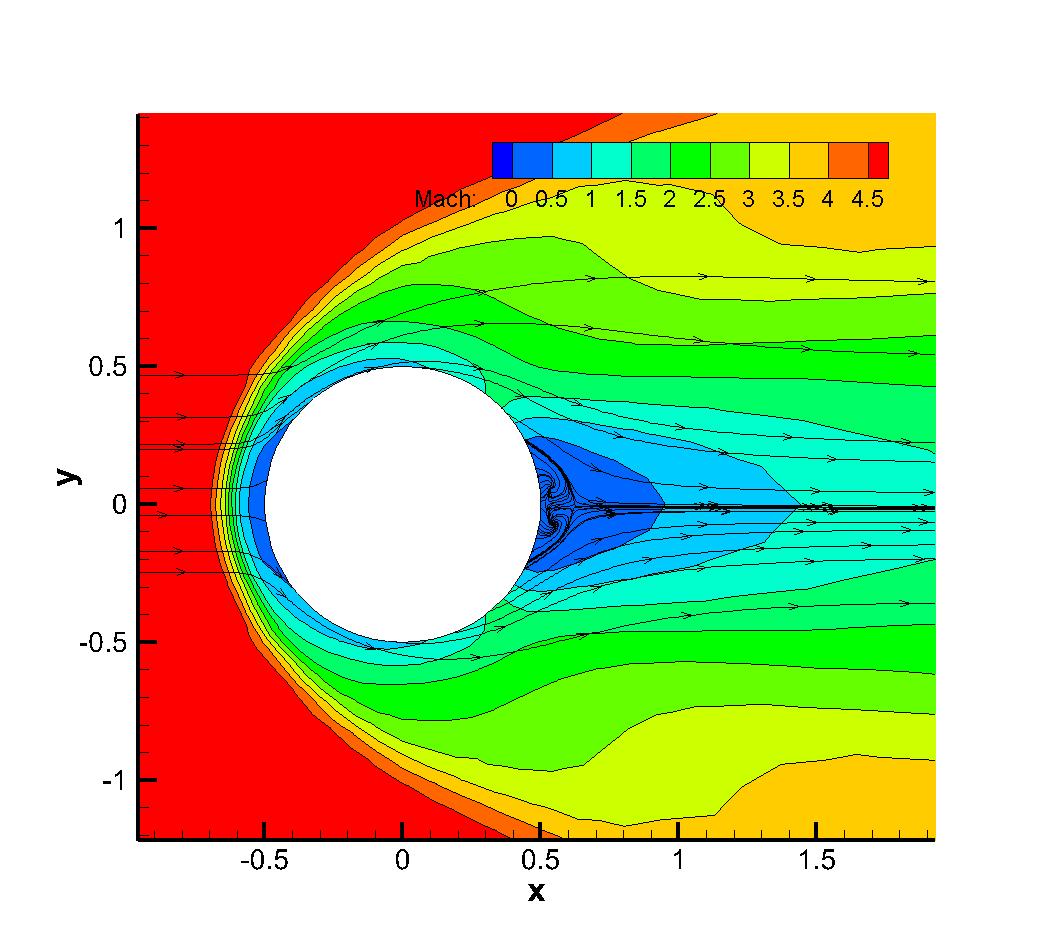}
		\includegraphics[width=0.32\textwidth]
	{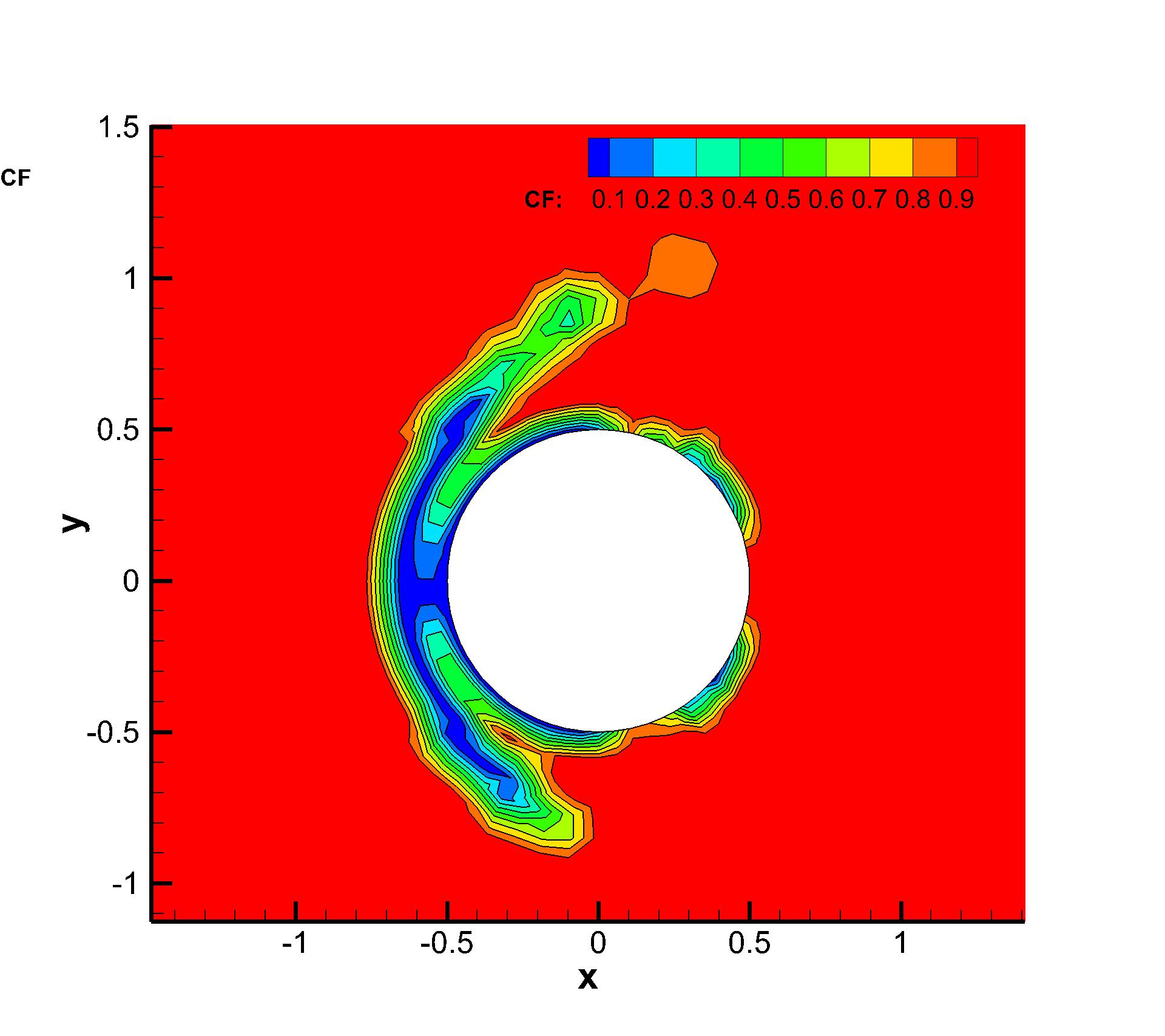}
	\caption{\label{sphere-ma5}
		Flow passing through a sphere Ma=5. Re=300. Left: the surface of the sphere is colored by pressure and the streamline is colored by temperature. Right: the CF distribution.}
\end{figure}

The $Ma=5$ flow passing through a sphere is considered here.
The Reynolds number is still set as $Re=300$. The Maxwell isothermal wall is adopted on the surface of the sphere.
The same computational mesh for case (a) is used here.
A primary flow field  calculated by the first-order kinetic method \cite{xu2014direct} is used as the initial field.
With the CF, the CGKS can safely pass this stringent test.
The Mach contours and streamlines in both 2-D and 3-D views are presented in
Fig.~\ref{sphere-ma5}. Only a tiny recirculation region is formed.  In addition, through the CF distribution in Fig.~\ref{sphere-ma5}, it can be observed that the cell-averaged slopes are mostly modified near the wall and the bow shock. The minimum CF is $3\times10^{-29}$.

\subsection{Transonic inviscid ONERA M6 wing}
\begin{figure}[htp]	
	\centering	
	\includegraphics[width=0.32\textwidth]{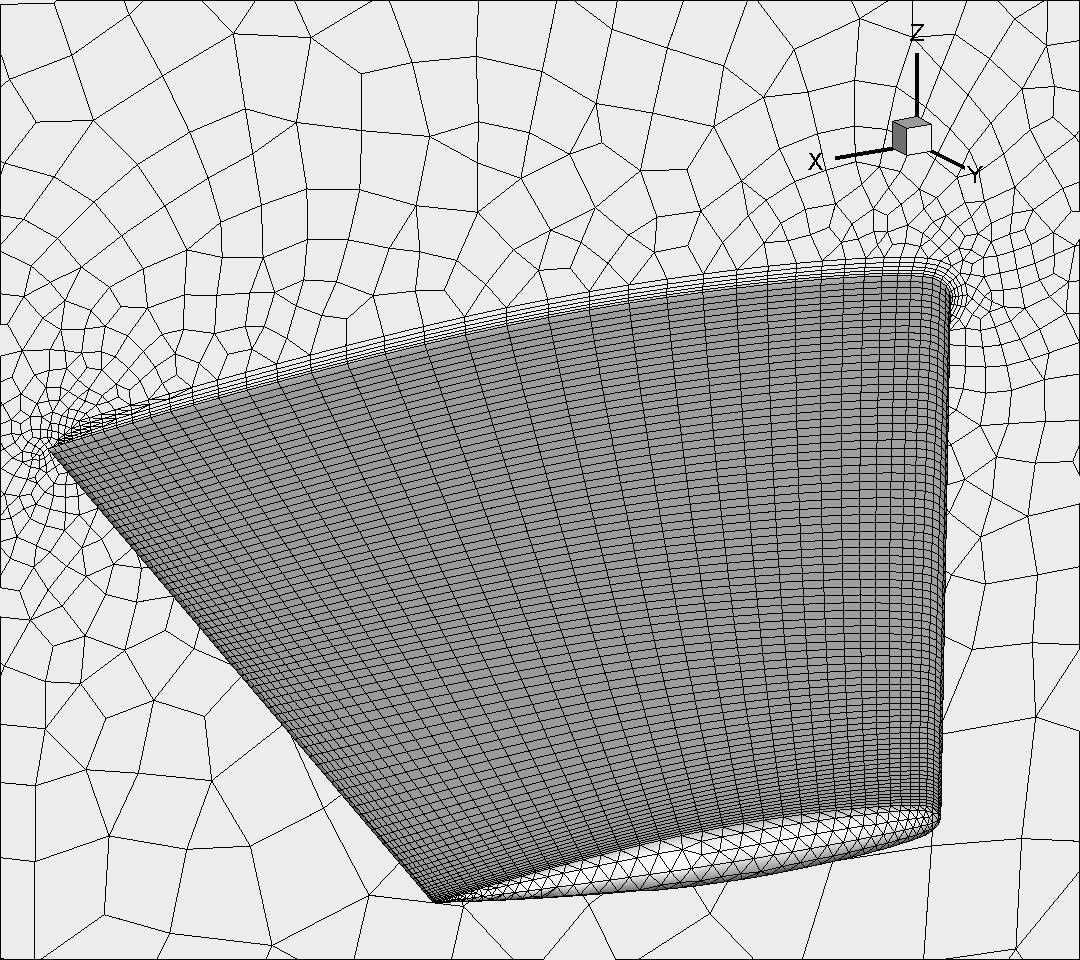}
	\includegraphics[width=0.32\textwidth]{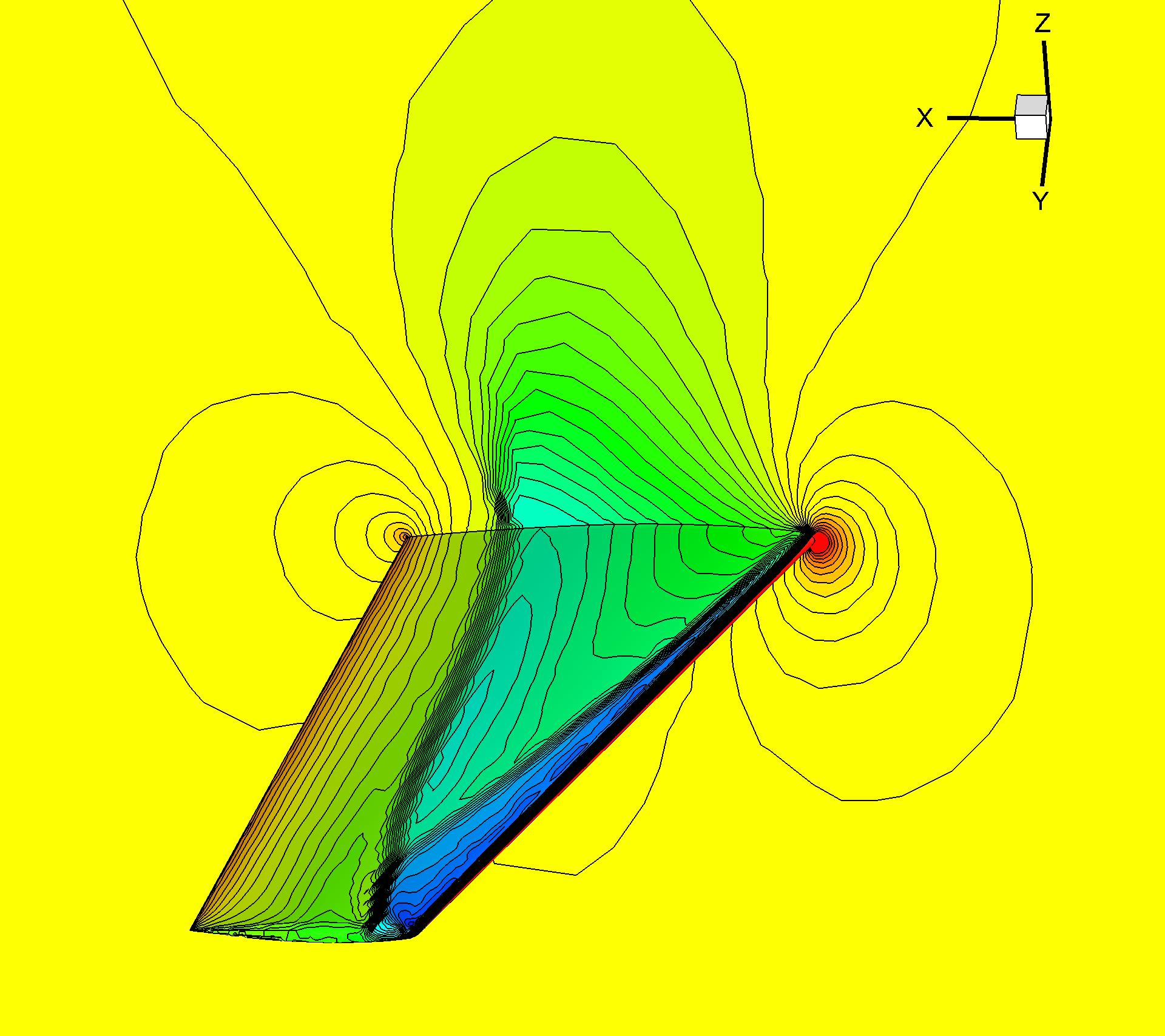}
	\includegraphics[width=0.32\textwidth]{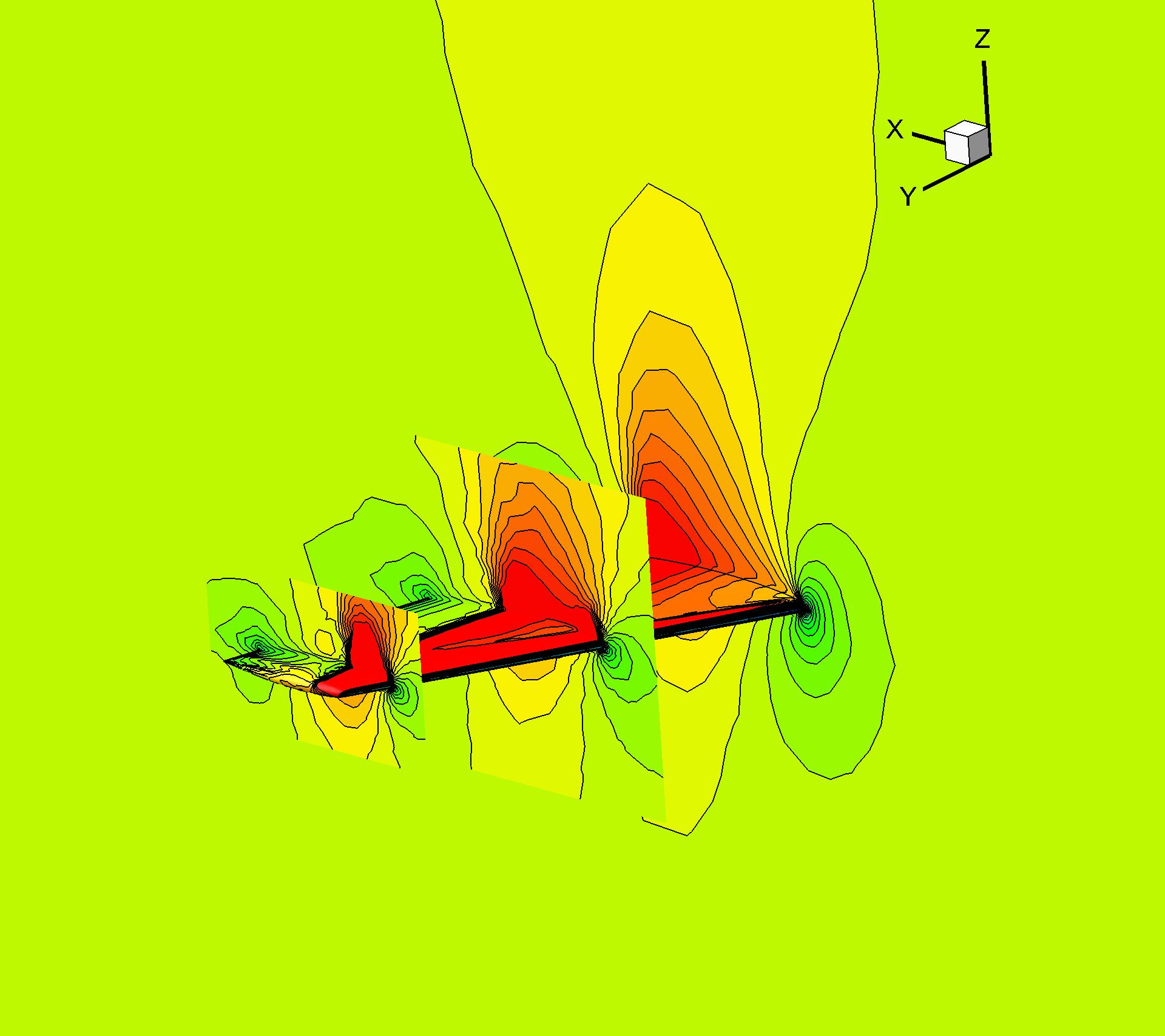}
	 \caption{\label{m6-wing}
		Transonic flow over an inviscid ONERA  M6 wing under Mesh I. Ma=0.8935. AOA=3.06$^{\circ}$. Mesh: 210,663 cells.}
\end{figure}

\begin{figure}[htp]	
	\centering	
	\includegraphics[width=0.32\textwidth]{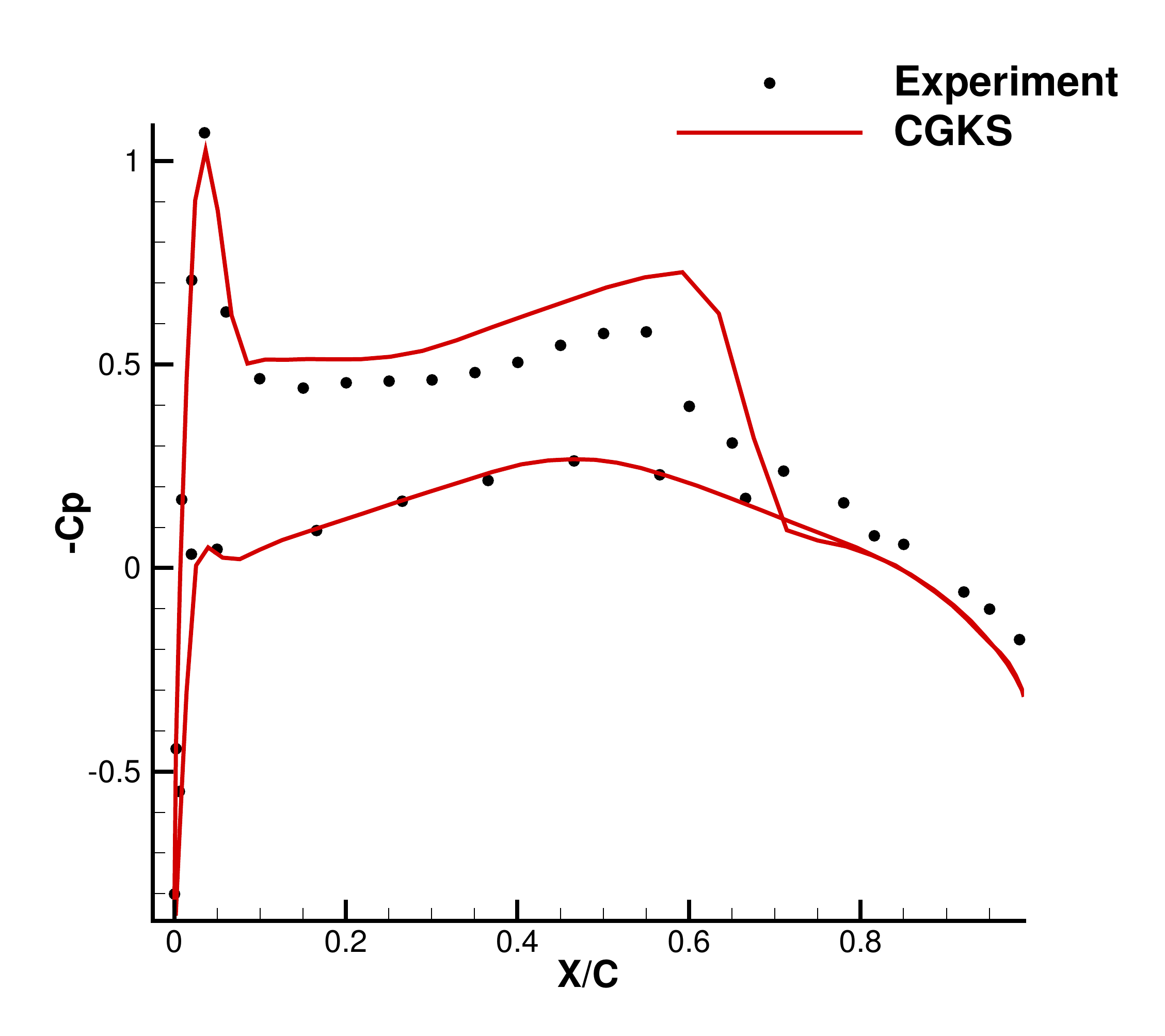}
	\includegraphics[width=0.32\textwidth]{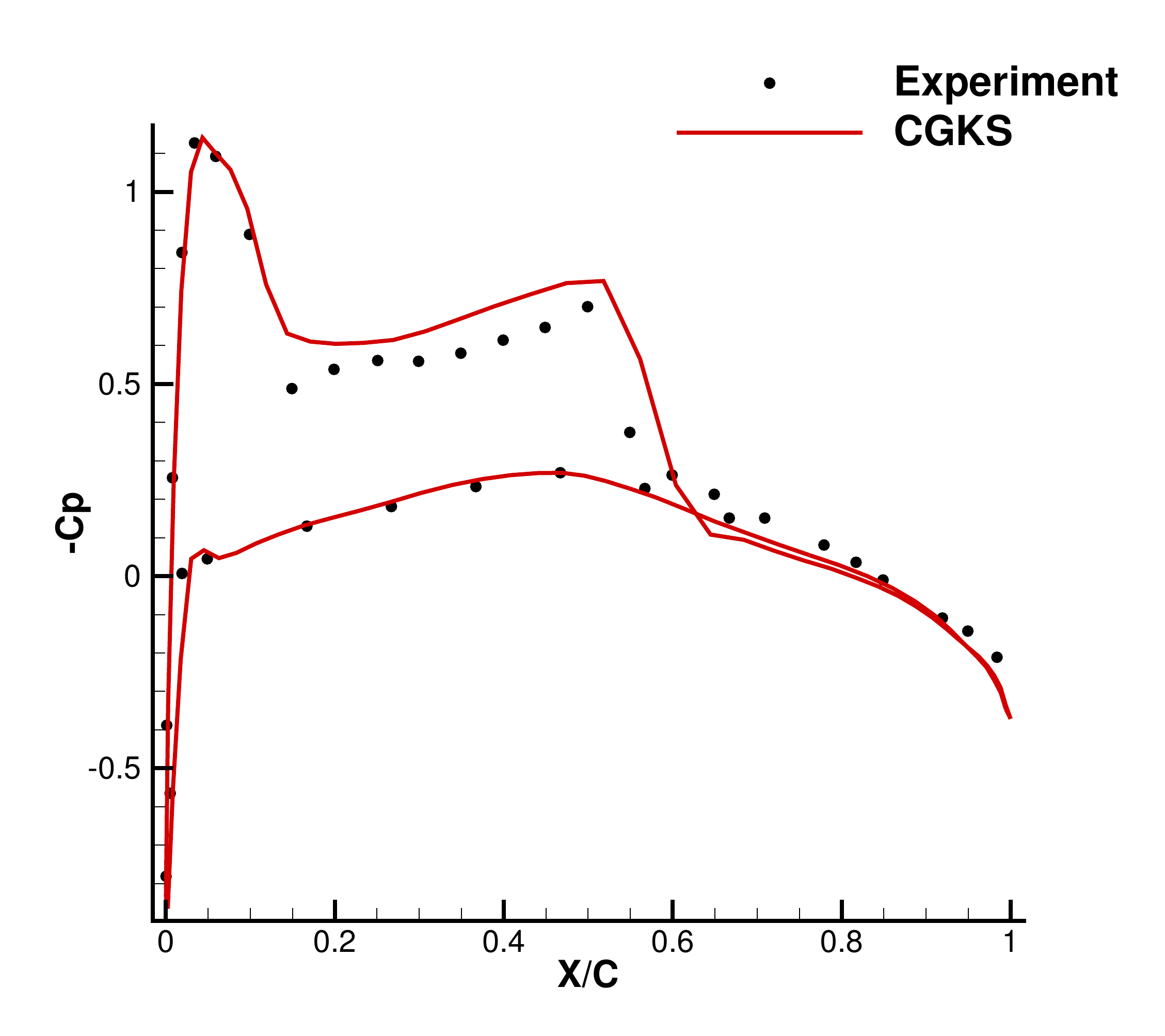}
	\includegraphics[width=0.32\textwidth]{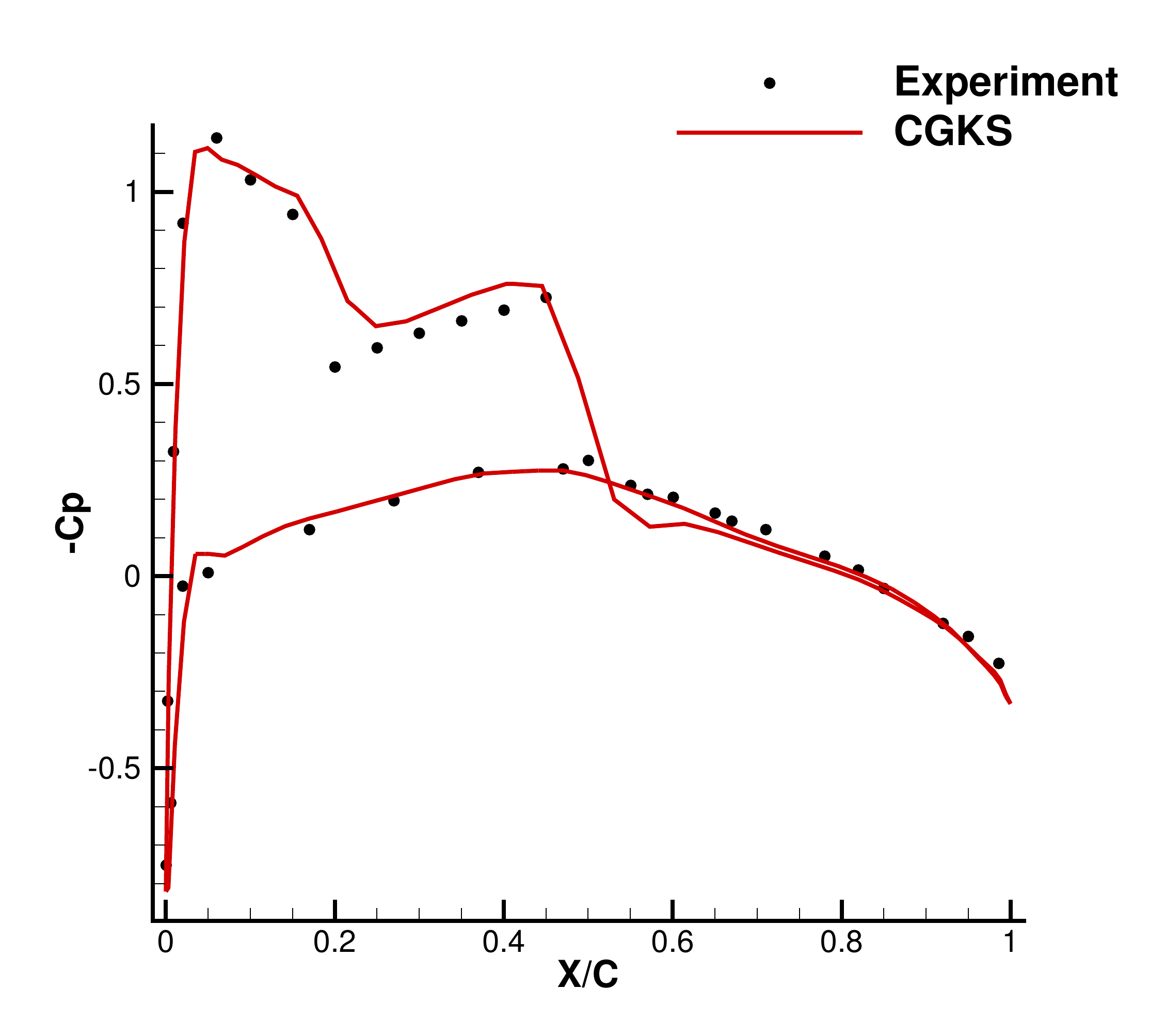}
	\includegraphics[width=0.32\textwidth]{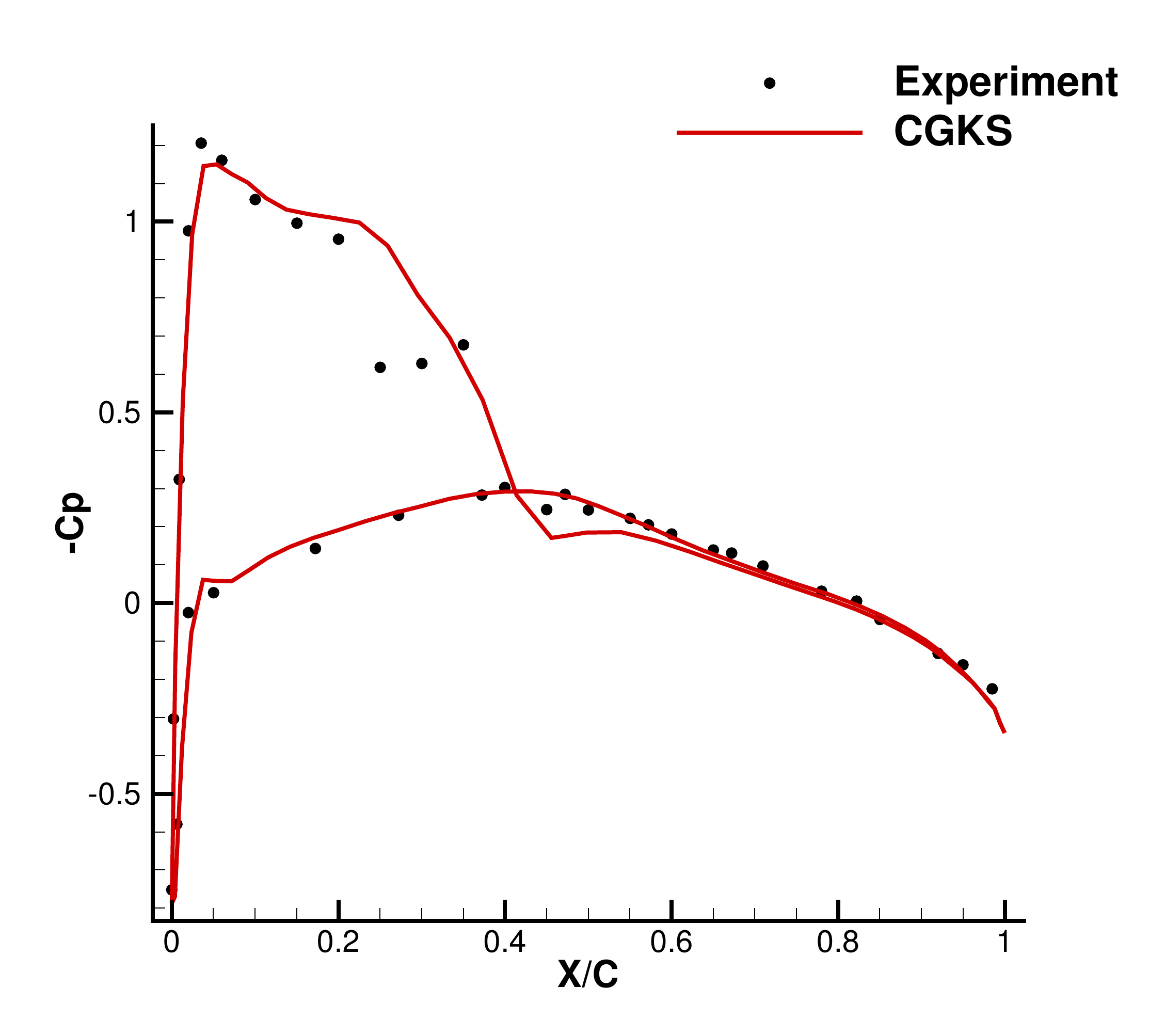}
	\includegraphics[width=0.32\textwidth]{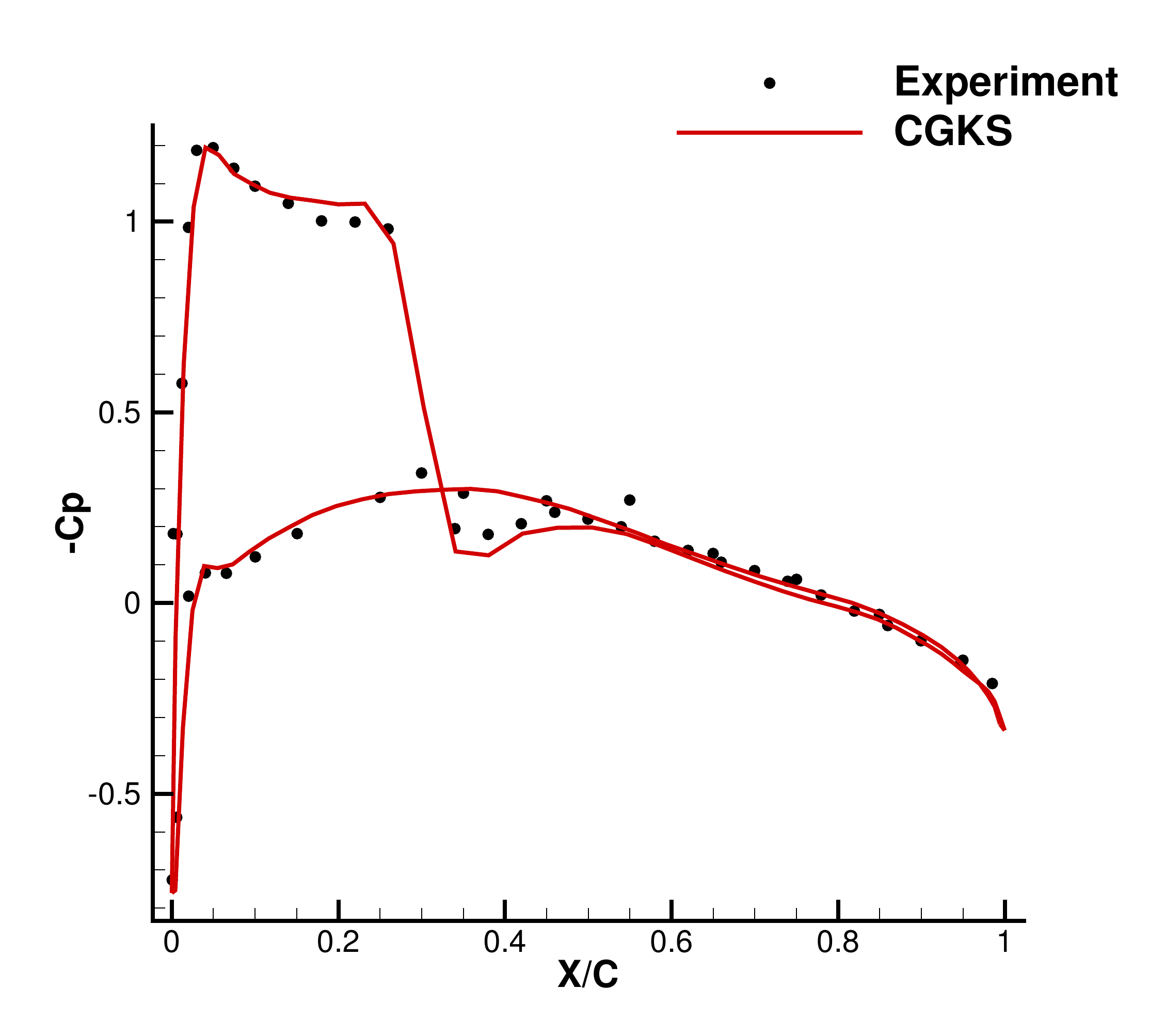}
	\includegraphics[width=0.32\textwidth]{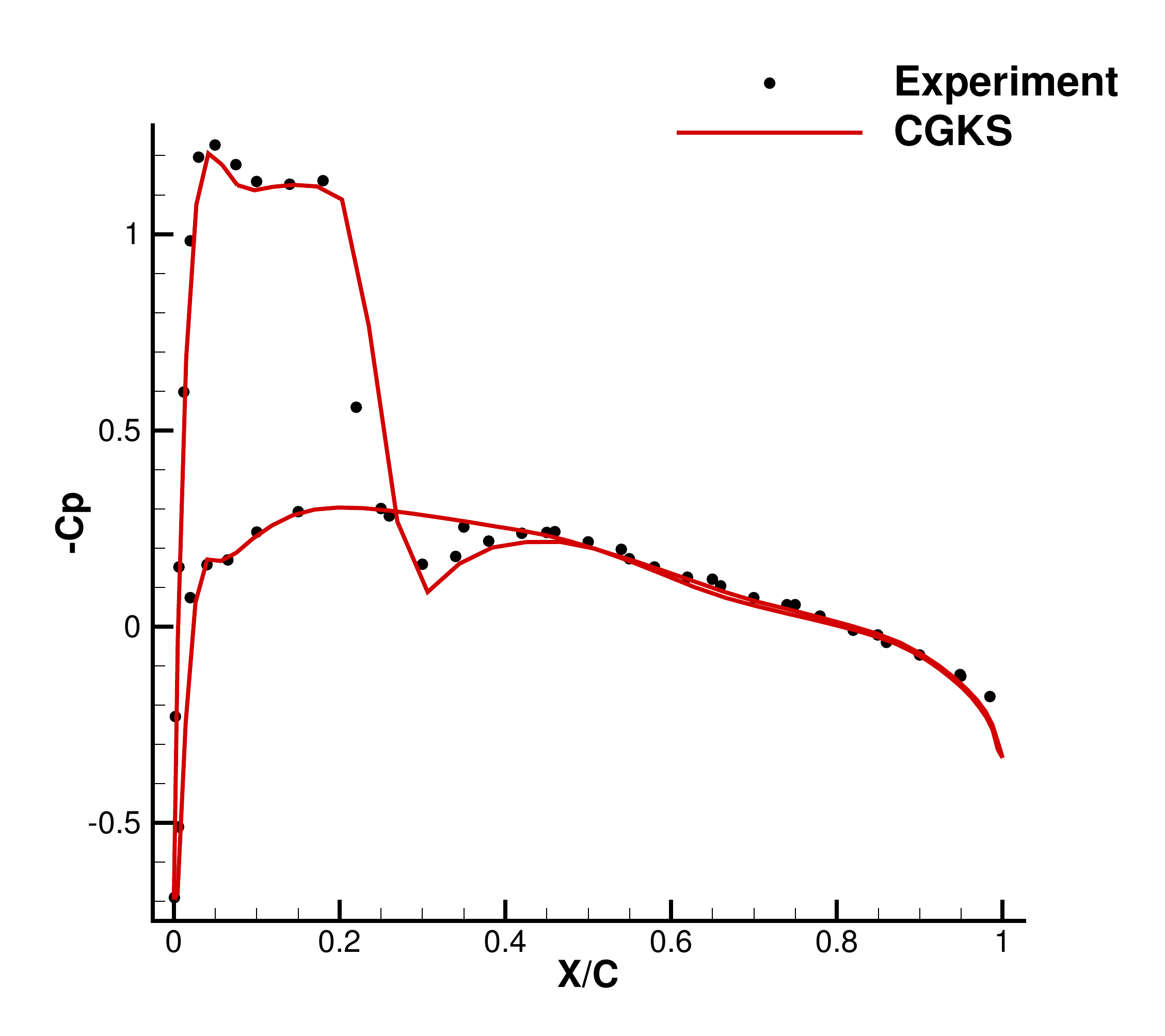}
	\vspace{-4mm} \caption{\label{m6-wing-cp}
		Pressure distributions for wing section at different semi-span locations Y/B on the ONERA  M6 wing. Ma=0.8935. AOA=3.06$^{\circ}$. Top: Y/B=0.20, 0.44, 0.65 from left to right. Bottom: Y/B=0.80, 0.90, 0.95 from left to right.}
\end{figure}

The transonic flow passing through the ONERA M6 wing is tested, as a validation case for compressible external flow.
Experimental data are provided in \cite{schmitt1979pressure} in a high Reynolds number.
Instead of the simulation with a high Reynolds number in the experiment \cite{schmitt1979pressure},
here an inviscid case reported in \cite{liu2020threedimension} is simulated. The incoming Mach number is Ma=0.8395 and the angle of attack is $AOA=3.06^{\circ}$.
The wing is hung on a slipwall, and the Riemann boundaries are applied 10 times of the root chord length away from the wing.
A hybrid unstructured mesh with a near-wall size $h\approx 2e-3$ is used in the computation, as shown in Fig.~\ref{m6-wing}.
The surface pressure distribution and Mach slices at different wing sections under are presented in Fig.~\ref{m6-wing}.
The ``Lambda'' shock is well resolved.
Quantitative comparisons on the pressure distributions for six different locations on the wing are given in Fig.~\ref{m6-wing-cp}, which show the CGKS can give a rough prediction for the real physical flow in this case.

\subsection{Hypersonic space vehicle}

\begin{figure}[htp]	
	\centering	
	\includegraphics[width=0.32\textwidth]
	{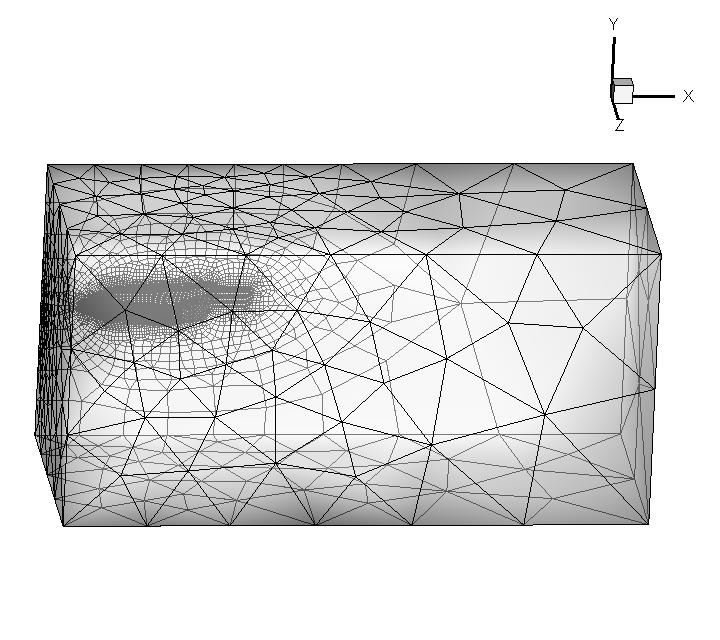}
	\includegraphics[width=0.32\textwidth]
	{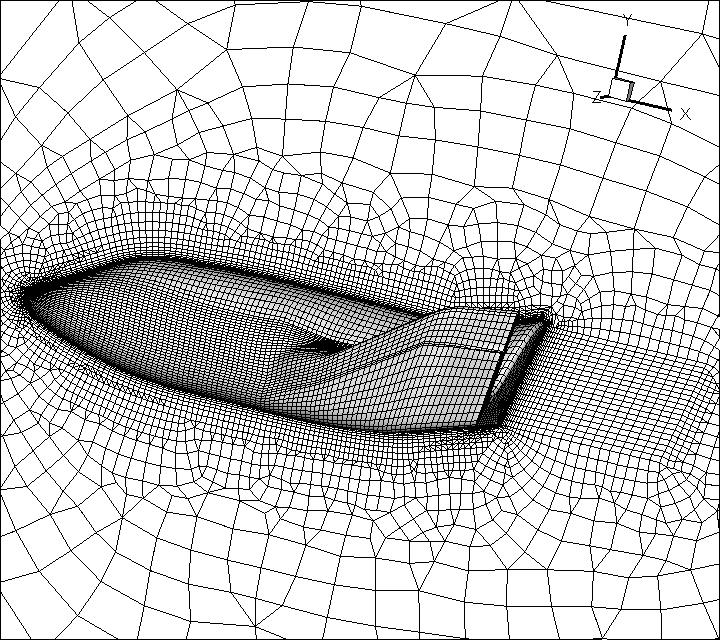}
		\includegraphics[width=0.32\textwidth]
	{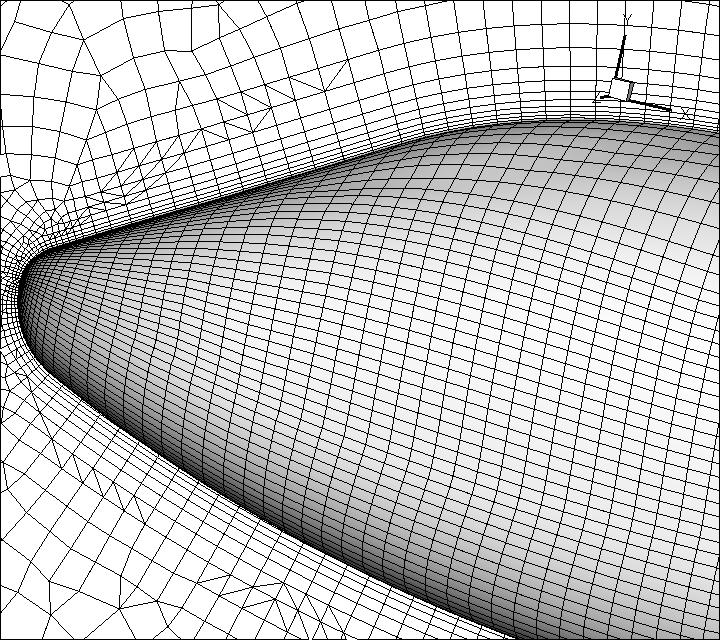}
	 \caption{\label{space-shuttle-mesh}
		Hypersonic flow over a space vehicle. Mesh number: 246,558. }
\end{figure}

\begin{figure}[htp]	
	\centering	
	\includegraphics[height=0.3\textwidth]
	{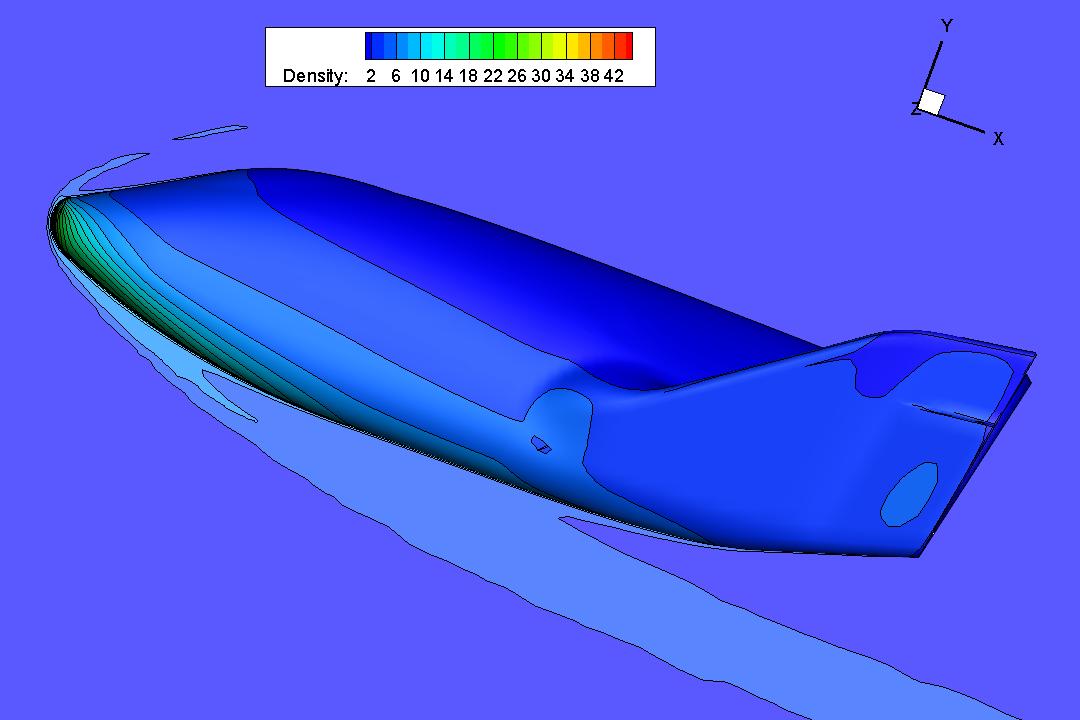}
	\includegraphics[height=0.3\textwidth]
	{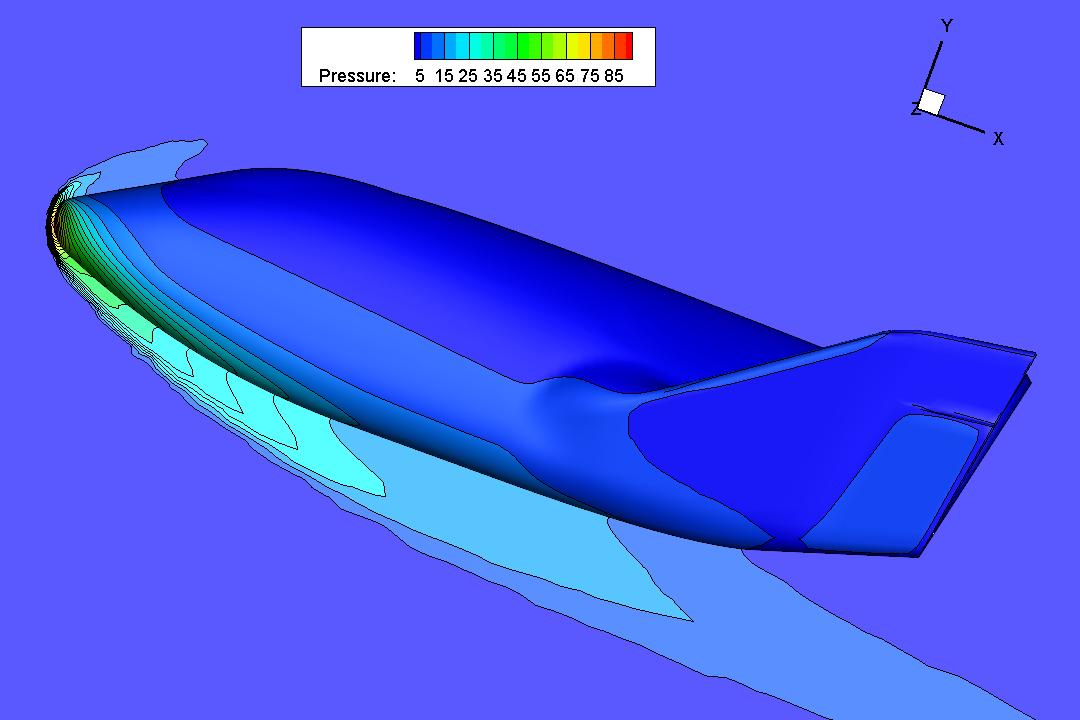}
	 \caption{\label{space-shuttle-pressure}
		Density and pressure distributions for the surface of the space vehicle. }
\end{figure}

\begin{figure}[htp]	
	\centering	
	\includegraphics[width=0.32\textwidth]
	{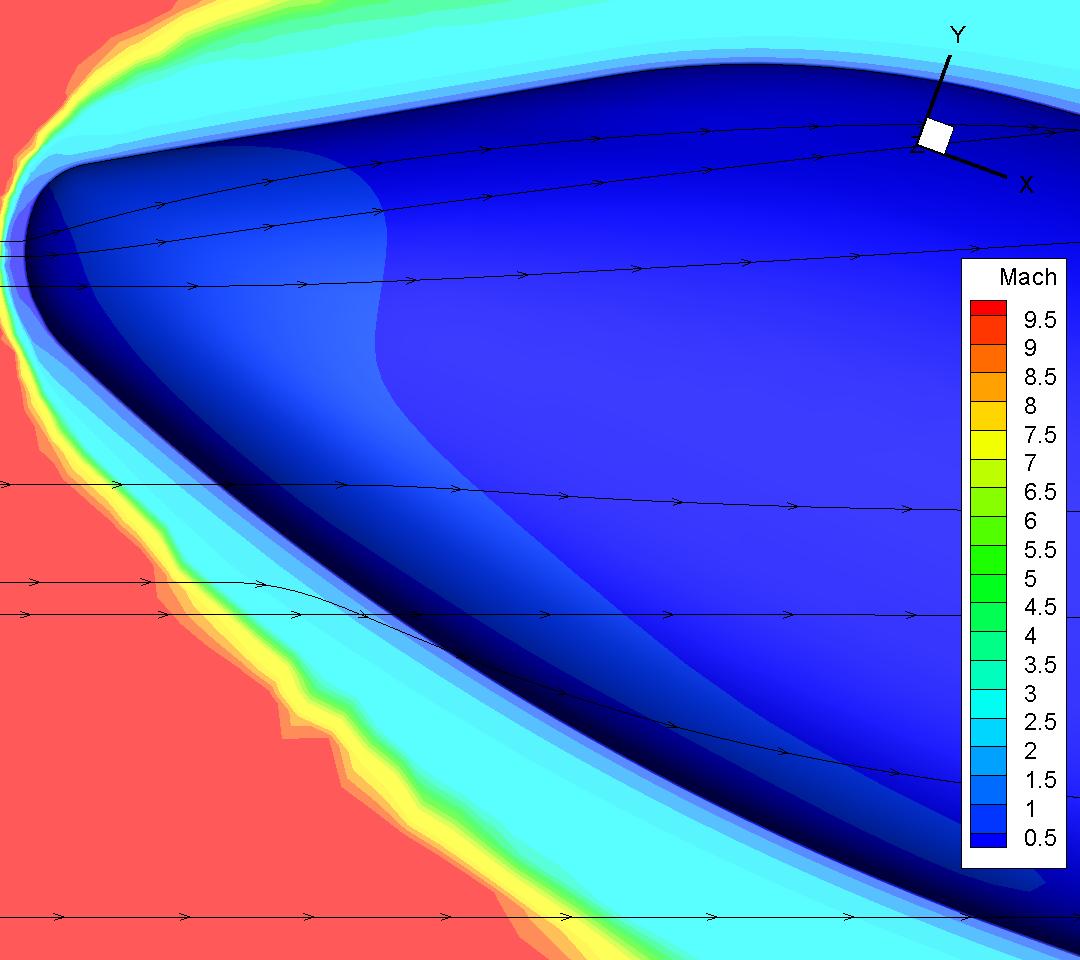}
	\includegraphics[width=0.32\textwidth]
	{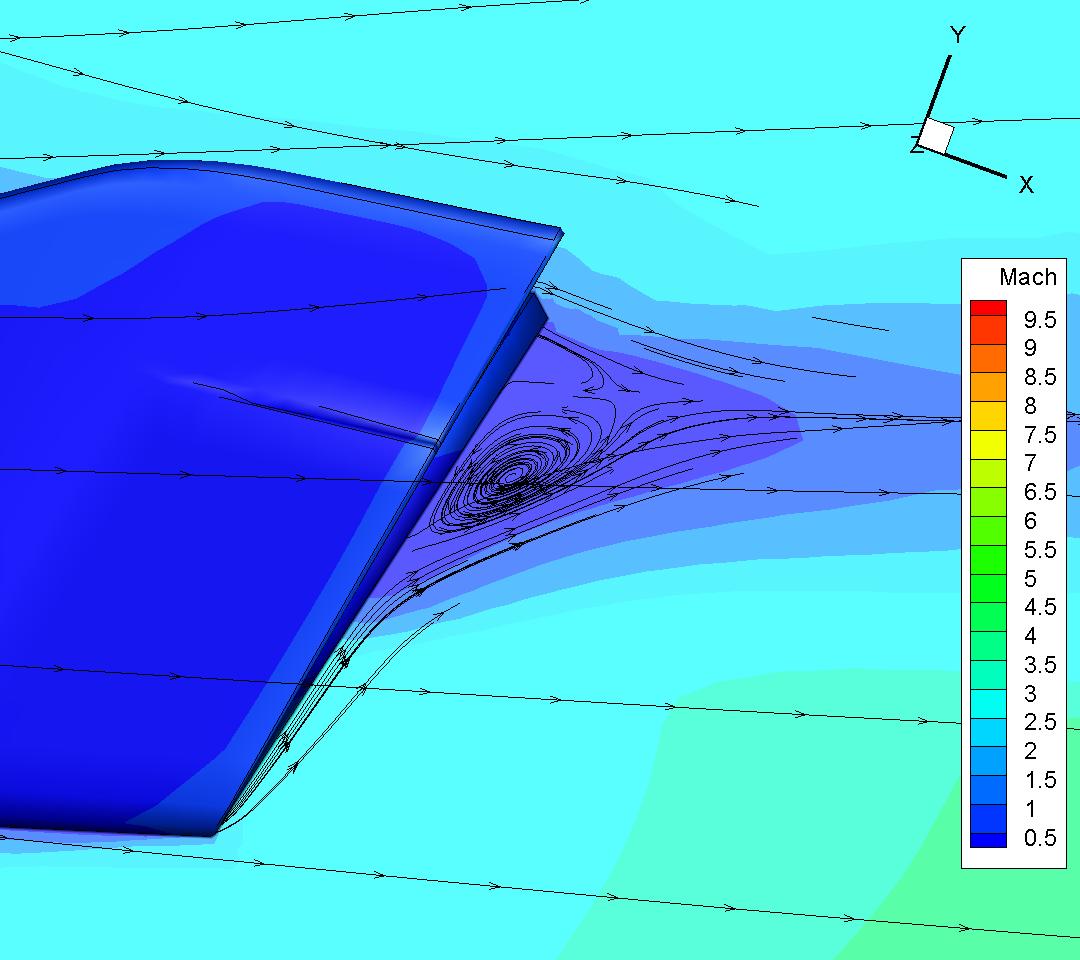}
		\includegraphics[width=0.32\textwidth]
	{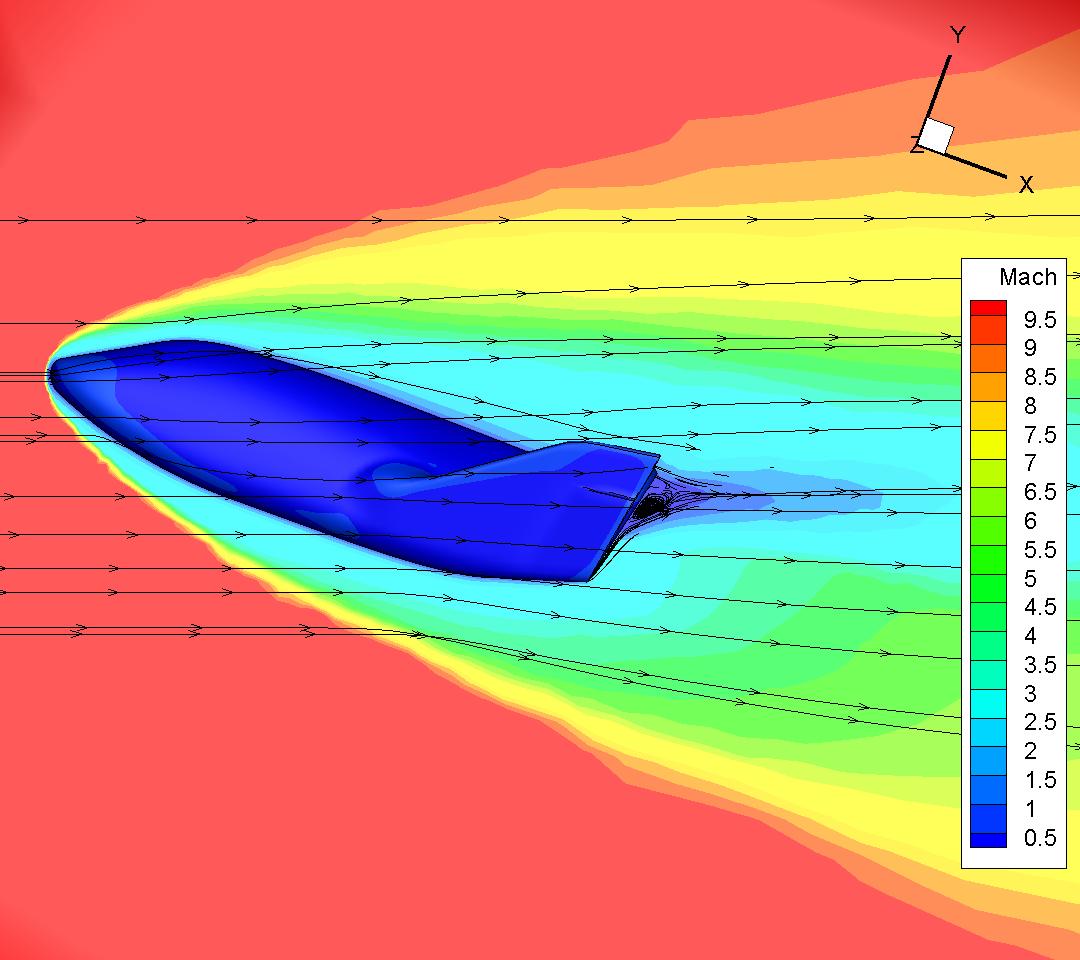}
	\caption{\label{space-shuttle-stream}
		Mach distributions and stream-lines for the hypersonic space vehicle. }
\end{figure}

A space vehicle model is considered to test the robustness of the proposed scheme for the hypersonic viscous flow.
The initial condition is taken as suggested in \cite{chen2020three}, which has  Ma=10 and AOA=20$^{\circ}$.
The Reynolds number is chosen as 14289 based on the vehicle's total length, which yields a Knudsen number Kn=$10^{-3}$. Thus, the flow is in a transition regime and the Maxwell isothermal boundary condition is applied on the vehicle's surface.
The surface mesh is given in Fig.~\ref{space-shuttle-mesh}.
The density and pressure distributions are shown in Fig.~\ref{space-shuttle-pressure}.
The Mach distributions and streamlines are also plotted in Fig.~\ref{space-shuttle-stream}.
The slip velocities can be observed on the surface of the space vehicle.

\section{Discussion and Conclusion}
The third-order CGKS is extended to 3-D mixed-element meshes in this paper.
The expected convergence order is obtained through the accuracy test and it can resolve the unsteady non-linear wave interactions in a high resolution, as shown in the Shu-Osher problem and the supersonic flow over a sphere.
Most importantly, by the first attempt of modifying the updated cell-averaged slopes using the so-called cell-averaged slope compression factor, the present scheme becomes very reliable for the stringent cases with strong shocks and complex geometry. Moreover, the reconstruction method becomes much simpler compared with the traditional WENO reconstruction methods \cite{hu1999weighted}, and the methods used in the previous CGKS \cite{ji2018compact,zhao2019compact,ji2021two}. It not only improves computational efficiency but also code portability.

\begin{figure}[htp]	
	\centering	
	\includegraphics[width=0.5\textwidth]
	{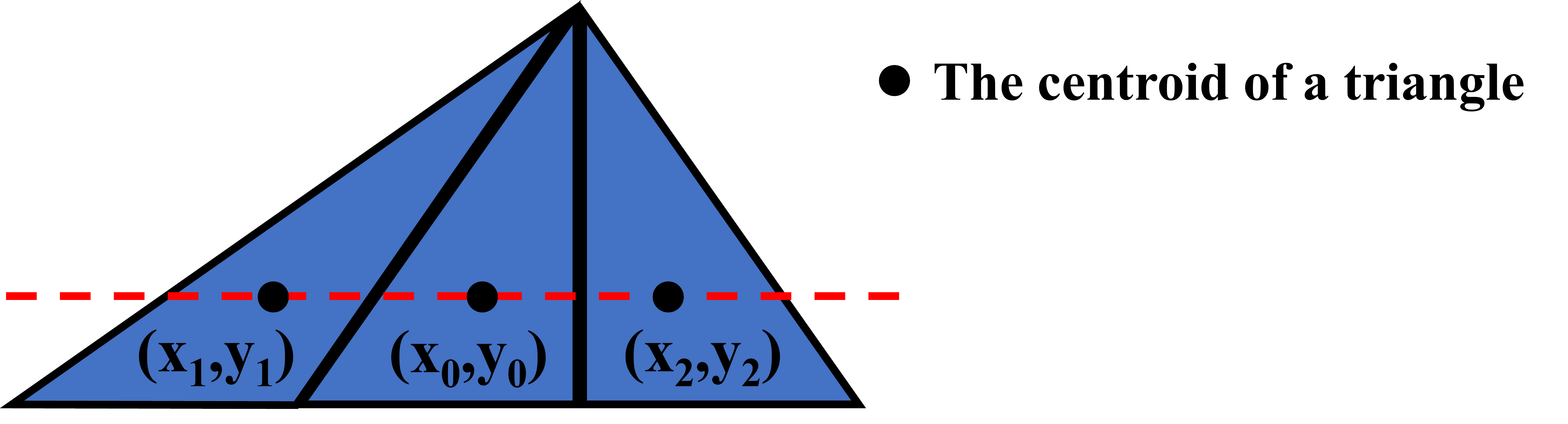}
	\caption{\label{bad-tri}
A singular 2-D case that the first-order polynomial cannot be constructed by the three cell-averaged values in the shown stencil. }
\end{figure}

The motivation of modifying cell-averaged slopes near numerical discontinuities was because CGKS, which worked well on 2D structural as well as unstructured meshes, became much less robust on 3D unstructured meshes, especially on tetrahedral mesh.
The main reason is that the WENO procedure can easily fail to give a reliable reconstructed value at the cell interface in 3-D case.
Although the WENO reconstruction is usually considered as a robust spatial reconstruction, it needs two prerequisites to work normally.
\begin{enumerate}
	\item There is no sub-cell discontinuity.
	\item The mesh quality is good.
\end{enumerate}
Prerequisite 1 cannot be satisfied without introducing the sub-cell resolution \cite{harten1989eno}.
Prerequisite 2 can be easily satisfied for 2-D structured and unstructured meshes. In contrast, the local mesh topology can be much deteriorated for 3-D unstructured mesh generating by modern commercial meshing software.
For example, the tetrahedral meshes near the boundary with a small radius of curvature can have a high volume ratio and the maximum included angle can be very close to 180$^\circ$. A simple 2-D example can be shown in Fig.~\ref{bad-tri} , where the first-order polynomial $p^1=Q_0+a_1 (x-x_0) +a_2 (y-y_0)$ cannot be determined since the y coordinates of the centroids $y_0=y_1=y_2$. Then if it is one of the sub-stencil in the traditional WENO reconstruction, the smoothness for such a sub-stencil cannot be correctly measured by the definition of the smooth indicator in Eq.~\eqref{smooth-indicator}.
Following this limiting case, suppose a situation where a targeted cell $\Omega_{0}$ has two sub-stencils, one is in smooth flow region but with bad mesh quality and one has a good mesh topology but crossing a weak shock, then the smooth indicators given by these two sub-stencils could be about the same and both incorrectly evaluated.
So the final WENO reconstruction gives an invalid slope, leading to a reduction in robustness.

Surely parameters related to the geometric quality can be chosen as weight functions in reconstruction to fix the geometric singularity to some degree.
The parameters can be volume ratio, mesh skewness factor, etc.
However, the selection of the parameters and proper weighting forms requires a lot of experience and perhaps machine learning.
Also, it will make the high-order method more complicated.

On the other hand, the first-order finite volume schemes with approximate Riemann solvers are positive-preserving under CFL condition, e.g., the first-order L-F scheme, and the first-order kinetic vector flux-splitting scheme \cite{tang1999gas}. The first-order GKS, even though not strictly  positive-preserving, but is robust enough for the hypersonic flow.
Since the reconstruction stencils in these schemes are only the targeted cell itself, they have great mesh adaptability.

\begin{figure}[htp]	
 	\centering	
 	\includegraphics[width=0.5\textwidth]
 	{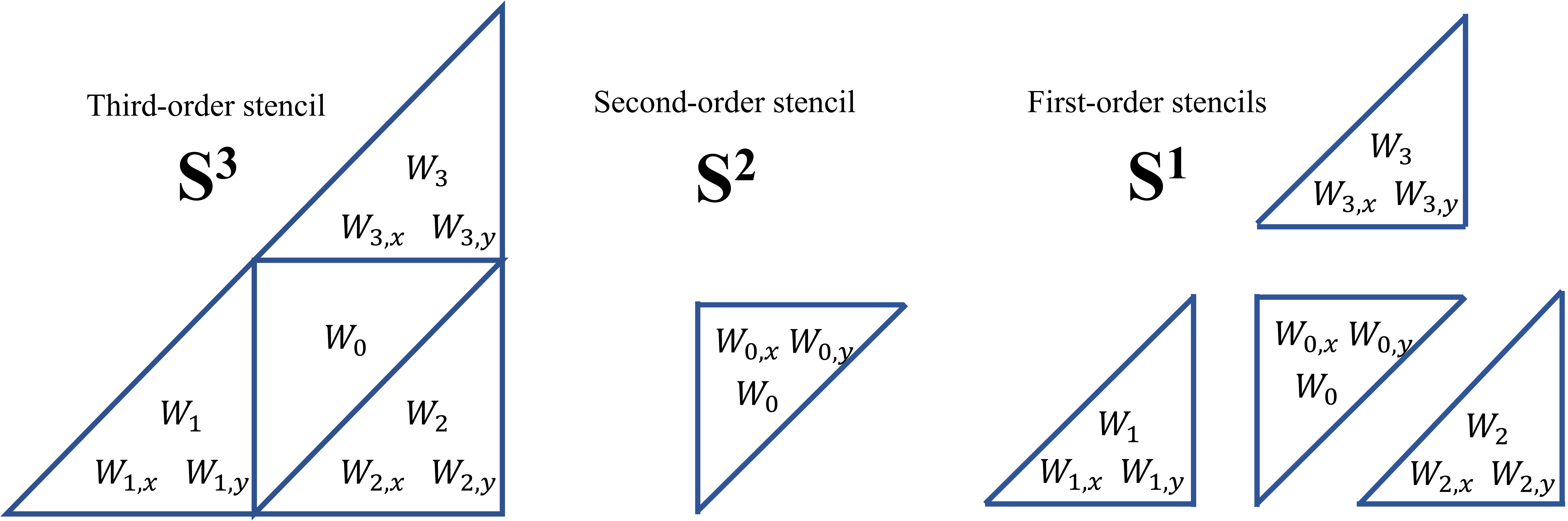}
 	\caption{\label{stencil}
 		An 2-D example for stencil selections for the compact multi-resolution WENO reconstruction. }
 \end{figure}

\begin{figure}[htp]	
	\centering	
	\includegraphics[width=0.8\textwidth]
	{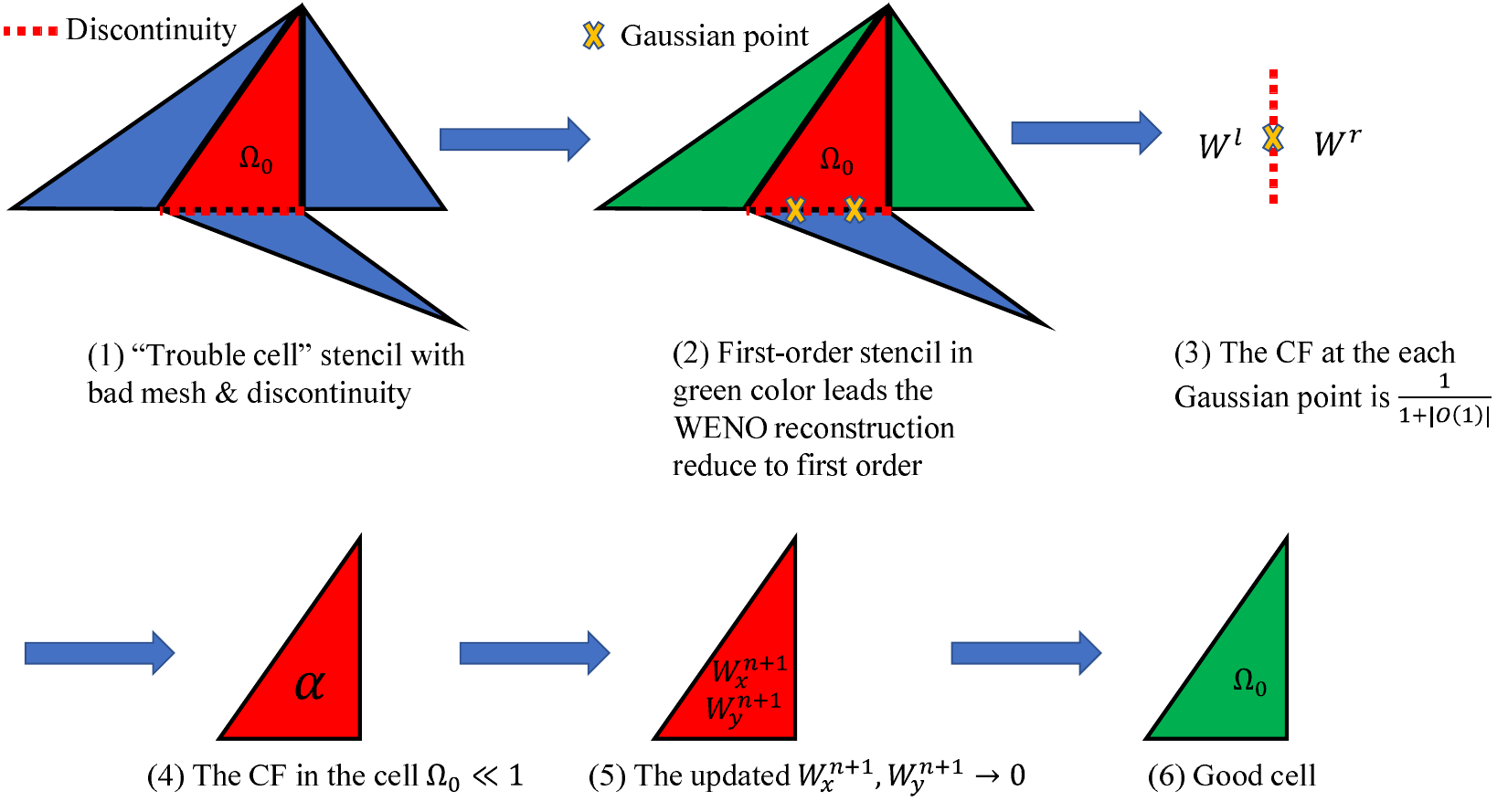}
	\caption{\label{why-cf-work}
		A 2-D example that how the compact GKS with the CF handles the bad mesh topology and discontinuity. }
\end{figure}

Thus, it inspires us to design a new reconstruction strategy, in which the high-order CGKS can truly reduce to first-order GKS when necessary.
In turn, three properties should be satisfied based on the WENO-type reconstruction and the CGKS framework,
\begin{enumerate}
 	\item The first-order sub-stencil (i.e., the zeroth polynomial  $p^0$  determined solely by the cell-averaged conservative variable on the targeted cell) must be included
 	\item The smooth indicator of the above sub-stencil should be as independent with grid quality as possible.
 	\item The smooth indicator of the above sub-stencil should be small enough if a discontinuity is located inside the targeted cell.
\end{enumerate}
The multi-resolution WENO \cite{zhu2018multiresolution} satisfies Property 1, which has also been adopted in the previous CGKS \cite{ji2021two}.
Since both cell-averaged conservative variables and their slopes are updated within a cell for CGKS, a first-order polynomial inside each cell can be always constructed, i.e., $p^1_i=\overline{W}_i +\overline{W}_{x,i} x_i$.
As reported in DG methods \cite{yang2019robust}, such type of polynomial has great mesh adaptability.
As a result, all these polynomials in the compact stencil can be used to evaluate the smooth indicator of $p^0$, which satisfies Property 2.
An illustration for stencil selections are shown in Fig.~\ref{stencil}.
Property 3 requires the sub-cell resolution of a scheme.
The numerical discontinuities through reconstruction around the targeted cell at $t^n$ step, either caused by physical discontinuities or mesh irregularity, will lead to a possible sub-cell discontinuity at $t^{n+1}$, and the divergence theorem in Eq.~\eqref{gauss-formula} for the cell-averaged gradients will be invalid.
Thus, it is reasonable to introduce the CF $\alpha$ proportional to the strength of the discontinuities, as shown in Eq.~\eqref{cf-face}, and compress the $\overline{W}_{x,i}^{n+1}$ accordingly. In this way, the smooth indicator will be smaller if the possible sub-cell discontinuity is stronger. The CF can be seen as an indicator of measuring the strength of the sub-cell discontinuity.
Finally, a delightful feedback mechanism is established, as shown in Fig.~\ref{why-cf-work}, in case of solution irregularity.
Different from most of limiters \cite{krivodonova2004shock,zhang2017positivity,clain2011high}, which obey the principle that finding the ``trouble cell'' first and removing it next, the current CGKS takes one more step that turning the ``trouble cell'' to become a ``good cell''.
In summary, the current CGKS, with the large explicit time step, great robustness and mesh adaptability, and high program portability, is confident for DNS or LES simulation for highly compressible flow with complex geometry.
The remaining problems include the construction of implicit CGKS with acceleration techniques and the extension to RANS simulation.

\section*{Acknowledgments}
The authors would like to thank Mr. Yipei Chen for the help on the setting of the hypersonic space vehicle case.
The current research is supported by National Numerical Windtunnel project  and  National Science Foundation of China (11772281, 91852114).

\section*{References}
\bibliographystyle{plain}%
\bibliography{jixingbib}

\begin{thebibliography}{10}

\bibitem{antoniadis2017assessment}
Antonis~F. Antoniadis, Panagiotis Tsoutsanis, and Dimitris Drikakis.
\newblock Assessment of high-order finite volume methods on unstructured meshes
  for {RANS} solutions of aeronautical configurations.
\newblock {\em Computers \& Fluids}, 146:86 -- 104, 2017.

\bibitem{bharti2006steady}
Ram~Prakash Bharti, R.~P. Chhabra, and V.~Eswaran.
\newblock Steady flow of power law fluids across a circular cylinder.
\newblock {\em The Canadian Journal of Chemical Engineering}, 84(4):406--421,
  2006.

\bibitem{BGK}
Prabhu~Lal Bhatnagar, Eugene~P Gross, and Max Krook.
\newblock A model for collision processes in gases. {I}. {Small} amplitude
  processes in charged and neutral one-component systems.
\newblock {\em Physical Review}, 94(3):511, 1954.

\bibitem{braza1986Numerical}
M.~Braza, P.~Chassaing, and H.~Ha Minh.
\newblock Numerical study and physical analysis of the pressure and velocity
  fields in the near wake of a circular cylinder.
\newblock {\em Journal of Fluid Mechanics}, 165:79–130, 1986.

\bibitem{chen2020three}
Yipei Chen, Yajun Zhu, and Kun Xu.
\newblock A three-dimensional unified gas-kinetic wave-particle solver for flow
  computation in all regimes.
\newblock {\em Physics of Fluids}, 32(9):096108, 2020.

\bibitem{cheng2017parallel}
Jian Cheng, Xiaodong Liu, Tiegang Liu, and Hong Luo.
\newblock A parallel, high-order direct discontinuous {Galerkin} method for the
  {Navier-Stokes} equations on {3D} hybrid grids.
\newblock {\em Communications in Computational Physics}, 21(5):1231–1257,
  2017.

\bibitem{clain2011high}
S.~Clain, S.~Diot, and R.~Loubère.
\newblock A high-order finite volume method for systems of conservation
  law--multi-dimensional optimal order detection ({MOOD}).
\newblock {\em Journal of Computational Physics}, 230(10):4028--4050, 2011.

\bibitem{coutanceau1977experimental}
Madeleine Coutanceau and Roger Bouard.
\newblock Experimental determination of the main features of the viscous flow
  in the wake of a circular cylinder in uniform translation. {Part 1. Steady
  flow}.
\newblock {\em Journal of Fluid Mechanics}, 79(2):231--256, 1977.

\bibitem{dumbser2010arbitrary}
Michael Dumbser.
\newblock {Arbitrary high order PNPM schemes on unstructured meshes for the
  compressible Navier--Stokes equations}.
\newblock {\em Computers \& Fluids}, 39(1):60--76, 2010.

\bibitem{grove1964experimental}
A.~S. Grove, F.~H. Shair, and E.~E. Petersen.
\newblock An experimental investigation of the steady separated flow past a
  circular cylinder.
\newblock {\em Journal of Fluid Mechanics}, 19(1):60–80, 1964.

\bibitem{harten1989eno}
Ami Harten.
\newblock {ENO} schemes with subcell resolution.
\newblock {\em Journal of Computational Physics}, 83(1):148--184, 1989.

\bibitem{hu1999weighted}
Changqing Hu and Chi-Wang Shu.
\newblock Weighted essentially non-oscillatory schemes on triangular meshes.
\newblock {\em Journal of Computational Physics}, 150(1):97--127, 1999.

\bibitem{Huynh2007FR}
Hung~T Huynh.
\newblock A flux reconstruction approach to high-order schemes including
  discontinuous {Galerkin} methods.
\newblock In {\em 18th AIAA Computational Fluid Dynamics Conference}, page
  4079, 2007.

\bibitem{ji2019high}
Xing Ji.
\newblock {\em High-order non-compact and compact gas-kinetic schemes}.
\newblock PhD thesis, Hong Kong Univeristy of Science and Technology, 2019.

\bibitem{ji2018compact}
Xing Ji, Liang Pan, Wei Shyy, and Kun Xu.
\newblock A compact fourth-order gas-kinetic scheme for the {Euler} and
  {Navier-Stokes} equations.
\newblock {\em Journal of Computational Physics}, 372:446 -- 472, 2018.

\bibitem{ji2021compact}
Xing Ji, Fengxiang Zhao, Wei Shyy, and Kun Xu.
\newblock Compact high-order gas-kinetic scheme for three-dimensional flow
  simulations.
\newblock {\em AIAA Journal}, 0(0):1--18, 0.

\bibitem{ji2021two}
Xing Ji, Fengxiang Zhao, Wei Shyy, and Kun Xu.
\newblock Two-step multi-resolution reconstruction-based compact gas-kinetic
  scheme on tetrahedral mesh.
\newblock {\em arXiv preprint arXiv:2102.01366}, 2021.

\bibitem{krivodonova2004shock}
L.~Krivodonova, J.~Xin, J.-F. Remacle, N.~Chevaugeon, and J.E. Flaherty.
\newblock Shock detection and limiting with discontinuous {Galerkin} methods
  for hyperbolic conservation laws.
\newblock {\em Applied Numerical Mathematics}, 48(3):323--338, 2004.
\newblock Workshop on Innovative Time Integrators for PDEs.

\bibitem{li2016twostage}
Jiequan Li and Zhifang Du.
\newblock A two-stage fourth order time-accurate discretization for
  {Lax--Wendroff} type flow solvers {I.} hyperbolic conservation laws.
\newblock {\em SIAM Journal on Scientific Computing}, 38(5):A3046--A3069, 2016.

\bibitem{li2014efficient}
Wanai Li.
\newblock {\em Efficient implementation of high-order accurate numerical
  methods on unstructured grids}.
\newblock Berlin, Heidelberg: Springer, 2014.

\bibitem{liu2020threedimension}
Yangyang Liu, Liming Yang, Chang Shu, and Huangwei Zhang.
\newblock Three-dimensional high-order least square-based finite
  difference-finite volume method on unstructured grids.
\newblock {\em Physics of Fluids}, 32(12):123604, 2020.

\bibitem{luo2010reconstructed}
Hong Luo, Luqing Luo, Robert Nourgaliev, Vincent~A Mousseau, and Nam Dinh.
\newblock {A reconstructed discontinuous Galerkin method for the compressible
  Navier--Stokes equations on arbitrary grids}.
\newblock {\em Journal of Computational Physics}, 229(19):6961--6978, 2010.

\bibitem{Nagata2016sphere}
T.~Nagata, T.~Nonomura, S.~Takahashi, Y.~Mizuno, and K.~Fukuda.
\newblock Investigation on subsonic to supersonic flow around a sphere at low
  {Reynolds} number of between 50 and 300 by direct numerical simulation.
\newblock {\em Physics of Fluids}, 28(5):056101, 2016.

\bibitem{pan2016unstructuredcompact}
Liang Pan and Kun Xu.
\newblock A third-order compact gas-kinetic scheme on unstructured meshes for
  compressible {Navier--Stokes} solutions.
\newblock {\em Journal of Computational Physics}, 318:327--348, 2016.

\bibitem{pan2020high}
Liang Pan and Kun Xu.
\newblock High-order gas-kinetic scheme with three-dimensional {WENO}
  reconstruction for the {Euler} and {Navier-Stokes} solutions.
\newblock {\em Computers \& Fluids}, 198:104401, 2020.

\bibitem{Pan2016twostage}
Liang Pan, Kun Xu, Qibing Li, and Jiequan Li.
\newblock An efficient and accurate two-stage fourth-order gas-kinetic scheme
  for the {Euler} and {Navier--Stokes} equations.
\newblock {\em Journal of Computational Physics}, 326:197--221, 2016.

\bibitem{schmitt1979pressure}
V~Schmitt.
\newblock Pressure distributions on the {ONERA} {M6-wing} at transonic mach
  numbers, experimental data base for computer program assessment.
\newblock {\em AGARD AR-138}, 1979.

\bibitem{multi-derivative}
David~C Seal, Yaman G{\"u}{\c{c}}l{\"u}, and Andrew~J Christlieb.
\newblock High-order multiderivative time integrators for hyperbolic
  conservation laws.
\newblock {\em Journal of Scientific Computing}, 60(1):101--140, 2014.

\bibitem{shu2016weno-dg-review}
Chi-Wang Shu.
\newblock High order {WENO} and {DG} methods for time-dependent
  convection-dominated {PDEs}: A brief survey of several recent developments.
\newblock {\em Journal of Computational Physics}, 316:598 -- 613, 2016.

\bibitem{shu1989efficient}
Chi-Wang Shu and Stanley Osher.
\newblock {Efficient implementation of essentially non-oscillatory
  shock-capturing schemes, II}.
\newblock In {\em Upwind and High-Resolution Schemes}, pages 328--374.
  Springer, 1989.

\bibitem{taneda1956experimental}
Sadatoshi Taneda.
\newblock Experimental investigation of the wakes behind cylinders and plates
  at low {Reynolds} numbers.
\newblock {\em Journal of the Physical Society of Japan}, 11(3):302--307, 1956.

\bibitem{tang1999gas}
Tao Tang and Kun Xu.
\newblock Gas-kinetic schemes for the compressible {Euler} equations:
  {Positivity}-preserving analysis.
\newblock {\em Zeitschrift für angewandte Mathematik und Physik ZAMP},
  50(2):258--281, 1999.

\bibitem{tritton1959experiments}
David~J Tritton.
\newblock Experiments on the flow past a circular cylinder at low {Reynolds}
  numbers.
\newblock {\em Journal of Fluid Mechanics}, 6(4):547--567, 1959.

\bibitem{wang2017thesis}
Qian Wang.
\newblock {\em Compact High-Order Finite Volume Method on Unstructured Grids}.
\newblock PhD thesis, Tsinghua University, 6 2017.

\bibitem{wang2017towards}
ZJ~Wang, Y~Li, F~Jia, GM~Laskowski, J~Kopriva, U~Paliath, and R~Bhaskaran.
\newblock Towards industrial large eddy simulation using the {FR/CPR} method.
\newblock {\em Computers \& Fluids}, 156:579--589, 2017.

\bibitem{GKS-2001}
Kun Xu.
\newblock A gas-kinetic {BGK} scheme for the {Navier--Stokes} equations and its
  connection with artificial dissipation and {Godunov} method.
\newblock {\em Journal of Computational Physics}, 171(1):289--335, 2001.

\bibitem{xu2014direct}
Kun Xu.
\newblock {\em Direct Modeling for Computational Fluid Dynamics: Construction
  and Application of Unified Gas-Kinetic Schemes}.
\newblock World Scientific, 2014.

\bibitem{yang2019robust}
Xiaoquan Yang, Jian Cheng, Hong Luo, and Qijun Zhao.
\newblock Robust implicit direct discontinuous {Galerkin} method for simulating
  the compressible turbulent flows.
\newblock {\em AIAA Journal}, 57(3):1113--1132, 2019.

\bibitem{yu2014accuracy}
Meilin Yu, Z.J. Wang, and Yen Liu.
\newblock On the accuracy and efficiency of discontinuous {Galerkin}, spectral
  difference and correction procedure via reconstruction methods.
\newblock {\em Journal of Computational Physics}, 259:70 -- 95, 2014.

\bibitem{zhang2019direct}
Fan Zhang, Jian Cheng, and Tiegang Liu.
\newblock A direct discontinuous {Galerkin method for the incompressible
  Navier--Stokes} equations on arbitrary grids.
\newblock {\em Journal of Computational Physics}, 380:269--294, 2019.

\bibitem{zhang2017positivity}
Xiangxiong Zhang.
\newblock On positivity-preserving high order discontinuous {Galerkin} schemes
  for compressible {Navier-Stokes} equations.
\newblock {\em Journal of Computational Physics}, 328:301--343, 2017.

\bibitem{zhao2019compact}
Fengxiang Zhao, Xing Ji, Wei Shyy, and Kun Xu.
\newblock Compact higher-order gas-kinetic schemes with spectral-like
  resolution for compressible flow simulations.
\newblock {\em Advances in Aerodynamics}, 1(1):13, 2019.

\bibitem{zhao2020compact}
Fengxiang Zhao, Xing Ji, Wei Shyy, and Kun Xu.
\newblock A compact high-order gas-kinetic scheme on unstructured mesh for
  acoustic and shock wave computations.
\newblock {\em arXiv preprint arXiv:2010.05717}, 2020.

\bibitem{zhu2018multiresolution}
Jun Zhu and Chi-Wang Shu.
\newblock A new type of multi-resolution {WENO} schemes with increasingly
  higher order of accuracy.
\newblock {\em Journal of Computational Physics}, 375:659--683, 2018.

\bibitem{zhu2020new}
Jun Zhu and Chi-Wang Shu.
\newblock A new type of third-order finite volume multi-resolution {WENO}
  schemes on tetrahedral meshes.
\newblock {\em Journal of Computational Physics}, 406:109212, 2020.

\end{thebibliography}

\end{document}